\documentclass[11pt]{article}
\usepackage[english]{babel}
\usepackage{latexsym}
\usepackage{epsfig}
\usepackage{graphicx}
\usepackage{comment}
\usepackage{amsfonts,amsmath}
\usepackage{latexsym}
\usepackage{color}
\newsymbol\blacktriangleright 1349
\newsymbol\blacktriangleleft 134A
\newsymbol\vartriangle 134D

\newtheorem{thm}{Theorem}

\newtheorem{corol}[thm]{Corollary}
\newtheorem{defin}{Definition}
\newtheorem{lemma}[thm]{Lemma}

\newtheorem{exa}{Example}
\newenvironment{example}{\begin{exa} \rm }{\hfill$\Box$ \end{exa}}

\newenvironment{proof}{\begin{trivlist}\item[] \mbox{\it Proof. }}
{\hfill$\Box$ \end{trivlist}}

\definecolor{lgrey}{gray}{.55}%
\definecolor{hgrey}{gray}{.3}%

\def\dv{\mathbb}
\def\caP{{\cal P}(N)}
\def\chain{{\cal C}}
\def\cone{\mbox{\rm cone}\,}
\def\conv{\mbox{\rm conv}\,}
\def\cor{C} 
\def\cors{\mbox{\tiny $C$}} 
\def\Cov{\mbox{\rm Cov}\,} 
\def\cover{\prec\!\!\cdot\,\,}
\def\cl{\mbox{\rm cl}\,}
\def\dist{\mbox{\rm dist}\,}
\def\diag{\Delta} 
\def\tr{\mbox{\rm tr}} 
\def\ne{\mbox{\rm ne}} 
\def\elemN{\mathcal{E}(N)}
\def\permA{\varepsilon} 
\def\permB{\eta} 
\def\permC{\sigma} 
\def\permD{\pi} 
\def\permE{\zeta} 
\def\permF{\xi} 
\def\permG{\varsigma} 
\def\perN{\Pi(N)}
\def\perNs{\mbox{\tiny$\Pi(N)$}}
\def\preced{\circ} 
\def\braid{\bullet} 
\def\topol{\textcolor{lgrey}{\bullet}} 
\def\outside{\star} 
\def\norfan{\diamond} 
\def\lGalo{\triangleleft} 
\def\rGalo{\triangleright} 
\def\lGprec{\lhd} 
\def\rGprec{\rhd} 
\def\lGbraid{\blacktriangleleft} 
\def\rGbraid{\blacktriangleright} 
\def\lGtopol{\textcolor{lgrey}{\blacktriangleleft}} 
\def\rGtopol{\textcolor{lgrey}{\blacktriangleright}} 
\def\rankv{\rho} 
\def\trans{\tau} 
\def\diff{\vartriangle\!\!} 
\def\inte{\mbox{\rm int}\,} 
\def\relint{\mbox{\rm relint}\,} 
\def\face{\mathsf{F}} 
\def\enuN{\Upsilon (N)} 
\def\ext{\mbox{\rm ext}\,}
\def\gam{w} 
\def\gamA{m} 
\def\gamB{r} 
\def\gamC{q} 
\def\Inv{\mbox{\rm Inv}\,} 
\def\latt{L} 
\def\mrg{\varphi} 
\def\nor{{\mathsf N}} 
\def\norm{{\cal N}} 
\def\supmoN{\lozenge(N)}  
\def\tight{\mathcal{T}} 
\def\tightB{\mathcal{S}} 
\def\EnPart{\mathbf{\Upsilon}} 
\def\FanPos{\mathfrak{F}} 
\def\TiStr{\mathbb{T}} 
\def\InStr{\mathcal{I}} 
\def\PerSG{G} 
\def\indset{\mathfrak{I}} 

\def\calF{{\cal F}}
\def\calX{{\cal X}}
\def\calY{{\cal Y}}
\def\calD{{\cal D}}
\def\contr{\angle\,}  
\begin{document}
\addtolength{\baselineskip}{1mm}
\parskip 0mm

\title{On combinatorial descriptions of faces of\\
the cone of supermodular functions\thanks{This research, unfortunately, has not been
supported financially by a grant GA\v{C}R.}}
\author{{\Large Milan Studen\'{y}}\\[0.3ex]
Czech Academy of Sciences,\\
Institute of Information Theory and Automation}
\date{\today} 
\maketitle

\begin{abstract}
This paper is to relate five different ways of combinatorial description
of non-empty faces of the cone of supermodular functions on the power set of a finite basic set $N$.
Instead of this cone we consider its subcone of {\em supermodular games\,}; it is also a polyhedral cone and
has the same (= isomorphic) lattice of faces. This step allows one to associate supermodular games with
certain polytopes in ${\dv R}^{N}$, known as {\em cores\/} (of these games) in context of cooperative game
theory, or {\em generalized permutohedra\/} in context of polyhedral geometry. Non-empty faces of the supermodular cone
then correspond to normal fans of those polytopes. This (basically) geometric way of description of faces
of the cone then leads to the combinatorial ways of their description.
The first combinatorial way is to identify the faces with certain partitions of the set of enumerations of $N$,
known as {\em rank tests\/} in context of algebraic statistics.
The second combinatorial way is to identify faces with certain collections of posets on $N$,
known as (complete) {\em fans of posets\/} in context of polyhedral geometry.
The third combinatorial way is to identify the faces with certain coverings of the power set of $N$,
introduced relatively recently in context of cooperative game theory under name {\em core structures\/}.
The fourth combinatorial way is to identify the faces with certain formal {\em conditional independence structures\/},
introduced formerly in context of multivariate statistics under name {\em structural semi-graphoids}.
The fifth way is to identify the faces with certain {\em subgraphs\/} of the permutohedral graph,
whose nodes are enumerations of $N$. We prove the equivalence of those six ways of description of non-empty faces
of the supermodular cone. This result also allows one to describe the faces of the polyhedral cone of
(rank functions of) {\em polymatroids\/} over $N$ and the faces of the submodular cone
over $N$: this is because these cones have the same face lattice as the supermodular cone over $N$.
The respective polytopes are known as {\em base polytopes\/} in context of optimization and
(poly)matroid theory.
\end{abstract}


\section{Introduction: motivation}\label{sec.introduction}
This is to explain the motivation for the research effort discussed here.
A future ambitious goal is to characterize extreme rays of the cone
of {\em standardized supermodular functions\/} in such a way that
one is able to determine or compute those rays in case of $6$ or more
basic variables; compare with \cite{CC24}. Our motivation for this goal comes from the following
two areas of interest.
\begin{itemize}
\item In the research on {\em conditional independence\/} (CI)
{\em structures\/} \cite{Stu05} the above mentioned extreme rays allow
one to characterize 
structures known as {\em structural semi-graphoids}.
\item In the context of {\em information theory\/}, in particular
in the research on inequalities valid for the {\em entropy function\/},
some methods are based on the knowledge of the extreme rays of the cone of {\em polymatroids}, which is affinely equivalent to the supermodular cone.
\end{itemize}
Besides that, there are other areas which can benefit
from a positive solution to the above mentioned task. Let us mention the
research on {\em cooperative coalitional games\/} \cite{vNM44},
the research on special polytopes known as {\em generalized permutohedra\/}
\cite{PRW08}, the area of {\em imprecise probabilities\/} \cite{Walley91},
and the area of {\em combinatorial optimization\/} \cite{Fuj91}.
\smallskip

The above mentioned task to determine extreme (standardized)
supermodular functions is difficult because the number of these functions
seems to grow at least super-exponentially with the number of basic variables.
Indeed, while in case of 3 basic variables one has only 5 extreme rays and in
case of 4 variables 39 rays, long computation in case of 5 variables resulted in 117978
extreme rays \cite{SBK00}, while the case of 6 variables
is not manageable by present computational software packages for linear programming tasks.
\smallskip

On the other hand, since the (extreme) supermodular functions (obtained in the case of a low number of basic variables)
show certain symmetries and admit a combinatorial interpretation, an alternative approach may appear to be feasible.
The idea is to associate those extreme rays with appropriate combinatorial objects and try to characterize those combinatorial objects.
If this works then one can possibly replace linear programming computation of extreme rays by an alternative computational procedure based on the combinatorial interpretation of the rays.
\smallskip

What is discussed in this treatment is the first gradual step on the way
to the intended alternative computational procedure. We prove that the extreme rays of the supermodular cone do have relevant combinatorial interpretation; actually, they have several mutually equivalent combinatorial descriptions. In addition to that, the combinatorial
interpretations are not limited to the extreme rays of the cone, they even work for all non-empty faces of the supermodular cone.
The discussed combinatorial objects have already appeared in the literature
but their essential role in our context has not been properly acknowledged.

More specifically, one version of these objects are {\em rank tests\/} discussed in Morton's thesis \cite{Mor07}
and in his later joint 2009 paper \cite{MPSSW09}. Another version of these objects are {\em fans of posets\/}
discussed in 2008 paper \cite{PRW08}. Yet another version of these objects are {\em core structures\/} introduced in 2010 in context of game theory \cite{KVV10}, which also
seem to correspond to ``collections of sublattices" discussed in 1991 by Fujishige \cite[\S\,3.3(d)]{Fuj91}. The formerly mentioned {\em structural semi-graphoids\/} \cite{Stu05} can also be viewed
as combinatorial objects of this kind although they do not seem to offer clear computational advantages.
Related geometric objects are {\em generalized permutohedra\/} discussed by Postnikov and others \cite{PRW08} and
the normal fans of these polytopes.
\smallskip

Note in this context that it was surprising for the author to recognize
combinatorial structures corresponding to extreme rays of the supermodular
cone because of several former unsuccessful attempts to classify
those rays using combinatorial tools; see Section~\ref{ssec.dead-ends}.

\section{Preliminaries: four lattices}\label{sec.preliminary}
Let $N$ be a finite non-empty {\em basic set\/} of variables; the number of its elements will
be denoted by $n := |N|\geq 1$. In this paper we intentionally regard the basic set $N$ as an
\underline{unordered} set. The subtle reasons for that, not evident at first sight, are explained in Section~\ref{ssec.why-enumeration}. On contrary to that, given a natural number $n\in {\dv N}$, the symbol $[n]\,:=\,\{1,\ldots, n\}$ will denote an \underline{ordered} set of integers between $1$ and $n$.

The power set of $N$ will be denoted by $\caP := \{ A: A\subseteq N\}$.
The {\em diagonal\/} in the Cartesian product $N\times N$ 
will be denoted by $\diag:=\{\, (u,u)\,:\ u\in N\}$.
Given $\ell\in N$, the symbol $\ell$ will also be used to denote the one-element set (= singleton) $\{ \ell\}$.
Given $A,B\subseteq N$, the symbol $AB$ will be a shortened notation for the union $A\cup B$.
The symbol $\delta(\star\star)$, where $\star\star$\, is a predicate (= statement),
will occasionally denote a zero-one function whose value is $+1$ if the statement $\star\star$ is valid and $0$ otherwise.
The {\em incidence vector\/} $\chi_{A}$ for a subset $A\subseteq N$,
is thus a vector $[x_{\ell}]_{\ell\in N}$ in ${\dv R}^{N}$
specified by $x_{\ell}:=\delta(\,\ell\in A\,)$ for $\ell\in N$.

In the main part of the paper we assume that the reader is familiar with several concepts
and basic facts both from lattice theory and polyhedral geometry. Nevertheless, for reader's convenience,
the definitions of those concepts (and relevant basic facts) are recalled in Appendix~\ref{sec.app.lattice} and Appendix~\ref{sec.app.polyhedral}.

\subsection{Enumerations and permutohedron}\label{ssec.enum-permutohedron}
Ordered lists of (all) elements of $N$ (without repetition) will be called {\em enumerations\/} (of $N$).
Formally, an enumeration of $N$ is a one-to-one (= bijective) mapping $\permA : [n]\to N$ from $[n]$ onto $N$.
We will use the record $|\,\permA(1)|\ldots|\permA(n)\,|$ to specify such an enumeration.
Every enumeration of $N$ can be interpreted as a total order (= toset) on $N$. In fact, there are two
ways to do that: one can either interpret $\permA$ in the ascending way as $\permA(1)\prec \ldots\prec \permA(n)$
or one may interpret it in the descending way as $\permA(1)\succ \ldots\succ \permA(n)$.

The set of (all) enumerations of 
$N$ will be denoted by $\enuN$.
Every enumeration can equivalently be described by the corresponding {\em rank vector\/}, which
is the vector $[\permA_{-1}(\ell)]_{\ell\in N}$ in ${\dv R}^{N}$. Formally, the rank vector (for an enumeration $\permA$) is
the inverse $\permA_{-1}: N\to [n]$ of the mapping $\permA : [n]\to N$. We are going
to denote it by $\rankv_{\permA}$ in a geometric context (of ${\dv R}^{N}$).

The {\em permutohedron\/} in ${\dv R}^{N}$, denoted by $\perN$, is defined as the convex hull
of the rank vectors for (all) enumerations of $N$. Formally,
$\perN :=\conv (\,\{\rankv_{\permA}\,:\ \permA\in\enuN\}\,)$. It is a subset of
the affine space $\{\, y\in {\dv R}^{N}\,:\ \sum_{\ell\in N} y_{\ell}={n+1\choose 2}\,\}$; therefore,
its dimension is $n-1$. The polyhedral description of the permutohedron in ${\dv R}^{N}$ is as follows:
$$
\perN ~=~ \{\, z\in {\dv R}^{N}\ :\ \sum_{\ell\in N} z_{\ell}= {n+1\choose 2} ~~\&~~
\sum_{\ell\in S} z_{\ell}\geq {|S|+1\choose 2}\quad \mbox{for every $S\subseteq N$}\,\}.
$$
This fact follows from a later stronger result on vertices of the core polytope
from Section~\ref{ssec.marginal-vectors} \cite{Sha72,SK16}.
In particular, the facets of\/ $\perN$ correspond to non-empty proper subsets of $N$.

\subsection{Face lattice of the permutohedron}\label{ssec.face-lattice}
Recall from Appendix~\ref{ssec.app.face-lattice} that the set of (all) faces of any polytope,
ordered by inclusion, is a finite lattice, known to be both atomistic and coatomistic.
In addition to that, it is a graded lattice, where the dimension (of faces) plays the role
of a height function.

In this subsection we present an elegant {\em combinatorial\/} description of (non-empty) faces of
the permutohedron $\perN$, mentioned in \cite[p.\,18]{Zie95}, whose proof can be found in \cite[\S\,1]{BS96}.
Note that the proof is based on the fact that the permutohedron is a simple polytope.

Every non-empty face of the permutohedron $\perN$ corresponds to an ordered partition
of the basic set $N$ into non-empty blocks. We will use the symbol $|\,A_{1}|\ldots|A_{m}\,|$ to
denote such a partition into $m$ blocks, $1\leq m\leq n$. Note that this notation is consistent
with our notation for enumerations. In fact, the ordered partition serves as a compact
description of the set of those enumerations of $N$ which encode the {\em vertices\/} of the
corresponding face of $\perN$: these are the respective rank vectors.
Specifically, it is the set of those enumerations $\permA$ of $N$ which are consistent with $|\,A_{1}|\ldots|A_{m}\,|$:
these are the ordered lists $\permA$ of elements of $N$ where first (all) the elements of $A_{1}$ are listed,
then (all) the elements of $A_{2}$ follow and so on.
Thus, any face of\/ $\perN$ also corresponds to a special subset of\/ $\enuN$.
The dimension of the face described by $|\,A_{1}|\ldots|A_{m}\,|$ is $n-m$ then.
Later, in Section~\ref{ssec.poset-lattice}, we interpret such (non-empty) face-associated subsets of $\enuN$
as a special case of subsets of\/ $\enuN$ assigned to partial orders on $N$.

In particular, (if $|N|\geq 2$ then) the facets of\/ $\perN$ are described by (ordered) partitions
of $N$ into two non-empty blocks  and correspond to proper subsets of $N$, as noted in Section~\ref{ssec.enum-permutohedron}.

\subsection{A lattice based on posets}\label{ssec.poset-lattice}
The subject of our interest will be particular partitions of the set $\enuN$ of enumerations of $N$ whose components correspond to posets on $N$. The subsets of $\enuN$ which correspond to posets are non-empty and
form a (join) semi-lattice which can be turned into a lattice by adding the empty set (of enumerations). In this section we define this special (extended) lattice by the method of Galois connections, as described in
Appendix~\ref{ssec.app.Galois}.

To this end we introduce a binary relation $\preced$ between elements of the set $X:=\enuN$ of enumerations
of $N$ and the elements of the Cartesian product $Y:=N\times N$. More specifically, if $\permA :[n]\to N$ is an enumeration of $N$ and $(u,v)\in N\times N$ then we consider $\permA$ to be in an incidence relation $\preced$ with the pair $(u,v)$ iff {\em $u$ precedes $v$\/} in $\permA$\,:
$$
\permA\,\preced\, (u,v) \quad :=\quad \permA_{-1}(u)\leq \permA_{-1}(v)\,,
\quad \mbox{written also in the form $u\preceq_{\permA} v$.}
$$
We will denote Galois connections based on this precedence relation $\preced$
using larger triangles $\rGprec$ and $\lGprec$, in order to distinguish them
from other important Galois connections.
Thus, any enumeration $\permA\in\enuN$ of $N$ is assigned (by forward Galois connection) a toset
$$
T_{\permA} := (\{\permA\})^{\rGprec} =
\{\, (u,v)\in N\times N\, :\ u\preceq_{\permA} v\,\},
\qquad \mbox{denoted alternatively by $\permA^{\rGprec}$}.
$$
We also introduce special notation for the case that {\em $u$ strictly precedes $v$\/} in $\permA$\,:
$$
u\prec_{\permA} v ~~ \mbox{will mean that~ $\permA_{-1}(u) < \permA_{-1}(v)$.}
$$
Any subset $T$ of\/ $Y=N\times N$ is assigned (by backward Galois connection) an enumeration set
$$
{\cal L}(T):= T^{\lGprec}= \{\, \permA\in\enuN\,:\ \forall\,(u,v)\in T~~ u\preceq_{\permA} v\,\}\,,
$$
which is the set of {\em linear extesions\/} for $T$ (provided $T$ is a subset of a poset on $N$, as otherwise ${\cal L}(T)$ turns out to be empty).
If\/ $T$ is disjoint with the diagonal $\diag$ on $N\times N$ then it could be represented (= depicted) by a directed graph having $N$
as the set of nodes, where an arrow $u\to v$ from $u$ to $v$ is made iff $(u,v)\in T$.

\begin{lemma}\label{lem.poset-lattice}\rm
Let $({\calX}^{\preced},\subseteq )$ and $({\calY}^{\preced},\subseteq )$ denote the finite lattices defined using the Galois connections $\rGprec$ and $\lGprec$ based on the incidence relation $\preced$ between $X$ and $Y$ introduced above.
\begin{itemize}
\item[(i)] If $n=|N|\geq 2$ then\, [\,$T\in {\calY}^{\preced}$ and\, $T\subset N\times N$\,]
iff\/ $T$ is a poset on $N$.
\item[(ii)] One has [\,$S\in {\calX}^{\preced}$ and $S\neq\emptyset$\,] iff
a (uniquely determined) poset $(N,\preceq)$ exists such that $S=\{\,\permA\in\enuN\,:\
\permA_{-1}(u)<\permA_{-1}(v)~~ \mbox{whenever $u\prec v$}\,\}$, i.e.\  $S={\cal L}(T)$ for a poset~$T$ on $N$.
\item[(iii)] Any $S\subseteq\enuN$ representing a non-empty face of\/ $\Pi (N)$
(as in Section\,\ref{ssec.face-lattice})  belongs to ${\calX}^{\preced}$.
\item[(iv)] The lattice $({\calX}^{\preced},\subseteq )$ is
atomistic, coatomistic, and graded.
\end{itemize}
Coatoms in ${\calX}^{\preced}$ are precisely the sets $S_{u\prec v}\,:=\,\{\,\permA\in\enuN \,:\ u\prec_{\permA}v\,\}$,
where $u,v\in N$, $u\neq v$.
\end{lemma}

Note in this context that the class of sets ${\calX}^{\preced}\setminus\{\emptyset\}$ can alternatively
be introduced by means of another ``strict" incidence relation $\permA\odot(u,v)$
defined by $\permA_{-1}(u)<\permA_{-1}(v)$.
Except the degenerate case $|N|=1$, the resulting lattice $({\calX}^{\odot},\subseteq)$
coincides with $({\calX}^{\preced},\subseteq)$ while the elements of\/ ${\calY}^{\odot}$
are the corresponding irreflexive relations on $N$ instead.

\begin{proof}
For the necessity in (i) assume $T\in{\calY}^{\preced}$ and\/ $T\neq N\times N$. Therefore,
$$
T=S^{\rGprec}:=\{\, (u,v)\in N\times N :\
\forall\,\permA\in S~~ \permA_{-1}(u)\leq\permA_{-1}(v)\,\}
$$
for some $S\subseteq\enuN$ (see Section~\ref{ssec.app.Galois}).
The reflexivity and transitivity of\/ $T$ follows directly from its definition on basis of $S$.
Then $T$ is anti-symmetric as otherwise $S=\emptyset$ gives $T=N\times N$.

For the sufficiency in (i) assume that $T\subset N\times N$ is a poset. As recalled in Appendix~\ref{ssec.app.trans-DAG}, $T\setminus\diag$ is then acyclic and coincides with the set of arrows of
a transitive directed acyclic graph $G$ over $N$. Hence, the set
\begin{eqnarray*}
S= T^{\lGprec} &=& \{\, \permA\in\enuN\,:\ \forall\,(u,v)\in T~~ \permA_{-1}(u)\leq\permA_{-1}(v)\,\}\\
&=& \{\, \permA\in\enuN\,:\  \permA_{-1}(a)<\permA_{-1}(b)\quad
\mbox{whenever $a\to b$ in $G$}\,\}
\end{eqnarray*}
of enumerations of $N$ which are consonant with $G$ is non-empty. Evidently, $T\subseteq S^{\rGprec}$ and to show $S^{\rGprec}\subseteq T$ one needs, for any $(u,v)\not\in T$, to find an enumeration $\permA\in S$ with $\permA_{-1}(u)>\permA_{-1}(v)$.
Since $S\neq\emptyset$, this is evident in case $v\to u$ in $G$.
In case of a pair $(u,v)$ of non-adjacent distinct nodes of $G$ one can consider
the set $A\subseteq N$ of ancestors of $v$ in $G$. Then $v\in A$, $u\not\in A$ by transitivity of $G$, and
the induced subgraph $G_{A}$ is acyclic. Hence, an enumeration of $A$ exists which is
consonant with $G_{A}$ and this can be extended to an enumeration $\permA$ of $N$ consonant with $G$.
Thus, $\permA\in S$ with $\permA_{-1}(v)<\permA_{-1}(u)$ was obtained.
Altogether, we have shown $T=S^{\rGprec}$, which means that $T\in {\calY}^{\preced}$.

As concerns (ii), this holds trivially in case $n=1$. In case $n\geq 2$ apply (i) and realize the fact that $({\calX}^{\preced},\subseteq )$ and $({\calY}^{\preced},\subseteq )$
are anti-isomorphic; see Appendix~\ref{ssec.app.Galois}. Specifically, the correspondence $T\in {\calY}^{\preced}\mapsto T^{\lGprec}\in {\calX}^{\preced}$ maps bijectively the partial orders $\preceq$ on $N$ (= posets on~$N$) onto non-empty sets in ${\calX}^{\preced}$.

Note (iii) holds trivially if $n=1$. In case $n\geq 2$ consider a non-empty
face of $\Pi (N)$ represented by an ordered partition
$|\,A_{1}|\ldots|A_{m}\,|$, $m\geq 1$, of $N$. It defines a partial order
$$
T:=\diag\cup \{\, (u,v)\in N\times N\,:\ u\in A_{i} ~\mbox{and}~ v\in A_{j}~~ \mbox{for some $i<j$}\,\,\}
$$
on $N$. Easily, ${\cal L}(T)$ coincides with the set of $\permA\in\enuN$ 
consistent with $|\,A_{1}|\ldots|A_{m}\,|$.
\smallskip

As concerns (iv), this holds trivially in case $n=1$. In case $n\geq 2$ consider $\permA\in\enuN$ with its toset $T_{\permA}$ and observe that $(T_{\permA})^{\lGprec}=\{\,\permA\,\}$.
Hence, the atoms of $({\calX}^{\preced},\subseteq )$ are just singleton subsets of\/ $\enuN$. Since any $S\in {\calX}^{\preced}$ is the union of singletons, the lattice $({\calX}^{\preced},\subseteq )$ is atomistic.

The fact that $({\calX}^{\preced},\subseteq )$ is coatomistic follows from the fact that $({\calY}^{\preced},\subseteq )$ is atomistic. To show the latter fact realize that, for
every $(u,v)\in (N\times N)\setminus\diag$, the set $\diag\cup \{\,(u,v)\,\}$ is a poset on~$N$.
Thus, using (i), the atoms of $({\calY}^{\preced},\subseteq )$ correspond to
singleton subsets of $(N\times N)\setminus\diag$. This allows one to observe, again
using (i), that $({\calY}^{\preced},\subseteq )$ is atomistic.

Note that $({\calX}^{\preced},\subseteq )$ is graded iff
$({\calY}^{\preced},\subseteq )$ is graded. To verify the latter fact one can use the characterization of\/ ${\calY}^{\preced}$ from (i).
The anti-isomorphism of $({\calX}^{\preced},\subseteq )$ and
$({\calY}^{\preced},\subseteq )$ implies that the coatoms of $({\calY}^{\preced},\subseteq )$
are the tosets $T_{\permA}$, $\permA\in\enuN$, where every
$T_{\permA}\setminus\diag$ has the cardinality ${n\choose 2}$.
As explained in Appendix~\ref{ssec.app.graded},
to reach our goal it is enough to show that any maximal
chain in $({\calY}^{\preced}\setminus\{Y\},\subseteq )$
has the length ${n\choose 2}$, that is, ${n\choose 2}+1$ elements.

This follows by repeated application of the following {\em sandwiche
principle}. Given two transitive irreflexive relations  $G^{\prime}\subset G$ on $N$ (= strict versions of posets on $N$) there exists a transitive irreflexive relation $G^{\prime\prime}$ on $N$
with $G^{\prime}\subseteq G^{\prime\prime}\subset G$ and $|G\setminus G^{\prime\prime}|=1$.
To verify the principle realize that, for the strict order $\prec$ defined by $G$, there exists $(u,v)\in N\times N$ with $u\cover v$ (in $G$) and $(u,v)\not\in G^{\prime}$.
Indeed, otherwise the covering relation $\cover$ for $G$ is contained in $G^{\prime}$ implying that its transitive closure, which is $G$, is also contained, yielding a contradictory conclusion $G\subseteq G^{\prime}$.
Take a pair $(u,v)$ with $u\cover v$ and $(u,v)\not\in G^{\prime}$ and put $G^{\prime\prime}:=G\setminus\{\,(u,v)\,\}$.
The transitive closure of\/ $G^{\prime\prime}$ is contained in (a transitive relation) $G$. Nevertheless, $(u,v)$ cannot be in the transitive closure of\/ $G^{\prime\prime}$, because this contradicts the assumption that $u\cover v$ in $G$. Thus, $G^{\prime\prime}$ is both transitive and irreflexive, which
verifies the sandwich principle.

The claim about coatoms in ${\calX}^{\preced}$ follows from the description of atoms in $({\calY}^{\preced},\subseteq)$:
these have the form $A_{(u,v)}=\diag\cup \{\, (u,v)\,\}$ for $u,v\in N$, $u\neq v$, and $A_{(u,v)}^{\lGprec}=\{\, (u,v)\,\}^{\lGprec}=S_{u\prec v}$.
\end{proof}

Note in this context what is the closure operation corresponding to the lattice $({\calY}^{\preced},\subseteq)$.
Given $T\subseteq Y=N\times N$, one can show that
$$
T^{\lGprec\rGprec} ~=~
\left\{
\begin{array}{cl}
\tr(T\cup\diag) & \mbox{if $T\setminus\diag$ is acyclic},\\
N\times N & \mbox{otherwise},
\end{array}
\right.
~~ \mbox{where $\tr(T)$ is the {\em transitive closure\/} of\/ $T$.}
$$
Recall that, in the former case, $\tr(T\cup\diag)\setminus\diag$ is both transitive and acyclic.

\subsubsection{Graphical characterization of poset-based sets of enumerations}\label{sssec.poset-lattice-characterization}
There is an elegant characterization of sets in the lattice $({\calX}^{\preced},\subseteq )$ in graphical terms (if $n\geq 2$).
We say that enumerations $\permA,\permB\in\enuN$ differ by an {\em adjacent transposition\/} if
\begin{eqnarray}
\lefteqn{\hspace*{-2cm}\exists\, 1\leq i<n ~~~\mbox{such that}~~ \permA(i)=\permB(i+1),~ \permA(i+1)=\permB(i),} \label{eq.adjacent-transpos}\\
&\mbox{and}&
\permA(k)=\permB(k)\qquad \mbox{for remaining $k\in [n]\setminus\{i,i+1\}$}\,. \nonumber
\end{eqnarray}
We are going to interpret the set of enumerations of $N$ as an undirected {\em permutohedral graph\/}
in which $\enuN$ is the set of nodes and edges are determined by adjacent transpositions. %
The reason for this terminology is as follows. Despite the given definition of adjacency
for enumerations is fully combinatorial,
it has a natural geometric interpretation: one can show that $\permA,\permB\in\enuN$ are adjacent
in this graph iff the segment $[\rankv_{\permA},\rankv_{\permB}]\subseteq {\dv R}^{N}$ between the respective rank vectors is a face (= a geometric edge) of the permutohedron $\perN$.
Note that this observation can also be derived from the characterization of faces of the permutohedron
$\perN$ recalled in Section~\ref{ssec.face-lattice}.

The point is that poset-based subsets of\/ $\enuN$ can equivalently be characterized solely in terms of this particular graph. The following definitions apply to any connected undirected graph with
the set of nodes $U$. A {\em geodesic\/} between nodes $\permA,\permB\in U$ is a walk between $\permA$ and $\permB$ (in the graph) which has the shortest possible length among such walks.
Note that it is necessarily a path (= nodes are not repeated in it)
and that several geodesics may exist between two nodes. The (graphical) {\em distance\/} $\dist(\permA,\permB)$ between two nodes $\permA,\permB\in U$
is the length of a geodesic between them. We say that a node $\permC\in U$ {\em is between nodes\/} $\permA\in U$ and $\permB\in U$ if $\permC$ belongs to some geodesic between $\permA$ and $\permB$; an equivalent condition is that $\dist(\permA,\permC)+\dist(\permC,\permB)= \dist(\permA,\permB)$.
A set $S\subseteq U$ is {\em geodetically convex\/} if, for any $\permA,\permB\in S$, all nodes between them belong to $S$. It is evident that the intersection of geodetically convex sets is a geodetically convex set.

To prove our result saying that (if $n\geq 2$) sets in the lattice ${\calX}^{\preced}$ coincide with geodetically convex sets in the permutohedral graph it is instrumental to introduce a particular {\em edge-labeling\/}  (= edge coloring)
for this graph. Specifically, we are going to label its edges by two-element subsets $\{u,v\}$ of
our basic set $N$: if $\permA,\permB\in\enuN$ differ by an adjacent transposition
\eqref{eq.adjacent-transpos} then we label the edge between $\permA$ and $\permB$ by the set
$\{u,v\}:=\{\permA(i),\permA(i+1)\}=\{\permB(i),\permB(i+1)\}$.

Observe that every two-element subset $\{u,v\}$ of $N$ defines an {\em automorphism\/} of the
permutohedral graph. Specifically, the transposition $\trans_{uv}:N\to N$ of $u$ and $v$,
defined by $\trans_{uv}(u)=v$, $\trans_{uv}(v)=u$, and $\trans_{uv}(t)=t$ for $t\in N\setminus\{u,v\}$,
is a permutation on $N$ and every enumeration $\permG :[n]\to N$ can be assigned its composition
$\permG\trans_{uv}$ with $\trans_{uv}:N\to N$. The mapping $\permG\mapsto \permG\trans_{uv}$
is a self-inverse transformation of\/ $\enuN$ preserving the edges in the permutohedral graph. In fact, it maps
bijectively the set
$$
S_{u\prec v}\,=\,\{\,\permG\in\enuN \,:\ u\prec_{\permG}v\,\}
\quad\mbox{to its complement}\quad
S_{v\prec u}\,=\,\{\,\permG\in\enuN \,:\ v\prec_{\permG}u\,\}\,.
$$
Given $\permA,\permB\in\enuN$, a two-element subset $\{u,v\}$ of $N$ will be called
an {\em inversion between\/} $\permA\in\enuN$ and $\permB\in\enuN$ if the mutual orders of $u$ and $v$ in $\permA$ and $\permB$ differ:
\begin{eqnarray*}
\Inv[\permA,\permB] ~:=~
\left\{\,\, \{u,v\}\subseteq N\, :~
\mbox{either $[\, u\prec_{\permA}v ~\&~ v\prec_{\permB}u\,]$ ~or~
$[\, v\prec_{\permA}u ~\&~ u\prec_{\permB}v\,]$}\,\,\right\}\,.
\end{eqnarray*}
The next lemma brings basic observations on the permutohedral graph.

\begin{lemma}\label{lem.permut-graph}\rm
Consider $\permA,\permB\in\enuN$. Then
\begin{itemize}
\item[(i)] $\dist(\permA,\permB)=|\,\Inv[\permA,\permB]\,|$ is the number of inversions
between $\permA$ and $\permB$. In particular, $\permA$ and $\permB$ are adjacent
in the permutohedral graph iff they differ by one inversion only. The label on the edge between $\permA$ and $\permB$ is then the only inversion between them.
\item[(ii)]
A walk in the permutohedral graph is a geodesic (between $\permA$ and $\permB$) iff no label 
on its edges is repeated. If this is the case then the set of labels on the geodesic coincides
with the set $\Inv[\permA,\permB]$ of inversions between $\permA$ and $\permB$.
\item[(iii)]
An enumeration $\permC\in\enuN$ is between $\permA$ and $\permB$ iff $\Inv[\permA,\permC]\cap\Inv[\permC,\permB]=\emptyset$.
Another equivalent condition is $T_{\permA}\cap T_{\permB}\subseteq T_{\permC}$, that is,
$[\, u\prec_{\permA} v ~\&~  u\prec_{\permB} v\,] ~\Rightarrow ~ u\prec_{\permC} v$.
\end{itemize}
\end{lemma}

\begin{proof}
Assume $n\geq 2$ as otherwise all claims are trivial.
Let us start with the (particular) observation in (i) that $\permA,\permB\in\enuN$ are adjacent
iff they differ by one inversion only. The necessity follows immediately from \eqref{eq.adjacent-transpos}
and for the sufficiency assume that $\{u,v\}\in\Inv[\permA,\permB]$ is unique, which means
that any other pair of elements of $N$ has the same mutual order in $\permA$ and in $\permB$.
Assume specifically that $u\prec_{\permA}v$ and $v\prec_{\permB}u$ as otherwise we exchange
$\permA$ for $\permB$.
Thus, one has $u=\permA(i)$ and $v=\permA(j)$ with $i<j$.
One has $j=i+1$ as otherwise there is $t\in N$ with $u\prec_{\permA}t\prec_{\permA}v$ implying
a contradictory conclusion $u\prec_{\permB}t\prec_{\permB}v\prec_{\permB}u$.
Moreover, for any $t\in N\setminus\{u,v\}$, one has $t\prec_{\permA}u \,\Leftrightarrow\, t\prec_{\permB}u$, and
$v\prec_{\permA}t \,\Leftrightarrow\, v\prec_{\permB}t$. These facts enforce $v=\permB(i)$ and $u=\permB(i+1)$.
The fact that the mutual orders in $\permA$ and $\permB$ coincide among $t\in N$ with $t\prec_{\permA}u$
and also among $t\in N$ with $v\prec_{\permA}t$ then implies $\permA(k)=\permB(k)$ for $k\in[n]\setminus\{i,i+1\}$.
Thus, \eqref{eq.adjacent-transpos} holds.
Note that the label $\{u,v\}$ on the edge between enumerations $\permA$ and $\permB$ is just the only
inversion between them.

The latter observation implies that, for every walk between enumerations $\permA,\permB\in\enuN$
and every $\{u,v\}\in\Inv[\permA,\permB]$, there is an edge of the walk labeled by $\{u,v\}$, as otherwise
$u$ and $v$ would have the same mutual order in $\permA$ and in $\permB$. Hence, $|\,\Inv[\permA,\permB]\,|$
is a lower bound for the length of a walk between $\permA$ and $\permB$.

Therefore, to verify (i) it is enough to prove, for $\permA,\permB\in\enuN$, by induction on $|\,\Inv[\permA,\permB]\,|$, that there exists a walk between them which has the length $|\,\Inv[\permA,\permB]\,|$.
If there is no inversion between $\permA$ and $\permB$ then
$\permA(1)\prec_{\permA} \permA(2)\, \ldots \prec_{\permA} \permA(n)$ gives
$\permA(1)\prec_{\permB} \permA(2)\, \ldots \prec_{\permB} \permA(n)$ and $\permA=\permB$,
that is, a walk of the length $0=|\,\Inv[\permA,\permB]\,|$ exists.
If $|\,\Inv[\permA,\permB]\,|\geq 1$ then necessarily $\permA\neq\permB$ and
$\exists\, 1\leq i<n$ such that one has $\permB(i+1)\prec_{\permA} \permB(i)$, for
otherwise $\permB(1)\prec_{\permA} \permB(2)\, \ldots \prec_{\permA} \permB(n)$ would give
a contradictory conclusion $\permB=\permA$. Define $\permC\in\enuN$ as the enumeration
obtained from $\permB$ by the respective adjacent transposition: $\permC(i)=\permB(i+1)=:u$,
$\permC(i+1)=\permB(i)=:v$, and $\permC(k)=\permB(k)$ for $k\in[n]\setminus\{i,i+1\}$.
Since $\{u,v\}$ is an inversion between $\permA$ and $\permB$,
one easily observes $\Inv[\permA,\permC]=\Inv[\permA,\permB]\setminus\{\, \{u,v\}\,\}$.
The induction premise implies the existence of a walk between $\permA$ and $\permC$ of the
length $|\,\Inv[\permA,\permC]\,|$ which can be prolonged by the edge between $\permC$ and $\permB$.
That verifies the induction step.

To prove the necessity in (ii) consider a geodesic $\mathfrak{g}$ between $\permA,\permB\in\enuN$
and assume for a contradiction that a label $\{u,v\}\subseteq N$ is repeated on its edges.
The last claim in (i) implies that the label cannot be repeated on consecutive edges of $\mathfrak{g}$:
if an edge between $\permD\in\enuN$ and $\permE\in\enuN$ has the same label as the
consecutive edge between $\permE$ and $\permC\in\enuN$ then $\permD=\permC$ and $\mathfrak{g}$ can be shortened,
contradicting the assumption that it is a geodesic.

Therefore, one can assume without loss of generality
that the first occurrence $\{u,v\}$ on the way from $\permA$ to $\permB$ is
between enumerations $\permD$ and $\permE$, with $\permD$ closer to $\permA$, and the next
occurrence of $\{u,v\}$ is between enumerations $\permC$ and $\permF$, with $\permF$ closer to $\permB$
(see Figure~\ref{fig.1}(a) for illustration). Assume specifically $u\prec_{\permA} v$, which enforces $\permD,\permF\in S_{u\prec v}$ and
$\permE,\permC\in S_{v\prec u}$. Let $\mathfrak{g}^{\prime}$ be the section of $\mathfrak{g}$
between $\permE$ and $\permC$. We know $\mathfrak{g}^{\prime}\subseteq S_{v\prec u}$ and, because the map
$\permG\mapsto\permG\trans_{uv}$ is an automorphism of the permutohedral graph, the section $\mathfrak{g}^{\prime}$
has a copy $\mathfrak{g}^{\prime\prime}$ in $S_{u\prec v}$, which is between $\permD$ and $\permF$
(see Figure~\ref{fig.1}(b) for illustration).
Therefore, the section of $\mathfrak{g}$ between $\permD$ and $\permF$ can be replaced by $\mathfrak{g}^{\prime\prime}$,
yielding a walk of a shorter length and contradicting the assumption that $\mathfrak{g}$ is a geodesic.

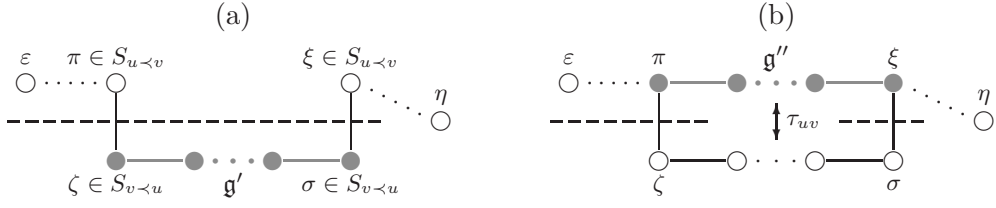
\begin{figure}[t]
\setlength{\unitlength}{0.8mm}
\begin{center}
\begin{picture}(79,33)
\put(5,20){\circle{3}}
\put(20,20){\circle{3}}
\put(20,7){\color{lgrey}\circle*{3}}
\put(33,7){\color{lgrey}\circle*{3}}
\put(46,7){\color{lgrey}\circle*{3}}
\put(59,7){\color{lgrey}\circle*{3}}
\put(59,20){\circle{3}}
\put(74,13.5){\circle{3}}
\put(5,24){\makebox(0,0){\footnotesize $\permA$}}
\put(20,3){\makebox(0,0){\footnotesize $\permE\in S_{v\prec u}$}}
\put(39.5,3){\makebox(0,0){\small $\mathfrak{g}^{\prime}$}}
\put(59,3){\makebox(0,0){\footnotesize $\permC\in S_{v\prec u}$}}
\put(20,24){\makebox(0,0){\footnotesize $\permD\in S_{u\prec v}$}}
\put(59,24){\makebox(0,0){\footnotesize $\permF\in S_{u\prec v}$}}
\put(74,17.5){\makebox(0,0){\footnotesize $\permB$}}
\put(39.5,31){\makebox(0,0){\rm (a)}}
\put(2,13.5){\line(1,0){2}}
\put(5,13.5){\line(1,0){2}}
\put(8,13.5){\line(1,0){2}}
\put(11,13.5){\line(1,0){2}}
\put(14,13.5){\line(1,0){2}}
\put(17,13.5){\line(1,0){2}}
\put(20,13.5){\line(1,0){2}}
\put(23,13.5){\line(1,0){2}}
\put(26,13.5){\line(1,0){2}}
\put(29,13.5){\line(1,0){2}}
\put(32,13.5){\line(1,0){2}}
\put(35,13.5){\line(1,0){2}}
\put(38,13.5){\line(1,0){2}}
\put(41,13.5){\line(1,0){2}}
\put(44,13.5){\line(1,0){2}}
\put(47,13.5){\line(1,0){2}}
\put(50,13.5){\line(1,0){2}}
\put(53,13.5){\line(1,0){2}}
\put(56,13.5){\line(1,0){2}}
\put(59,13.5){\line(1,0){2}}
\put(62,13.5){\line(1,0){2}}
\put(20,18.1){\line(0,-1){9.2}}
\put(59,18.1){\line(0,-1){9.2}}
\put(8.5,20){\circle*{0.4}}
\put(10.5,20){\circle*{0.4}}
\put(12.5,20){\circle*{0.4}}
\put(14.5,20){\circle*{0.4}}
\put(16.5,20){\circle*{0.4}}
\put(62.5,18.75){\circle*{0.4}}
\put(64.5,17.75){\circle*{0.4}}
\put(66.5,16.75){\circle*{0.4}}
\put(68.5,15.75){\circle*{0.4}}
\put(70.5,14.75){\circle*{0.4}}
\put(36.5,7){\color{lgrey}\circle*{0.8}}
\put(39.5,7){\color{lgrey}\circle*{0.8}}
\put(42.5,7){\color{lgrey}\circle*{0.8}}
\thicklines
\put(21.9,7){\color{lgrey}\line(1,0){9.2}}
\put(47.9,7){\color{lgrey}\line(1,0){9.2}}
\end{picture}
~~~~~
\begin{picture}(79,33)
\put(5,20){\circle{3}}
\put(20,20){\color{lgrey}\circle*{3}}
\put(33,20){\color{lgrey}\circle*{3}}
\put(46,20){\color{lgrey}\circle*{3}}
\put(59,20){\color{lgrey}\circle*{3}}
\put(20,7){\circle{3}}
\put(33,7){\circle{3}}
\put(46,7){\circle{3}}
\put(59,7){\circle{3}}
\put(74,13.5){\circle{3}}
\put(5,24){\makebox(0,0){\footnotesize $\permA$}}
\put(20,3){\makebox(0,0){\footnotesize $\permE$}}
\put(39.5,24){\makebox(0,0){\small $\mathfrak{g}^{\prime\prime}$}}
\put(59,3){\makebox(0,0){\footnotesize $\permC$}}
\put(20,24){\makebox(0,0){\footnotesize $\permD$}}
\put(59,24){\makebox(0,0){\footnotesize $\permF$}}
\put(74,17.5){\makebox(0,0){\footnotesize $\permB$}}
\put(44,13.5){\makebox(0,0){\footnotesize $\trans_{uv}$}}
\put(39.5,31){\makebox(0,0){\rm (b)}}
\put(2,13.5){\line(1,0){2}}
\put(5,13.5){\line(1,0){2}}
\put(8,13.5){\line(1,0){2}}
\put(11,13.5){\line(1,0){2}}
\put(14,13.5){\line(1,0){2}}
\put(17,13.5){\line(1,0){2}}
\put(20,13.5){\line(1,0){2}}
\put(23,13.5){\line(1,0){2}}
\put(26,13.5){\line(1,0){2}}
\put(50,13.5){\line(1,0){2}}
\put(53,13.5){\line(1,0){2}}
\put(56,13.5){\line(1,0){2}}
\put(59,13.5){\line(1,0){2}}
\put(62,13.5){\line(1,0){2}}
\put(20,18.1){\line(0,-1){9.2}}
\put(59,18.1){\line(0,-1){9.2}}
\put(39.5,16.5){\vector(0,-1){6}}
\put(39.5,10.5){\vector(0,1){6}}
\put(8.5,20){\circle*{0.4}}
\put(10.5,20){\circle*{0.4}}
\put(12.5,20){\circle*{0.4}}
\put(14.5,20){\circle*{0.4}}
\put(16.5,20){\circle*{0.4}}
\put(62.5,18.75){\circle*{0.4}}
\put(64.5,17.75){\circle*{0.4}}
\put(66.5,16.75){\circle*{0.4}}
\put(68.5,15.75){\circle*{0.4}}
\put(70.5,14.75){\circle*{0.4}}
\put(36.5,7){\circle*{0.6}}
\put(39.5,7){\circle*{0.6}}
\put(42.5,7){\circle*{0.6}}
\put(36.5,20){\color{lgrey}\circle*{0.8}}
\put(39.5,20){\color{lgrey}\circle*{0.8}}
\put(42.5,20){\color{lgrey}\circle*{0.8}}
\put(21.9,7){\line(1,0){9.2}}
\put(47.9,7){\line(1,0){9.2}}
\thicklines
\put(21.9,20){\color{lgrey}\line(1,0){9.2}}
\put(47.9,20){\color{lgrey}\line(1,0){9.2}}
\end{picture}
\caption{A picture illustrating the proof of Lemma~\ref{lem.permut-graph}.\label{fig.1}}
\end{center}
\end{figure}

To prove the sufficiency in (ii) assume that $\mathfrak{g}$ is a walk between $\permA$ and $\permB$
in which no label is repeated. We have already observed that any inversion $\{u,v\}\in\Inv[\permA,\permB]$
must occur as a label on $\mathfrak{g}$. Nevertheless, no two-element set $\{u,v\}\subseteq N$ which
is {\em not an inversion\/} between $\permA$ and $\permB$ can be a label on $\mathfrak{g}$:
indeed, otherwise, since only one edge of $\mathfrak{g}$ is labeled by $\{u,v\}$ then, the mutual orders
of $u$ and $v$ in $\permA$ and $\permB$ must differ, which contradicts the assumption
that $\{u,v\}$ is not an inversion between $\permA$ and $\permB$. Thus, we have verified
both that the length of $\mathfrak{g}$ is $|\,\Inv[\permA,\permB]\,|=\dist(\permA,\permB)$
and the claim that $\Inv[\permA,\permB]$ is the set of labels on $\mathfrak{g}$.
\smallskip

The first claim in (iii) then follows from (ii). If $\permC$ belongs to a geodesic between $\permA$
and $\permB$ then the labels on it between $\permA$ and $\permC$ must differ from the labels between $\permC$ and $\permB$.
Conversely, if $\Inv[\permA,\permC]$ and $\Inv[\permC,\permB]$ are disjoint then the concatenation
of a geodesic between $\permA$ and $\permC$ with a geodesic between $\permC$ and $\permB$ yields a geodesic
between $\permA$ and $\permB$.

Assuming $\Inv[\permA,\permC]\cap\Inv[\permC,\permB]=\emptyset$, the conditions
$u\preceq_{\permA}v$ and $u\preceq_{\permB}v$ together imply $u\preceq_{\permC}v$
as otherwise $\{u,v\}$ belongs to the intersection of the sets of inversions.
Conversely, if $T_{\permA}\cap T_{\permB}\subseteq T_{\permC}$ then
any hypothetic simultaneous inversion $\{u,v\}$ in $\Inv[\permA,\permC]\cap\Inv[\permC,\permB]$
must be compared equally in $\permA$ and $\permB$, and, therefore, in $\permC$, which is a contradiction.
\end{proof}

Two special concepts for a non-empty {\em geodetically convex\/} set\/
$\emptyset\neq S\subseteq\enuN$ in the permutohedral graph are needed. The set of {\em inversions in $S$\/} is
$$
\Inv (S) ~:=~ \left\{\,\, \{u,v\}\subseteq N\, :\
\mbox{there are adjacent $\permA,\permB\in S$ labeled by $\{u,v\}$}\,\, \right\},
$$
while the {\em covering relation\/} for $S$ is the following binary relation on $N$:
\begin{eqnarray*}
\lefteqn{\Cov (S) ~:=~ \{\,\, (u,v)\in N\times N\, :~ }\\
&&\qquad ~ \mbox{there are adjacent $\permA\in S$ and $\permB\in \enuN\setminus S$ labeled by $\{u,v\}$, where $u\prec_{\permA}v$}\,\, \}.
\end{eqnarray*}
Note that the latter terminology is motivated by the fact that $\Cov (S)$ appears to coincide with the covering relation $\cover$ for the poset on $N$ assigned to $S$ (see later Section~\ref{sssec.posets-remarks}).
These two notions are substantial in the proof of the next crucial result.

\begin{thm}\label{thm.poset-characterization}\em
If $n=|N|\geq 2$ then, given $S\subseteq\enuN$, $S\in {\calX}^{\preced}$ iff $S$ is geodetically convex.
\end{thm}

\begin{proof}
{\bf I.}\, To verify the necessity use Lemma~\ref{lem.poset-lattice}(iv) saying that
$({\calX}^{\preced},\subseteq)$ is a coatomistic lattice. Since the intersection of geodetically convex
sets is geodetically convex and so is $S=\enuN$ it is enough to verify that
the coatoms in $({\calX}^{\preced},\subseteq)$ are geodetically convex. These are sets of the form
$S_{u\prec v}$ for pairs $(u,v)\in N\times N$, $u\neq v$. We apply Lemma~\ref{lem.permut-graph}(ii):
given a geodesic between $\permA,\permB\in S_{u\prec v}$, the fact $\{u,v\}\not\in \Inv[\permA,\permB]$
implies that $\{u,v\}$ is not a label on an edge of the geodesic, and, thus, every node on the geodesic
belongs to $S_{u\prec v}$.

{\bf II.}\, The first step to verify the sufficiency is to observe that,
for every $(u,v)\in\Cov (S)$ and $\permC\in S$, one has $u\prec_{\permC} v$.
Suppose specifically that an edge between enumerations $\permA\in S$ and $\permB\in \enuN\setminus S$ is
labeled by $\{u,v\}$ and $u\prec_{\permA}v$, which implies $v\prec_{\permB}u$ (see Figure~\ref{fig.2}(a) for illustration).
Consider a geodesic $\mathfrak{g}$ between $\permC$ and $\permA$. It belongs to $S$, by the
assumption that $S$ is geodetically convex.

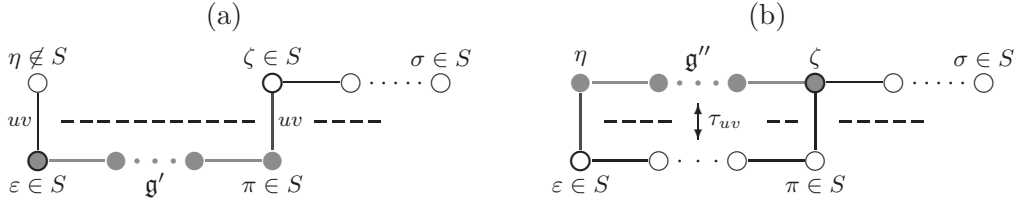
\begin{figure}[t]
\setlength{\unitlength}{0.8mm}
\begin{center}
\begin{picture}(79,33)
\put(7,20){\circle{3}}
\put(7,7){\color{lgrey}\circle*{3}}
\put(20,7){\color{lgrey}\circle*{3}}
\put(33,7){\color{lgrey}\circle*{3}}
\put(46,7){\color{lgrey}\circle*{3}}
\put(59,20){\circle{3}}
\put(74,20){\circle{3}}
\put(7,24){\makebox(0,0){\footnotesize $\permB\not\in S$}}
\put(7,3){\makebox(0,0){\footnotesize $\permA\in S$}}
\put(26.5,3){\makebox(0,0){\small $\mathfrak{g}^{\prime}$}}
\put(46,3){\makebox(0,0){\footnotesize $\permD\in S$}}
\put(46,24){\makebox(0,0){\footnotesize $\permE\in S$}}
\put(74,24){\makebox(0,0){\footnotesize $\permC\in S$}}
\put(4,13.5){\makebox(0,0){\scriptsize $uv$}}
\put(49,13.5){\makebox(0,0){\scriptsize $uv$}}
\put(38,31){\makebox(0,0){\rm (a)}}
\put(62.5,20){\circle*{0.4}}
\put(64.5,20){\circle*{0.4}}
\put(66.5,20){\circle*{0.4}}
\put(68.5,20){\circle*{0.4}}
\put(70.5,20){\circle*{0.4}}
%
\put(11,13.5){\line(1,0){2}}
\put(14,13.5){\line(1,0){2}}
\put(17,13.5){\line(1,0){2}}
\put(20,13.5){\line(1,0){2}}
\put(23,13.5){\line(1,0){2}}
\put(26,13.5){\line(1,0){2}}
\put(29,13.5){\line(1,0){2}}
\put(32,13.5){\line(1,0){2}}
\put(35,13.5){\line(1,0){2}}
\put(38,13.5){\line(1,0){2}}
\put(41,13.5){\line(1,0){2}}
\put(53,13.5){\line(1,0){2}}
\put(56,13.5){\line(1,0){2}}
\put(59,13.5){\line(1,0){2}}
\put(62,13.5){\line(1,0){2}}
\put(7,18.1){\line(0,-1){9.2}}
\put(46,18.1){\line(0,-1){9.2}}
\put(47.9,20){\line(1,0){9.2}}
\put(23.5,7){\color{lgrey}\circle*{0.8}}
\put(26.5,7){\color{lgrey}\circle*{0.8}}
\put(29.5,7){\color{lgrey}\circle*{0.8}}
\thicklines
\put(8.9,7){\color{lgrey}\line(1,0){9.2}}
\put(34.9,7){\color{lgrey}\line(1,0){9.2}}
\put(46,18.1){\color{hgrey}\line(0,-1){9.2}}
\put(7,7){\circle{3}}
\put(46,20){\circle{3}}
\end{picture}
~~~~~
\begin{picture}(79,33)
\put(20,7){\circle{3}}
\put(33,7){\circle{3}}
\put(46,7){\circle{3}}
\put(7,20){\color{lgrey}\circle*{3}}
\put(20,20){\color{lgrey}\circle*{3}}
\put(33,20){\color{lgrey}\circle*{3}}
\put(46,20){\color{lgrey}\circle*{3}}
\put(59,20){\circle{3}}
\put(74,20){\circle{3}}
\put(7,24){\makebox(0,0){\footnotesize $\permB$}}
\put(7,3){\makebox(0,0){\footnotesize $\permA\in S$}}
\put(26.5,24){\makebox(0,0){\small $\mathfrak{g}^{\prime\prime}$}}
\put(46,3){\makebox(0,0){\footnotesize $\permD\in S$}}
\put(46,24){\makebox(0,0){\footnotesize $\permE$}}
\put(74,24){\makebox(0,0){\footnotesize $\permC\in S$}}
\put(31,13.5){\makebox(0,0){\footnotesize $\trans_{uv}$}}
\put(38,31){\makebox(0,0){\rm (b)}}
\put(62.5,20){\circle*{0.4}}
\put(64.5,20){\circle*{0.4}}
\put(66.5,20){\circle*{0.4}}
\put(68.5,20){\circle*{0.4}}
\put(70.5,20){\circle*{0.4}}
%
\put(11,13.5){\line(1,0){2}}
\put(14,13.5){\line(1,0){2}}
\put(17,13.5){\line(1,0){2}}
\put(20,13.5){\line(1,0){2}}
\put(38,13.5){\line(1,0){2}}
\put(41,13.5){\line(1,0){2}}
\put(50,13.5){\line(1,0){2}}
\put(53,13.5){\line(1,0){2}}
\put(56,13.5){\line(1,0){2}}
\put(59,13.5){\line(1,0){2}}
\put(62,13.5){\line(1,0){2}}
%
\put(46,18.1){\line(0,-1){9.2}}
\put(47.9,20){\line(1,0){9.2}}
\put(8.9,7){\line(1,0){9.2}}
\put(34.9,7){\line(1,0){9.2}}
\put(26.5,16.5){\vector(0,-1){6}}
\put(26.5,10.5){\vector(0,1){6}}
\put(23.5,7){\circle*{0.6}}
\put(26.5,7){\circle*{0.6}}
\put(29.5,7){\circle*{0.6}}
\put(23.5,20){\color{lgrey}\circle*{0.8}}
\put(26.5,20){\color{lgrey}\circle*{0.8}}
\put(29.5,20){\color{lgrey}\circle*{0.8}}
\thicklines
\put(8.9,20){\color{lgrey}\line(1,0){9.2}}
\put(34.9,20){\color{lgrey}\line(1,0){9.2}}
\put(7,7){\circle{3}}
\put(46,20){\circle{3}}
\put(7,18.1){\color{hgrey}\line(0,-1){9.2}}
\end{picture}
\caption{A picture illustrating the proof of Theorem~\ref{thm.poset-characterization}.\label{fig.2}}
\end{center}
\end{figure}

Assume for a contradiction that $v\prec_{\permC} u$ giving
$\{u,v\}\in\Inv [\permA,\permC]$ and, by Lemma~\ref{lem.permut-graph}(ii), a unique
edge of $\mathfrak{g}$ is labeled by $\{u,v\}$. Let it be the one between $\permE\in S$ and $\permD\in S$, with
$\permE$ closer to $\permC$. Thus, the section $\mathfrak{g}^{\prime}$ of $\mathfrak{g}$ between
$\permD$ and $\permA$ is in $S_{u\prec v}$. Since the map
$\permG\mapsto\permG\trans_{uv}$ is an automorphism of the permutohedral graph the section $\mathfrak{g}^{\prime}$
has a copy $\mathfrak{g}^{\prime\prime}$ in $S_{v\prec u}$, which is between $\permE$ and $\permB$
(see Figure~\ref{fig.2}(b) for illustration).
The walk composed of $\mathfrak{g}^{\prime\prime}$ and the edge between $\permB$ and $\permA$
has the same length as the section of $\mathfrak{g}$ between $\permE$ and $\permA$ then.
Therefore, it is a geodesic between $\permE$ and $\permA$ and the assumption that $S$
is geodetically convex implies a contradictory conclusion $\permB\in S$.
Thus, one necessarily has $u\prec_{\permC} v$.

{\bf III.}\, The next step it to show that, if $\emptyset\neq S\subseteq\enuN$ is a non-empty geodetically convex set
and $\permE\in\enuN\setminus S$ then there is $(u,v)\in\Cov (S)$ with
$v\prec_{\permE} u$. Consider a walk $\mathfrak{g}$ from $\permE$ to some $\permA\in S$ of
the least possible length among walks from $\permE$ to $S$. By definition, $\mathfrak{g}$ is a geodesic between $\permE$
and $\permA$. Let $\permB$ the last node of $\mathfrak{g}$ before $\permA$. Then necessarily
$\permB\in\enuN\setminus S$; let $\{u,v\}$ be the label of the edge between $\permB$ and $\permA$.
Without loss of generality assume $u\prec_{\permA} v$ (otherwise we exchange $u$ for $v$), which gives $(u,v)\in\Cov(S)$. By Lemma~\ref{lem.permut-graph}(ii) applied to $\mathfrak{g}$, $\{u,v\}$ is an inversion
between $\permE$ and $\permA$, which gives $v\prec_{\permE} u$.

{\bf IV.}\, As $|N|\geq 2$ is assumed, one has $\emptyset\in {\calX}^{\preced}$ and to verify
the sufficiency it is enough to show $S\in {\calX}^{\preced}$
for a non-empty geodetically convex set $\emptyset\neq S\subseteq\enuN$.
By Step II., $S\subseteq\Cov (S)^{\lGprec}$.
By Step~III., $\enuN\setminus S ~\subseteq~ \enuN\setminus (\Cov(S)^{\lGprec})$ gives $\Cov (S)^{\lGprec}\subseteq S$.
Thus, we have shown $S=\Cov (S)^{\lGprec}$, which means $S\in {\calX}^{\preced}$ (see Appendix~\ref{ssec.app.Galois}).
\end{proof}

\subsubsection{Remarks on the semi-lattice of posets}\label{sssec.posets-remarks}
This sub-section contains some advanced observations and can be skipped without losing understanding of the rest of the paper.
Given a poset $T\subseteq N\times N$ on $N$ and the
respective set $S={\cal L}(T)=T^{\lGprec}$ of enumerations (in ${\calX}^{\preced}$), one can show, which is omitted in this paper, that,
for any two-element subset $\{u,v\}\subseteq N$, exclusively one of the following conditions is true:
\begin{itemize}
\item $\{u,v\}\in\Inv (S)$,
\item $u$ and $v$ are comparable in $T$, that is, either $(u,v)\in T$ or $(v,u)\in T$.
\end{itemize}
In other words, $\{u,v\}\in\Inv (S)$ iff $u$ and $v$ are incomparable in $T$.
One can additionally show that, for an ordered pair $(u,v)\in N\times N$, one has $(u,v)\in\Cov(S)$ iff
$u\cover v$ in the poset $T$. These observations justify our terminology introduced before Theorem~\ref{thm.poset-characterization}.

One can show using these facts that, if $n\geq 2$, the function
$S\in {\calX}^{\preced}\setminus \{\emptyset\}\mapsto |\,\Inv(S)\,|$, extended  by a convention
$|\,\Inv(\emptyset)\,|:=-1$, is a {\em height function\/} for the graded lattice $({\calX}^{\preced},\subseteq)$.
The reader may tend to think that this function coincides, for $\emptyset\neq S\in {\cal X}^{\preced}$, with the ``diameter"
function, which assigns the number $\mbox{\rm diam}\,(S) :=\max\,\{\, \dist(\permA,\permB)\,:\ \permA,\permB\in S\,\}$
to every $\emptyset\neq S\subseteq\enuN$.
This is indeed the case if $|N|\leq 5$, but in general
one only has $\mbox{\rm diam}\,(S)\leq |\,\Inv(S)\,|$.

This phenomenon closely relates to the concept of {\em dimension\/} of a poset $T$ \cite[\S\,1.1]{Wie17},
which can equivalently be defined by the formula $\mbox{\rm dim}\,(T) \,:=\, \min\,\{\, |S|\, :\ T=S^{\rGprec}\,\}$.
Note in this context that one can show that, for every pair of distinct enumerations $\permA,\permB\in\enuN$,
the set of enumerations between $\permA$ and $\permB$ belongs to ${\calX}^{\preced}$ and corresponds to a
poset of dimension 2, defined as the intersection of $T_{\permA}$ and $T_{\permB}$.
A classic result \cite{Hir51} in the theory of poset dimension says that $\mbox{\rm dim}\,(T) \leq
\lfloor \frac{n}{2}\rfloor$ for a poset $T$ on an $n$-element set $N$. 
This upper bound is tight owing to a classic construction \cite{DM41} of a poset of given dimension.
Thus, the simplest example of a poset of dimension 3 exists in case $|N|=6$. It also gives
the simplest example of a poset-based set $S\subseteq\enuN$ of enumerations for which
$\mbox{\rm diam}\,(S)< |\,\Inv(S)\,|$.

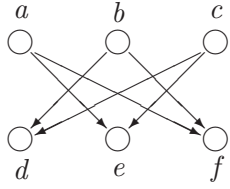
\begin{figure}[t]
\setlength{\unitlength}{1mm}
\begin{center}
\begin{picture}(40,24)
\put(7,5){\circle{3}}
\put(7,18){\circle{3}}
\put(20,5){\circle{3}}
\put(20,18){\circle{3}}
\put(33,5){\circle{3}}
\put(33,18){\circle{3}}
\put(7,22){\makebox(0,0){\it a}}
\put(20,22){\makebox(0,0){\it b}}
\put(33,22){\makebox(0,0){\it c}}
\put(7,1){\makebox(0,0){\it d}}
\put(20,1){\makebox(0,0){\it e}}
\put(33,1){\makebox(0,0){\it f}}
\put(8.5,16.5){\vector(1,-1){10}}
\put(31.5,16.5){\vector(-1,-1){10}}
\put(21.5,16.5){\vector(1,-1){10}}
\put(18.5,16.5){\vector(-1,-1){10}}
\put(9,16.5){\vector(2,-1){22}}
\put(31,16.5){\vector(-2,-1){22}}
\end{picture}
\end{center}
\caption{Directed Hasse diagram of the poset from Example~\ref{exa.dim-poset}.\label{fig.3}}
\end{figure}

\begin{exa}\label{exa.dim-poset}\rm
Put $N:=\{a,b,c,d,e,f\}$ and $T:=\diag\cup\{\, (a,e),(a,f),(b,d),(b,f),(c,d),(c,e) \,\}$;
see Figure~\ref{fig.3} for the directed Hasse diagram of this poset.
One can show that the set $S:=T^{\lGprec}$ is the (disjoint) union of (face-associated) subsets of $\enuN$:
$$
|\,abc\,|\,def\,|, \quad |\,ab\,|\,f\,|\,c\,|\,de\,|,  \quad |\,ac\,|\,e\,|\,b\,|\,d\,f\,|,
\quad |\,bc\,|\,d\,|\,a\,|\,ef\,|\,.
$$
Hence, the inversions in $S$ are just incomparable pairs in $T$, that is, $|\,\Inv(S)\,|=9$.
On the other hand, for $\permA\in S$, if $\permA\not\in |\,ab\,|\,f\,|\,c\,|\,de\,|$ then $c\prec_{\permA} f$. This observation allows one to conclude, using Lemma~\ref{lem.permut-graph}(i),
that any geodesic between elements of $S\setminus\, |\,ab\,|\,f\,|\,c\,|\,de\,|$ has the length at most $8$.
The implication $\permA\in S\setminus\, |\,ac\,|\,e\,|\,b\,|\,d\,f\,| ~\Rightarrow~
b\prec_{\permA} e$ and the implication $\permA\in S\setminus\, |\,bc\,|\,d\,|\,a\,|\,ef\,| ~\Rightarrow~ a\prec_{\permA} d$ lead to analogous conclusions. This, combined with the choice $\permA := |a|c|e|b|f|d|$ and $\permB := |b|c|d|a|f|e|$,
for which $\dist (\permA,\permB)=|\,\Inv [\permA,\permB ]\,|=8$
allows one to observe that $\mbox{\rm diam}\,(S)=8$.
\end{exa}

The height function for $({\calX}^{\preced},\subseteq)$ was described above
in combinatorial terms, namely by means of inversions between enumerations.
On the other hand, elements of ${\calX}^{\preced}$ were characterized in Theorem~\ref{thm.poset-characterization}
solely in graphical terms. The reader may ask whether the height function can also be defined in this way.
This is indeed the case. The point is that the {\em combinatorial equivalence of edges\/} in the permutohedral graph,
being defined by the condition that the edges are labeled by the same two-element subset $\{u,v\}$ of $N$,
can equivalently be defined solely in graphical terms.
Specifically, given an edge between $\permA\in\enuN$ and $\permB\in\enuN$, consider the set
$$
S_{\permA :\permB}:=\{\, \permC\in\enuN\,:\ \dist(\permC,\permA)<\dist(\permC,\permB)\,\}
$$
of enumerations that are closer to $\permA$ than to $\permB$.
If the edge is labeled by
$\{u,v\}\subseteq N$ and $u\prec_{\permA} v$ then one can observe, using Lemma~\ref{lem.permut-graph}(ii) and the fact that
the mapping $\permG\mapsto \permG\trans_{uv}$ is an automorphism of the permutohedral graph,
that $S_{u\prec v}\subseteq S_{\permA :\permB}$. Since $S_{v\prec u}$ is the complement of $S_{u\prec v}$ in $\enuN$,
one has $S_{u\prec v}=S_{\permA :\permB}$ and $S_{v\prec u}= S_{\permB :\permA}$, meaning that
$\{S_{v\prec u},S_{u\prec v}\}=\{S_{\permA :\permB},S_{\permB :\permA}\}$ is a partition of $\enuN$, definable in graphical terms.
In particular, two edges have the same combinatorial label $\{u,v\}\subseteq N$ iff they
yield the same graphical bi-partition of $\enuN$.

This particular equivalence of edges also has geometric interpretation that the respective geometric edges of the permutohedron $\perN$ (= segments connecting the respective rank vectors) are {\em parallel}, or formulated
in other way, are perpendicular to the same hyperplane
$$
{\mathsf H}_{uv}:=\{x\in {\dv R}^{N}\, :\ x_{u}=x_{v}\,\}\,.
$$
To summarize: the edge-labeling introduced in Section~\ref{sssec.poset-lattice-characterization}
has deeper geometric meaning and the value of the height function for $({\calX}^{\preced},\subseteq)$ at non-empty
$S\in {\calX}^{\preced}$ is nothing but the number of different labels/colors of edges within $S$.
\smallskip

Our proof of Theorem~\ref{thm.poset-characterization} was inspired
by a result from \cite[Theorem\,9]{HN13}, where another characterization of poset-based subsets
of $\enuN$ was presented. It can be re-phrased as follows: (a non-empty connected) set $S\subseteq\enuN$
corresponds to a poset on $N$ iff, for every $(u,v)\in N\times N$ with $u\neq v$, exclusively one of the
following 3 conditions holds:
\begin{itemize}
\item $\{u,v\}\in\Inv (S)$,
\item $(u,v)$ belongs to the transitive closure of $\Cov (S)$,
\item $(v,u)$ belongs to the transitive closure of $\Cov (S)$.
\end{itemize}
Note that a tacit assumption that $S\subseteq\enuN$ is connected in the permutohedral
graph was omitted in the formulation of \cite[Theorem\,9]{HN13}, probably by a mistake.
\smallskip

The last remark is that if $|N|\geq 3$ then the lattice
based on posets is not a face lattice (of any polytope). This
follows from \cite[Theorem\,2.7(iii)]{Zie95}. Specifically, in case of $N=\{a,b,c\}$ the lattice has an interval of the length 2 which has only 3 elements instead of 4 (thus, it is not a ``diamond"):
take $\{\, |a|b|c| \,\}$, $\{\, |a|b|c| \,-\!\!\!\!-\,\, |b|a|c| \,\}$, and $\{\, |a|b|c|
\,-\!\!\!\!-\,\, |b|a|c| \,-\!\!\!\!-\,\, |b|c|a| \,\}$.

\subsection{Lattice of preposets viewed from different perspectives}\label{ssec.preposet-lattice}
Another important lattice we deal with in this paper is the lattice of {\em preposets\/} on $N$.
There are two possible interpretations of these preposets. One of them
is geometric, in terms of special cones in ${\dv R}^{N}$, as presented in \cite[\S\,3.4]{PRW08}.
The other interpretation is combinatorial, in terms of certain rings of subsets of $N$, and this
relates to well-known Birkhoff's representation theorem for finite distributive lattices \cite[\S\,III.3]{Bir95}.
Both interpretations can be introduced in terms Galois connections, as described in
Appendix~\ref{ssec.app.Galois}.

\subsubsection{Geometric view: braid cones}\label{sssec.preposet-cones}
To this end we define a binary relation $\braid$ between vectors in $X:={\dv R}^{N}$
and the elements of the Cartesian product $Y:=N\times N$. Specifically, given a vector $x\in {\dv R}^{N}$
and $(u,v)\in N\times N$, we define their incidence relation $\braid$ through the {\em comparison of respective vector components\/}:
$$
x\,\braid\, (u,v) \quad :=\quad x_{u}\leq x_{v}\,.
$$
We will use larger black triangles $\rGbraid$ and $\lGbraid$ to
denote Galois connections based on this relation~$\braid$.
Any subset $T$ of\/ $Y=N\times N$ is assigned (by backward Galois connection) a polyhedral cone
$$
\nor_{T}:= T^{\lGbraid}= \{\, x\in {\dv R}^{N}\,:\ \forall\,(u,v)\in T~~  x_{u}\leq x_{v}\,\}\,,
$$
called the {\em braid cone\/} for $T$. In particular, any enumeration $\permA\in\enuN$ of $N$ is assigned a particular full-dimensional cone in ${\dv R}^{N}$ as the braid cone for the toset $T_{\permA}$\,:
\begin{eqnarray}
\lefteqn{\hspace*{-23mm}\nor^{\permA} \,:=~ \{\, x\in {\dv R}^{N}\,:\ x_{\permA(1)} \leq x_{\permA(2)} \leq \ldots \leq x_{\permA(n)}\,\,\}\,} \label{eq.enum-cone}\\
 &=& \{\, x\in {\dv R}^{N}\,:\ \forall\, (u,v)\in N\times N \quad u\preceq_{\permA}v ~\Rightarrow~  x_{u}\leq x_{v}\,\,\}; \nonumber
\end{eqnarray}
these cones are called {\em Weyl chambers\/} by some authors \cite{Mor07,PRW08}.
One particular geometric relation between
braid cones in ${\dv R}^{N}$ plays an important role in our context:
the relation of a cone of being a {\em face} of another cone.
This relation appears to have a counterpart in the world of preorders:
given preposets\/ $R,T\subseteq N\times N$, we say that
$R$ is a {\em contraction\/} of\/ $T$ and write
\begin{eqnarray*}
\lefteqn{\hspace*{-18mm}T\contr R ~~\mbox{if}~~ [\,\exists\, Q\subseteq T^{op} \,:\ R=\tr(T\cup Q)\,],~~
\mbox{with $\tr(T)$ denoting the transitive closure of\/ $T$}}\\
&& \hspace*{47mm}\mbox{and $T^{op}$ the opposite (preposet) of\/ $T$.}
\end{eqnarray*}
The next lemma claims that the method of Galois connections yields the lattice of preposets.

\begin{lemma}\label{lem.braid-lattice}\rm
Let $({\calX}^{\braid},\subseteq )$ and $({\calY}^{\braid},\subseteq )$ denote the lattices defined using the
Galois connections $\rGbraid$ and $\lGbraid$ based on the incidence relation $\braid$ introduced above.
\begin{itemize}
\item[(i)] Then $T\in {\calY}^{\braid}$ iff\/ $T$ is a preposet on $N$.
\item[(ii)] The braid fan (see Appendix~\ref{ssec.app.normal.fan}) is a subset
of the lattice ${\calX}^{\braid}$ of braid cones.
\item[(iii)] Given a preposet $T\subseteq N\times N$, its braid cone $\nor_{T}$ is full-dimensional
iff $T$ is a poset.
\item[(iv)] Given preposets $R,T\subseteq N\times N$, the cone
$\nor_{R}$ is a face of $\nor_{T}$ iff $R$ is a contraction of\/ $T$.
\item[(v)] The lattice $({\calY}^{\braid},\subseteq )$ is
atomistic and coatomistic.
\end{itemize}
\end{lemma}

Thus, $({\calY}^{\braid},\subseteq )$ is the lattice of preposets on $N$. The closure operation corresponding to this Moore family is the reflexive-and-transitive closure: given $T\subseteq N\times N$,
$$
T^{\lGbraid\rGbraid} = \tr(T\cup\diag ),\quad
\mbox{where $\tr(T)$ denotes the transitive closure of\/ $T$.}
$$
Note that, if $n\geq 3$, then the lattice of preposets is not graded and there are braid cones other than
elements of the braid fan or the whole space ${\dv R}^{N}$.

\begin{proof}
For the necessity in (i) assume $T=\nor^{\rGbraid}=\{\,(u,v)\in N\times N\, :\ \forall\, x\in\nor\quad x_{u}\leq x_{v} \,\}$ for some $\nor\subseteq {\dv R}^{N}$ (see Appendix~\ref{ssec.app.Galois}). Evidently, $T$ is both reflexive and transitive.

For the sufficiency in (i) assume that $T$ is a preposet on $N$ and denote by ${\cal E}$ the set of equivalence classes of $E:=T\cap T^{op}$. Put $\nor:=\nor_{T}=T^{\lGbraid}$; it is immediate
that\/ $T\subseteq \nor^{\rGbraid}$. We will show
$(a,b)\in (N\times N)\setminus T ~\Rightarrow~ (a,b)\not\in \nor^{\rGbraid}$ to verify the other inclusion $\nor^{\rGbraid}\subseteq T$.

Indeed, the interpretation of $T$ in the form of a poset $\preceq$ on ${\cal E}$ (see Appendix~\ref{sec.app.lattice}),
which can additionally be viewed as a transitive directed acyclic graph $G$ over ${\cal E}$ (Appendix~\ref{ssec.app.trans-DAG}), implies that there are $A,B\in {\cal E}$ with $a\in A$, $b\in B$ and
$\neg[A\preceq B$], that is, $A$ is not an ancestor of $B$ in $G$. Therefore, an ordering $U_{1},\ldots, U_{m}$, $m\geq 1$, of elements of ${\cal E}$ exists, which is consonant with $G$
and in which $B$ strictly precedes $A$.
Define $x\in {\dv R}^{N}$ as the ``rank vector" of this block-enumeration:
put, for any $u\in N$, $x_{u}:=i$ where $u\in U_{i}$. Then $(u,v)\in T$ implies either $[\exists\,i\,: u,v\in U_{i}]$ or $[\exists\,i<j\,: u\in U_{i}\,\&\,v\in U_{j}]$, yielding $x_{u}\leq x_{v}$.
Hence, $x\in \nor_{T}=\nor$ but $x_{b}<x_{a}$, and $(a,b)\not\in \nor^{\rGbraid}$.

The claim in (ii) follows easily from the comparison of the definition of $\nor_{T}$ with the definition of the braid fan from Appendix~\ref{ssec.app.normal.fan}.

As concerns (iii), if $\dim(\nor_{T})=n$ then no pair of distinct $u,v\in N$ with $(u,v),(v,u)\in T$ exists, for otherwise $\nor_{T}$ is contained in a proper linear subspace of\/ ${\dv R}^{N}$. Hence, $T\cap T^{op}=\diag$
and\/ $T$ is anti-symmetric.
Conversely, if $T$ is poset then ${\cal L}(T)\neq\emptyset$ and $T\subseteq T_{\permA}$ for some enumeration $\permA\in\enuN$ and the anti-isomorphism of $({\calX}^{\braid},\subseteq )$ and $({\calY}^{\braid},\subseteq )$ yields $\nor^{\permA}\subseteq  \nor_{T}$.

For (iv) realize that a (non-empty) face of $\nor_{T}$ can be obtained by turning some of its facet-defining inequalities into equalities (Appendix~\ref{ssec.app.vertices-facets}),
which is equivalent to adding the respective converse inequalities.
In the language of preposets this corresponds to adding a part of $T^{op}$.

As concerns (v), realize that $\diag\cup\{(u,v)\}$, for distinct $u,v\in N$, are posets on $N$,
and, thus atoms in $({\calY}^{\braid},\subseteq)$. Clearly, any preposet on $N$ is the union of those, and, thus,
their supremum.
Given a partition $U,V\subseteq N$ into two non-empty blocks,
the set $(N\times N)\setminus (V\times U)$ is a preposet on $N$, and, thus, a coatom in $({\calY}^{\braid},\subseteq)$. Given a preposet $T$ on $N$, to show that $T$ is the intersection
of those, we fix $(a,b)\in (N\times N)\setminus T$ and repeat the consideration in the proof of (i).
Take the respective enumeration $U_{1},\ldots, U_{m}$, $m\geq 1$, of ${\cal E}$ in which
$B=:U_{\ell}$ strictly precedes $A$. Put $U:=\bigcup_{i\leq\ell} U_{i}$, $V:=\bigcup_{i>\ell} U_{i}$,
and observe $T\subseteq (N\times N)\setminus (V\times U)$ and $(a,b)\in V\times U$.
\end{proof}

\subsubsection{Combinatorial view: finite topologies}\label{sssec.preposet-topologies}

Define a binary relation $\topol$ between elements of the power set $X:=\caP$
and elements of $Y:=N\times N$. Specifically, if $D$ is a subset of $N$
and $(u,v)\in N\times N$ then we consider $D$ to be in the incidence relation $\topol$ with the pair $(u,v)$
if the set {\em $D$ respects coming $u$ before $v$\,}:
$$
D\,\topol\, (u,v) \quad := \quad [\,u\in D ~\vee~ v\not\in D\,]\,,
~~ \mbox{interpreted as $[\,v\in D \,\Rightarrow\, u\in D\,]$.}
$$
We will use larger gray triangles $\rGtopol$ and $\lGtopol$ to
denote Galois connections based on this relation~$\topol$.
Any subset $T$ of\/ $Y=N\times N$ is assigned (by backward Galois connection) a set system
$$
{\calD}_{T}:= T^{\lGtopol}= \{\, D\subseteq N\,:\ \forall\,(u,v)\in T~~  \mbox{if $v\in D$ then $u\in D$}\,\}~~\subseteq \caP\,,
$$
known as the class of {\em down-sets\/} for\/ $T$.
Note that any such a class $\calD_{T}$ satisfies $\emptyset,N\in\calD_{T}$ and is closed under set intersection and union: $D,E\in\calD_{T} \,\Rightarrow\, D\cap E,D\cup E\in\calD_{T}$.
Thus, it forms a ring of sets and can be interpreted as a (finite) {\em topology on $N$}, because in case of finite $N$, definition of a topology on $N$ reduces to the above requirements.

Any enumeration $\permA\,=\, |\,\permA(1)|\ldots|\permA(n)\,|$ of $N$ can be represented by a particular
topology $\chain_{\permA}$ on $N$, the one assigned to the toset $T_{\permA}$, namely the respective {\em  maximal chain\/} of subsets of $N$:
\begin{equation}
\chain_{\permA} ~:=~ \left\{~ \bigcup_{i=1}^{j} \permA (i)\quad :\ j=0,\ldots ,n~\right\}
~~=~~ \left\{~ \emptyset ,\, \{\permA(1)\}\,,\, \{\permA(1),\permA(2)\}\,,~ \ldots ~,\,
{N} ~\right\}\,.
\label{eq.max-chain}
\end{equation}

We show that $({\calY}^{\topol},\subseteq )$ coincides with the lattice of preposets $({\calY}^{\braid},\subseteq )$, which establishes implicitly a one-to-one correspondence between braid cones in ${\dv R}^{N}$ (= elements
of ${\calX}^{\braid}$) and topologies on $N$ (= elements of ${\calX}^{\topol}$).
The relation between topologies reflecting the contraction
relation between preposets (or corresponding to being-a-face relation between braid cones)
has slightly awkward description in terms of topologies.

\begin{lemma}\label{lem.topol-lattice}\rm
Let $({\calX}^{\topol},\subseteq )$ and $({\calY}^{\topol},\subseteq )$ denote the lattices defined using the
Galois connections $\rGtopol$ and $\lGtopol$ based on the incidence relation $\topol$ defined above.
\begin{itemize}
\item[(i)] Then $T\in {\calY}^{\topol}$ iff\/ $T$ is a preposet on $N$.
\item[(ii)] A set system $\calD\subseteq\caP$ belongs to the lattice ${\calX}^{\topol}$ iff
both $\emptyset,N\in\calD$ and the implication $D,E\in\calD \,\Rightarrow\, D\cap E,D\cup E\in\calD$
holds (= $\calD$ is a topology on $N$).
\item[(iii)] A preposet $T\subseteq N\times N$ is a poset iff
the topology $\calD_{T}$ distinguishes points: for any distinct $u,v\in N$
there is $D\in\calD_{T}$ such that either $[\,u\in D ~\&~ v\not\in D \,]$ or
$[\,v\in D ~\&~ u\not\in D \,]$. This happens exactly when there is $\permA\in\enuN$ such that
$\chain_{\permA}\subseteq\calD_{T}$.
\item[(iv)] Given a braid cone $\nor\subseteq {\dv R}^{N}$, the set system
$\calD_{\nor} :=\{\, D\subseteq N\, :\ -\chi_{D}\in\nor\,\}$ is the corresponding topology.
Given a topology $\calD$, the cone
$\nor_{\calD} := \cone(\{\chi_{N}\}\cup\{-\chi_{D}:D\in\calD\})$ is the corresponding braid cone.
\item[(v)] Given preposets $R,T\subseteq N\times N$, one has $T\contr R$   iff
the topology $\calD_{R}$ is a {\em reduction\/} of the topology  $\calD_{T}$ in the following sense:
for some (possibly empty) set of pairs $(u_{i},v_{i})\in N\times N$ with
$u_{i}\in\bigcap_{v_{i}\in D\in\calD_{T}} D$ one has
$\calD_{R}=\{\, D\in\calD_{T} \,:\ u_{i}\in D \,\Rightarrow\, v_{i}\in D~\mbox{for all $i$}\,\}$.
\end{itemize}
\end{lemma}

Note in this context that, by the fundamental representation theorem for finite distributive lattices, rings of subsets of a finite set (ordered by inclusion) are universal examples of such lattices.
What plays substantial role in the proof of this result is a particular correspondence between finite distributive lattices and finite posets, mediated by the concept of a down-set;
see \cite[Corollary 1 in \S\,III.3]{Bir95} or \cite[Theorem 1.4 in \S\,1.3.2]{Gra16}.
This is the correspondence we discuss here, but with one important technical difference. While
distinguishing points is immaterial for the purpose of representing a (distributive) lattice,
it matters for our purposes. We do need to extend the correspondence to preposets and topologies on the \underline{same} finite set.

\begin{proof}
For the necessity in (i) assume\, $T=\calD^{\rGtopol}:=\{\,(u,v)\in N\times N\, :\ \forall\, D\in\calD~~\,
v\in D \Rightarrow u\in D \,\}$ for some $\calD\subseteq\caP$ (see Appendix~\ref{ssec.app.Galois}). Evidently, $T$ is both reflexive and transitive.

For the sufficiency in (i) assume that $T$ is a preposet on $N$; it is immediate that $T\subseteq (\calD_{T})^{\rGtopol}$. For any $v\in N$ we put $I(v):=\{ w\in N\,:\ (w,v)\in T\,\}$; the reflexivity of\/ $T$ implies $v\in I(v)$, while transitivity of\/ $T$ allows one to observe $I(v)\in \calD_{T}$.
Given $(u,v)\in (\calD_{T})^{\rGtopol}$, these two facts about $I(v)$ and the definition of $\calD^{\rGtopol}$ imply $u\in I(v)$, which means $(u,v)\in T$. Therefore, $(\calD_{T})^{\rGtopol}\subseteq T$ and $T=(\calD_{T})^{\rGtopol}$ yields $T\in {\calY}^{\topol}$.

The necessity in (ii) is immediate: assume $T\subseteq N\times N$ and realize that the definition of\/ $\calD_{T}$ allows one easily to verify the conditions in (ii) for $\calD_{T}$.

For the sufficiency in (ii) assume that $\calD\subseteq\caP$ is a topology on $N$ and put $T:=\calD^{\rGtopol}$; it is immediate that $\calD\subseteq T^{\lGtopol}=\calD_{T}$. For any $v\in N$ we put $D(v):=\bigcap_{\,v\in D\in\calD} D$.
Then $N\in\calD$ implies $v\in D(v)$ and, since $\calD$ is closed under intersection, $D(v)\in\calD$.
Thus, $E\subseteq\bigcup_{v\in E} D(v)$ for any $E\subseteq N$. Given $E\in\calD_{T}$, the goal is to verify
the other inclusion $\bigcup_{v\in E} D(v)\subseteq E$.\\[0.4ex]
Indeed, having $w\in \bigcup_{v\in E} D(v)$, there exists $v\in E$ with $w\in D(v)$. The latter condition means
$[\,\forall\, D\in\calD~~ v\in D \Rightarrow w\in D\,]$, which is nothing but $(w,v)\in \calD^{\rGtopol}=T$.
By the definition of $\calD_{T}$, the facts $E\in\calD_{T}$, $(w,v)\in T$, and $v\in E$ imply
$w\in E$, concluding the proof of the inclusion.\\[0.4ex]
Therefore, since $\emptyset\in\calD$ and $\calD$ is closed under union, one gets $E=\bigcup_{v\in E} D(v)\in\calD$.
This concludes the proof of $E\in \calD_{T}\,\Rightarrow\, E\in\calD$ and $\calD=\calD_{T}=T^{\lGtopol}$ yields $\calD\in {\calX}^{\topol}$.

As concerns (iii), if $\calD_{T}$ distinguishes points then there is no pair of distinct $u,v\in N$ with $(u,v),(v,u)\in T$, for otherwise $\forall\, D\in\calD_{T}~\, u\in D \Leftrightarrow v\in D$.
Thus, $T\cap T^{op}=\diag$ and\/ $T$ is anti-symmetric. Conversely, if $T$ is poset then $T\subseteq T_{\permA}$ for some enumeration $\permA\in\enuN$ and the anti-isomorphism of $({\calX}^{\topol},\subseteq )$ and $({\calY}^{\topol},\subseteq )$ yields $\chain_{\permA}\subseteq\calD_{T}$. Hence, $\calD_{T}$ distinguishes points.

For (iv) assume that $\nor=\nor_{T}$ and $\calD=\calD_{T}$ for a preposet $T$ on $N$. Given $D\subseteq N$ and $(u,v)\in T$, re-write $D\,\topol\, (u,v)$, that is, $[\,v\in D \Rightarrow u\in D\,]$ in the form
$[\, \chi_{D}(v)=1 \Rightarrow \chi_{D}(u)=1\,]$ and then as $-\chi_{D}(u)\leq -\chi_{D}(v)$, which means $(-\chi_{D})\braid (u,v)$.
This allows one to observe that $D\in\calD_{T}$ iff $-\chi_{D}\in\nor_{T}$, which yields the first claim in (iv).
Since $\nor_{T}$ is a cone, it also implies that $\cone(\{\chi_{N}\}\cup\{-\chi_{D}:D\in\calD_{T}\})\subseteq \nor_{T}$ in the second claim of (iv). Thus, to prove the second claim in (iv) it suffices to verify the other inclusion $\nor_{T}\subseteq \cone(\{\chi_{N}\}\cup\{-\chi_{D}:D\in\calD_{T}\})$.\\[0.3ex]
For this purpose, consider the set ${\cal E}$ of equivalence classes of $E:=T\cap T^{op}$ and view $T$ as a poset on ${\cal E}$ and, thus, as a transitive directed acyclic graph $G$ over ${\cal E}$
(see Appendices~\ref{sec.app.lattice} and \ref{ssec.app.trans-DAG}).
Observe that any vector $x\in\nor_{T}$ complies with $G$ in the following sense:
if $u\in U\in {\cal E}$, $v\in V\in {\cal E}$ and $U$ is an ancestor of $V$ in $G$ then $x_{u}\leq x_{v}$. Hence, $x$ is ``constant" on sets from ${\cal E}$\,: $x_{u}=x_{v}$ if $(u,v)\in E$.
This allows one to construct an ordering $U(1),\ldots, U(m)$, $m\geq 1$, of (all) elements of ${\cal E}$, which both respects $x$ in the sense that
$$
u\in U(i),~ v\in U(j),~ \mbox{and}~ i\leq j\quad \mbox{implies}\quad x_{u}\leq x_{v}
\,,
$$
and is consonant with $G$. Let us introduce $x(i)$ for $i=1,\ldots, m$ as the shared value $x_{u}$ for $u\in U(i)$ and write
$$
x= x(m)\cdot\chi_{N} \,+\, \sum_{i=1}^{m-1}\, (x(i+1)- x(i))\cdot (-\chi_{\,U(1)\ldots U(i)})\,,
$$
which equality can be verified by substituting any $u\in U(j)$ for $j\in [m]$.
Since $x(i+1)-x(i)\geq 0$ for $i\in[m-1]$, the vector $x$ belongs to the conic hull of $\{\,-\chi_{\,U(1)\ldots U(i)}\,:\ i\in [m-1]\,\}\cup\{\pm\chi_{N}\}$.
The consonancy of $U(1),\ldots, U(m)$ with $G$ implies that $U(1)\ldots U(i)\in\calD_{T}$ for $i\in [m]$, which thus yields the desired conclusion.

The reduction condition in (v) is basically re-writing of the definition of the contraction
$T\contr R$\, in terms of $\calD_{T}$ and $\calD_{R}$. Indeed, one has $(u,v)\in T=(\calD_{T})^{\rGtopol}$ iff $u\in\bigcap_{v\in D\in\calD_{T}} D$. Then $Q\subseteq T^{op}$ is a set of opposite pairs $(v_{i},u_{i})$
to $(u_{i},v_{i})$ with $u_{i}\in\bigcap_{v_{i}\in D\in\calD_{T}} D$. Then
\begin{eqnarray*}
\lefteqn{\hspace*{-3mm}\calD_{R}=\calD_{T\cup Q} = \{\, D\subseteq N\, :\ \forall\, (a,b)\in T\cup Q\quad
b\in D \Rightarrow a\in D\,\}}\\
&=& \{\, D\in\calD_{T}\, :\ \forall\, (a,b)\in Q\quad
b\in D \Rightarrow a\in D\,\} = \{\, D\in\calD_{T}\, :\ \forall\, i\quad
u_{i}\in D \Rightarrow v_{i}\in D\,\}\,,
\end{eqnarray*}
which what is written in (v).
\end{proof}

\subsection{A lattice of set systems composed of maximal chains}\label{ssec.weak-lattice}
Another (larger) lattice of subsets of\/ $\enuN$ is isomorphic to a certain class of set systems, ordered by inclusion. We also introduce this lattice by means of Galois connections (see Appendix~\ref{ssec.app.Galois}).
We formally define a binary relation $\outside$ between elements of the set $X:=\enuN$ of enumerations of $N$ and elements of the power set $Y:=\caP$ of $N$: given $\permA\in\enuN$ and $L\subseteq N$ we put
$$
\permA\,\outside\, L \quad :=\quad L\not\in \chain_{\permA}\,,
$$
which means that the set $L$ is outside the maximal chain of sets given by $\permA$.
Let $({\calX}^{\outside},\subseteq )$ denote the lattice of subsets of\/
$\enuN$ defined through the Galois connections based on the incidence relation \,$\outside$\, defined above.
This lattice is anti-isomorphic to its dual lattice ${\calY}^{\outside}$ of subsets of $Y=\caP$ and it is more convenient
to identify any $S\in {\calX}^{\outside}$ isomorphically with the complement ${\calD}:=\caP\setminus{\cal T}$
of the forward Galois connection image ${\cal T}=S^{\rGalo}\in {\calY}^{\outside}$. This means
$$
S\longmapsto {\calD}=\bigcup_{\permA\in S} \chain_{\permA} ~~\mbox{and conversely}~~
{\calD}\longmapsto S=\{\,\permA\in\enuN\, :\ \chain_{\permA}\subseteq {\calD}\,\}
$$
and $({\calX}^{\outside},\subseteq )$ is order-isomorphic with the collection of
set systems ${\calD}\subseteq\caP$, composed of maximal chains in $\caP$, ordered by inclusion. Note in this context
that the lattice $({\calX}^{\outside},\subseteq )$ can be shown to be both atomistic and coatomistic but it
is not a graded lattice for $n=3$.
The lattice $({\calX}^{\preced},\subseteq )$  based on posets from Section~\ref{ssec.poset-lattice} is naturally embedded into $({\calX}^{\outside},\subseteq )$.

\begin{lemma}\label{lem.chain-lattice}\rm
Let $T\subseteq N\times N$ be a poset on $N$. Then
\begin{itemize}
\item[(i)] $\bigcup_{\permA\in {\cal L}(T)} \chain_{\permA}=\calD_{T}$, and
\item[(ii)] ${\cal L}(T)=\{\,\permA\in\enuN\, :\ \chain_{\permA}\subseteq \calD_{T}\,\}$.
\end{itemize}
In particular, any $S\in{\calX}^{\preced}$ from the lattice defined in Section~\ref{ssec.poset-lattice} belongs to ${\calX}^{\outside}$.
\end{lemma}

Note in this context that the assignment $T\mapsto \calD_{T}$ from Section~\ref{sssec.preposet-topologies}
restricted to posets is one of standard constructions in combinatorics \cite[\S\,3.5]{Sta97}.
Specifically, Lemma~\ref{lem.chain-lattice}(ii) allows one to identify {\em linear extensions\/} for\/ $T$
on basis of $\calD_{T}$ and, thus, this two-step procedure forms the basis for common algorithms to compute the number of linear extensions for $T$.

\begin{proof}
To verify the inclusion $\bigcup_{\permA\in {\cal L}(T)} \chain_{\permA}\subseteq\calD_{T}$ in (i) consider a
set of the form $D=\bigcup_{i=1}^{j} \permA(i)$ for some $\permA\in {\cal L}(T)$. Given
$v=\permA(i)\in D$, $1\leq i\leq j$, and $(u,v)\in T$ one has $u\preceq_{\permA} v$, that is,
$\permA_{-1}(u)\leq\permA_{-1}(v)=i\leq j$ implying $u\in D$. Hence, by definition, $D\in \calD_{T}$.

To verify the inclusion $\calD_{T}\subseteq \bigcup_{\permA\in {\cal L}(T)} \chain_{\permA}$ in (i) we view $T\setminus\diag$ as a
directed acyclic graph $G$ over $N$ (see Appendix~\ref{ssec.app.trans-DAG}),
when $\calD_{T}$ is the system of ancestral sets in $G$. Given an ancestral set $A\subseteq N$ in $G$, an enumeration
of $A$ consonant with $G_{A}$ exists, which can be extended to an enumeration $\permA$ of $N$
consonant with $G$. Clearly, $D=\bigcup_{i=1}^{|A|} \permA(i)$, which verifies $D\in \bigcup_{\permA\in {\cal L}(T)} \chain_{\permA}$.

The implication $\permA\in {\cal L}(T) ~\Rightarrow~ \chain_{\permA}\subseteq \calD_{T}$ in (ii) is immediate from (i). To verify the converse implication consider $\permA\in\enuN$ with $\chain_{\permA}\subseteq \calD_{T}$. Given $(u,v)\in T$, consider $1\leq j\leq n$
with $v=\permA(j)$ and have
$D:=\bigcup_{i=1}^{j} \permA(i)\in \chain_{\permA}\subseteq \calD_{T}$.
The definition of $\calD_{T}$ and $v\in D$ imply $u\in D$.
As $\permA_{-1}(u)\leq j=\permA_{-1}(v)$ gives
$u\preceq_{\permA} v$, completing the proof of $\permA\in {\cal L}(T)$.

The last claim is trivial for empty $S$ (if $n\geq 2$).
Given $\emptyset\neq S\in {\calX}^{\preced}$, by Lemma~\ref{lem.poset-lattice}(ii), $S={\cal L}(T)$
for some poset $T\subseteq N\times N$ on $N$ and using (ii) one has $S\in {\calX}^{\outside}$.
\end{proof}

While Lemma~\ref{lem.chain-lattice} describes the embedding of $({\calX}^{\preced},\subseteq )$
in terms of subsets of $\enuN$, one can also describe the situation in terms of subsets of $\caP$.

\begin{corol}\label{cor.set-systems}\rm
Let $\calD\subseteq\caP$ be a set system.
\begin{itemize}
\item[(i)] Then $\calD$ is a topology distinguishing points
iff $\calD\neq\emptyset$ and $S\in {\calX}^{\preced}$ exists with ${\calD}=\bigcup_{\permA\in S}\chain_{\permA}$, which
is then non-empty and uniquely determined as $S=\{\,\permA\in\enuN\, :\ \chain_{\permA}\subseteq {\calD}\,\}$.
\item[(ii)]If $\emptyset\neq S\in{\calX}^{\outside}$ is
such that $\bigcup_{\permA\in S}\chain_{\permA}$ is closed under intersection and union then $S\in{\calX}^{\preced}$.
\end{itemize}
\end{corol}

Note that ${\calX}^{\outside}\setminus\{\emptyset\}$ differs from ${\calX}^{\preced}\setminus\{\emptyset\}$ if $n\geq 3$:
take $N:=\{a,b,c\}$ and $S:=\{\, |a|b|c|, |c|b|a| \,\}$.
Then ${\calD}:=\bigcup_{\permA\in S}\chain_{\permA}=N\setminus\{\,b, ac\,\}$
and $S=\{\,\permA\in\enuN\, :\ \chain_{\permA}\subseteq {\calD}\,\}$
ensures $S\in{\calX}^{\outside}$.
However, $\calD$ is not a topology (although it distinguishes points) for which reason $S\not\in {\calX}^{\preced}$.
%

\begin{proof}
As concerns the necessity in (i), using Lemma~\ref{lem.topol-lattice}(ii)(iii),
observe that $\calD=\calD_{T}$ for some poset $T$ on $N$
and the set $S:=\{\,\permA\in\enuN\, :\ \chain_{\permA}\subseteq {\calD}\,\}$ is non-empty. Lemma~\ref{lem.chain-lattice}(ii)
then gives ${\cal L}(T)=\{\,\permA\in\enuN\, :\ \chain_{\permA}\subseteq \calD_{T}\,\}=S$, and,
by Lemma~\ref{lem.poset-lattice}(ii), $S\in {\cal X}^{\preced}$.
Finally, Lemma~\ref{lem.chain-lattice}(i) gives $\bigcup_{\permA\in S} \chain_{\permA}=\calD_{T}$.
The uniqueness of $S$ follows from the correspondence of ${\cal X}^{\outside}$
with set systems.

As concerns the sufficiency in (i), note that $S\in {\calX}^{\preced}$ with ${\calD}=\bigcup_{\permA\in S}\chain_{\permA}$ must be non-empty, as otherwise $\calD=\emptyset$. Then Lemma~\ref{lem.poset-lattice}(ii)
implies the existence of a poset $T$ over $N$ with ${\cal L}(T)=S$ and, by
Lemma~\ref{lem.chain-lattice}(i), $\calD_{T}= \bigcup_{\permA\in {\cal L}(T)} \chain_{\permA}=
\bigcup_{\permA\in S} \chain_{\permA}=\calD$. Lemma~\ref{lem.topol-lattice}(ii)(iii) then
yields the rest.

For (ii) assume $\emptyset\neq S\in{\calX}^{\outside}$ with ${\calD}:=\bigcup_{\permA\in S}\chain_{\permA}$ closed
under intersection and union. As $\calD$ is a topology distinguishing points
and one can apply (i) to get $\tilde{S}\in{\calX}^{\preced}$ with $\bigcup_{\permA\in \tilde{S}}\chain_{\permA}={\calD}$.
Observe that $\tilde{S}\in{\calX}^{\outside}$, by Lemma~\ref{lem.chain-lattice}.
The one-to-one correspondence between elements of ${\cal X}^{\outside}$
and the set systems $\calD$ then yields
$S=\tilde{S}\in{\calX}^{\preced}$.
\end{proof}

\subsection{Why enumerations instead of permutations}\label{ssec.why-enumeration}
This sub-section can be skipped without losing the understanding of the rest of the paper.
Its aim is to explain deeper reasons for our (notational) conventions.
\smallskip

We defined the set $\enuN$ as the collection of all enumerations of an \underline{unordered} set $N$
and introduced some geometric structure in it by interpreting its elements as vertices of a geometric
object in ${\dv R}^{N}$, namely of the permutohedron $\perN$. In particular, $\enuN$
can be viewed as the set of nodes of an undirected graph whose edges are defined geometrically:
they correspond to $1$-dimensional faces (= geometric edges) of the permutohedron.
There are also ``higher-dimensional" geometric relations among elements of $\enuN$,
for example belonging to the same $k$-dimensional face of the permutohedron $\perN$, $2\leq k< n$.

This geometric structure, namely the face lattice of\/ $\perN$ with an unordered set $N$,
is the relevant mathematical structure for our purposes. Let us emphasize particularly strongly that the set $\enuN$ is
{\em not a  group\/}! Indeed, there is no distinguished element within this set like the identity element within a group.

On the other hand, the elements of $\enuN$ are interpreted as vertices of a geometric object
in ${\dv R}^{N}$, namely of the permutohedron $\perN$, and this geometric object
{\em admits group actions\/} from the symmetric group on an $n$-element set, $n=|N|$.
The point is that two different group actions on $\enuN$ (= interpretations) are possible.
To explain that recall that the objects of our interest are (bijective) mappings between
an ordered {\em set of numbers\/} $[n]$ and an unordered {\em alphabet\/} (of symbols) $N$ with $n=|N|$.
Thus, one can distinguish
\begin{itemize}
\item an action from the group $S_{N}$ of permutations of the unordered alphabet $N$,
which is a composition with a bijective mapping $\alpha:N\to N$, and
\item an action from the group $S_{[n]}$ of permutations of the set of integers $[n]$,
which is a composition with a bijective mapping $\sigma:[n]\to [n]$.
\end{itemize}
Both these group actions give rise to self-transformations of $\enuN$, which behave in a different way,
but both of them are useful\/!
The permutations of the alphabet $N$ lead to linear self-transformations
of ${\dv R}^{N}$ preserving the geometric structure of\/ $\perN$, but these
are not all such transformations.
On the other hand, special permutations of the set $[n]$ of integers describe the local neighborhood structure
for any enumeration $\permA\in\enuN$. Specifically, these are the so-called {\em adjacent transpositions},
defined as the transpositions of numbers $i\leftrightarrow i+1$, where $i\in [n]\setminus\{n\}$.
In addition to that, another special ``reversing" permutation of $[n]$, defined by
$i\mapsto n+1-i$ for $i\in [n]$, leads to another linear self-transformation
of ${\dv R}^{N}$ preserving the geometry of $\perN$; it is different from actions of $S_{N}$ on $\enuN$ if $n\geq 3$, to be precise.
Note in this context that it is claimed in \cite{Cri11} (without a proof)
that this single transformation together with the transformations given by permutations of $N$
generate the whole symmetry group of\/ $\perN$ (= the group
of transformations preserving its geometric structure).
To summarize: our approach makes it possible to distinguish these
two different types of group actions very easily and in a natural way.
\smallskip

On the contrary, many authors dealing with this topic \cite{PRW08,MPSSW09,Cri11,MUWY18}
regard the basic set $N$ as an \underline{ordered} set from the very beginning. Typically, they simply 
put $N:=[n]$, where $n=|N|$, as they probably consider this to be a harmless convention, identifying the objects of interest with well-known mathematical objects, namely with {\em permutations} of a finite set.
Nevertheless, this step is equivalent to forcing a redundant, irrelevant and misleading mathematical structure on the
set of interest, namely the group structure\/! In its consequence, this convention leads to later difficulties with (correct) identification of objects of interest. The permutations of $[n]$ then
have too many meanings: they represent both vectors in the Euclidean space, that is,
$S_{[n]}\subseteq {\dv R}^{[n]}$, and two different kinds of self-transformations of this space $S_{[n]}$, namely the ``left" and ``right" actions of this group $S_{[n]}$ on itself\,!
This is surely messy for the reader and some of the above mentioned authors tried to overcome these identification problems by additional conventions, like introducing ``descent" vectors in \cite{MPSSW09}
or distinguishing left and right actions on $S_{[n]}$ in  \cite{Cri11}. Nevertheless, these additional conventions
are, in fact, enforced by the starting (unwise) convention that $N:=[n]$.

\section{Supermodularity and conditional independence}\label{sec.CI+supermod}
In this section we introduce concepts related to supermodular functions
and their conditional independence interpretation. Some particular notation will be utilized in this context.

Given $A\subseteq N$, the symbol $\delta_{A}$ will denote a zero-one
set function $\delta_{A}:\caP\to {\dv R}$ defined by $\delta_{A}(S):= \delta(\,A=S\,)$ for $S\subseteq N$. This is an identifier of the subset $A$ of $N$ in the space ${\dv R}^{\caP}$.\footnote{Our notation
distinguishes different vector interpretations of a set $A\subseteq N$: while $\delta_{A}$ is a vector in ${\dv R}^{\caP}$ (= a set function), the formerly introduced incidence vector $\chi_{A}$ is a vector in ${\dv R}^{N}$.}
This notational convention allows one to write formulas for elements of\/ ${\dv R}^{\caP}$.

For instance, given an ordered triplet $(A,B|C)$ of pairwise disjoint subsets of $N$, we put
$$
u_{(A,B|C)} := \delta_{ABC} +\delta_{C} -\delta_{AC} -\delta_{BC}\in {\dv Z}^{\caP}\,,\quad
\mbox{where canceling terms is possible.}
$$
This is a particular vector in ${\dv R}^{\caP}$ whose components are integers.
Given a set function $\gam:\caP\to {\dv R}$, the scalar product of\/ $\gam\in {\dv R}^{\caP}$ with $u_{(A,B|C)}$ in ${\dv R}^{\caP}$ is the respective ``supermodular" {difference expression}, occasionally denoted as follows:
$$
\diff\gam\, (A,B|C) ~:=~ \langle \gam,u_{(A,B|C)}\rangle =\, \gam(ABC)+\gam(C)-\gam(AC)-\gam(BC)\,.
$$
A triplet $(A,B|C)$ of pairwise disjoint subsets of $N$ will briefly be called a {\em triplet over N}.
It will be called {\em elementary\/}
if $|A|=1=|B|$, which means it has the form $(a,b|C)$, where $a,b\in N$ are distinct and $C\subseteq N\setminus ab$.

The class of elementary triplets over $N$ will be denoted by the symbol\, $\elemN$.
Elements of ${\dv Z}^{\caP}$ of the form $u_{(a,b|C)}$, $(a,b|C)\in\elemN$, will be
called {\em elementary imsets\/} over $N$. They admit conditional independence interpretation \cite{Stu05}
to be explained in Section~\ref{ssec.CI-iterpret}.

\subsection{Supermodular games and their core polytopes}\label{ssec.supermodular}
Given a finite non-empty basic set $N$, a {\em set function over $N$} is a
function $\gam:\caP\to {\dv R}$, that is, $\gam\in {\dv R}^{\caP}$. A (transferable-utility coalitional) {\em game\/} is modeled by a set function $\gam$ over $N$ satisfying $\gam(\emptyset)=0$.
A set function $\gam$ over $N$ is called {\em supermodular\/} if
$$
\gam(D\cup E) + \gam(D\cap E) ~\geq~ \gam(D) + \gam(E)\quad
\mbox{for every $D,E\subseteq N$.}
$$
It is a well-known fact that the class of supermodular functions over $N$ is a polyhedral
cone in ${\dv R}^{\caP}$. Its linearity space is the set of modular functions over $N$, which
has the dimension $n+1$, where $n=|N|$. The facet-defining inequalities for this cone
are precisely the inequalities $\langle \gam,u_{(a,b|C)}\rangle\geq 0$ for elementary imsets
over $N$ \cite[Corollary\,11]{KVV10}.

In the sequel, the symbol $\supmoN$ will denote the cone of {\em supermodular \underline{games} $\gam$}
delimited by the equality $\gam(\emptyset)=0$ and the (same) inequalities
$\langle \gam,u_{(a,b|C)}\rangle\geq 0$ for $(a,b|C)\in\elemN$.

Given a supermodular game $\gam$ over $N$, its {\em core\/} is a polyhedron
in ${\dv R}^{N}$ defined as follows:
$$
\cor(\gam) \,:=\, \{\, [z_{\ell}]_{\ell\in N}\in {\dv R}^{N} \,:\
\sum_{\ell\in N} z_{\ell}=\gam(N) ~~\&~~ \sum_{\ell\in S} z_{\ell}\geq \gam(S)\quad
\mbox{for every $S\subseteq N$}\,\}\,.
$$
It can be shown to be a bounded polyhedron, and, therefore, a polytope on ${\dv R}^{N}$.
A basic fact is that $\cor(\gam)\neq\emptyset$ for a supermodular game $\gam$ and, moreover,
the game $\gam$ is then {\em exact\/}, which means that the lower bounds defining the core
are tight. More formally, the game can be reconstructed from the
(vertices of the) core \cite[Theorem\,24(x),(xii)]{SK16}:
$$
\gam(S) \,=\, \min\, \{\, \sum_{\ell\in S} z_{\ell} \,:\ z\in \cor(\gam)\,\} \,=\,
\min\, \{\, \sum_{\ell\in S} y_{\ell} \,:\ y\in \ext(\cor(\gam))\,\}\quad
\mbox{for $S\subseteq N$}\,.
$$
The reader should be warned that the converse implication does not hold: a counter-example
exists in case $|N|=4$; see \cite[Example\,3]{SK16}.

\subsection{Vertices of the core polytope}\label{ssec.marginal-vectors}
A classic result in game theory by Shapley \cite{Sha72}
is the characterization of vertices of the core of a supermodular game over $N$ as certain vectors assigned to enumerations of $N$.

Given a \underline{game} $\gam:\caP\to {\dv R}$ over $N$ and an enumeration $\permA\in\enuN$,
the corresponding {\em marginal vector\/} is the vector $\mrg^{\gam}(\permA)$ in ${\dv R}^{N}$
determined by the values of $\gam$ on the respective maximal chain $\chain_{\permA}$
(see the formula \eqref{eq.max-chain} in Section~\ref{sssec.preposet-topologies})
as follows:
\begin{eqnarray}
\lefteqn{\mrg^{\gam}(\permA)=[y_{\ell}]_{\ell\in N},\quad
\mbox{where}~~ y_{\permA(i)} \,:=\,
\gam\, (\,\bigcup_{j=1}^{i} \permA(j)\,) \,-\,
\gam\, (\,\bigcup_{j=1}^{i-1} \permA(j)\,)
~~\mbox{for $i\in [n]$,}\label{eq.marg-vector}} \\
&\mbox{alternative writing is}&
\mrg^{\gam}(\permA)_{\ell} ~:=~
\gam\, (\bigcup_{j\leq\,\permA_{-1}(\ell)} \permA(j)\,) \,-\,
\gam\, (\bigcup_{j<\,\permA_{-1}(\ell)} \permA(j)\,)
\quad
\mbox{for $\ell\in N$.}
\nonumber
\end{eqnarray}
Given a game $\gam:\caP\to {\dv R}$ over $N$ and a vector $y\in{\dv R}^{N}$,
the respective {\em tightness class\/} for $\gam$ and $y$ is the collection
of subsets of $N$ at which $y$ is tight for $\gam$, denoted by
\begin{equation}
\tight^{\gam}_{y} ~:=~ \{\, S\subseteq N\,:\ \sum_{\ell\in S} y_{\ell} =\gam(S)\,\,\}\,.
\label{eq.tight-class}
\end{equation}
It is nearly immediate, for $\permA\in\enuN$ and $y\in {\dv R}^{N}$, that
\begin{equation}
\permA\in\mrg^{\gam}_{-1}(\{y\})
\quad\Leftrightarrow\quad
y=\mrg^{\gam}(\permA)
\quad\Leftrightarrow\quad
\chain_{\permA}\subseteq \tight^{\gam}_{y}\,.
\label{eq.marg-chain}
\end{equation}
Thus, it follows from the definitions in Section~\ref{ssec.weak-lattice} that
$\mrg^{\gam}_{-1}(\{y\})\in {\calX}^{\outside}$ for any game $\gam$ over $N$ and $y\in {\dv R}^{N}$.
In fact, one can show, for $S\subseteq\enuN$, that $S\in{\calX}^{\outside}$ iff there exists a game $\gam$ over $N$ and $y\in {\dv R}^{N}$ with $S=\mrg^{\gam}_{-1}(\{y\})$; the proof is, however, omitted in this paper.

The relation \eqref{eq.marg-vector} defines a {\em marginal-vector mapping\/} $\mrg^{\gam}:\enuN\to {\dv R}^{N}$
which need not be injective. The definition can be illustrated by the following example.

\begin{figure}
\setlength{\unitlength}{0.65mm}
\begin{center}
\scalebox{0.7}{\begin{picture}(90,95)
%
\put(30,79){\makebox(0,-8){\tiny $\{a,b,c\}$}}
\put(30,79){\makebox(0,4){\rm\normalsize $2$}}
\put(30,79){\makebox(0,-4.5){\line(1,0){17}}}
\put(6,55){\oval(18,12)}
\put(6,55){\makebox(0,-8){\tiny $\{a,b\}$}}
\put(6,55){\makebox(0,4){\rm\normalsize 1}}
\put(6,55){\makebox(0,-4.5){\line(1,0){17}}}
\put(30,55){\oval(18,12)}
\put(30,55){\makebox(0,-8){\tiny $\{a,c\}$}}
\put(30,55){\makebox(0,4){\rm\normalsize 1}}
\put(30,55){\makebox(0,-4.5){\line(1,0){17}}}
%
\put(54,55){\makebox(0,-8){\tiny $\{b,c\}$}}
\put(54,55){\makebox(0,4){\rm\normalsize 1}}
\put(54,55){\makebox(0,-4.5){\line(1,0){17}}}
\put(67,18){\makebox(0,0){\small $\permA(1)=c;\, y_{c}=0-0$}}
\put(79,43){\makebox(0,0){\small $\permA(2)=b;\, y_{b}=1-0$}}
\put(69,68){\makebox(0,0){\small $\permA(3)=a;\, y_{a}=2-1$}}
\put(6,31) {\oval(18,12)}
\put(6,31) {\makebox(0,-8){\tiny $\{a\}$}}
\put(6,31) {\makebox(0,4){\rm\normalsize $0$}}
\put(6,31){\makebox(0,-4.5){\line(1,0){17}}}
\put(30,31){\oval(18,12)}
\put(30,31){\makebox(0,-8){\tiny $\{b\}$}}
\put(30,31){\makebox(0,4){\rm\normalsize $0$}}
\put(30,31){\makebox(0,-4.5){\line(1,0){17}}}
%
\put(54,31){\makebox(0,-8){\tiny $\{c\}$}}
\put(54,31){\makebox(0,4){\rm\normalsize $0$}}
\put(54,31){\makebox(0,-4.5){\line(1,0){17}}}
%
\put(30,7){\makebox(0,-8){\tiny $\emptyset$}}
\put(30,7){\makebox(0,4){\rm\normalsize $0$}}
\put(30,7){\makebox(0,-4.5){\line(1,0){17}}}
\put(27,73){\line(-3,-2){18}}
\put(30,73){\line( 0,-1){12}}
%
\put(6,49){\line( 0,-1){12}}
\put(9,49){\line( 3,-2){18}}
\put(27,49){\line(-3,-2){18}}
\put(33,49){\line( 3,-2){18}}
\put(51,49){\line(-3,-2){18}}
%
\put(30,25){\line( 0,-1){12}}
\put(9,25){\line( 3,-2){18}}
\thicklines
\put(30,79){\oval(18,12)}
\put(54,55){\oval(18,12)}
\put(54,31){\oval(18,12)}
\put(30,7){\oval(18,12)}
\put(54,49){\line( 0,-1){12}}
\put(51,25){\line(-3,-2){18}}
\put(33,73){\line( 3,-2){18}}
\end{picture}
}
\end{center}
\caption{A picture illustrating the definition of a marginal vector from Example~\ref{exa.marg-vector}.\label{fig.4}}
\end{figure}
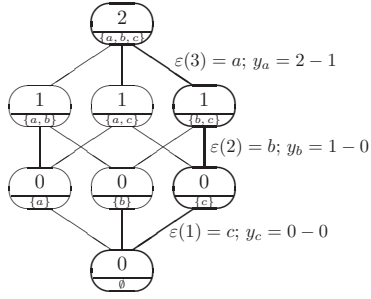

\begin{example}\label{exa.marg-vector}\rm
Put $N:=\{a,b,c\}$ and $\gam:= 2\cdot\delta_{N}+\sum_{S\subseteq N:|S|=2} \delta_{S}$.
Consider $\permA :=|c|b|a|$ and compute $\mrg^{\gam}(\permA)=[y_{a},y_{b},y_{c}]$.
See Figure~\ref{fig.4} for illustration, where the respective maximal chain $\chain_{\permA}$ is in bold.
As $\permA(1)=c$ one has $y_{c}\stackrel{\eqref{eq.marg-vector}}{=}\gam(\permA(1))-\gam(\emptyset)=\gam(c)-\gam(\emptyset)=0$, analogously $\permA(2)=b$ gives $y_{b}=\gam(\,\permA(1)\permA(2)\,)-\gam(\permA(1))=\gam(bc)-\gam(c)=1$,
and $\permA(3)=a$ gives $y_{a}=\gam(abc)-\gam(bc)=1$. Altogether, $\mrg^{\gam}(\permA)=[y_{a},y_{b},y_{c}]=[1,1,0]$. Observe that for $\permB :=|c|a|b|$ one also has
$\mrg^{\gam}(\permB)=[y_{a},y_{b},y_{c}]=[1,1,0]=\mrg^{\gam}(\permA)$.
\end{example}

Here is a classic result, proved as \cite[Corollary\,3]{SK16}; see also \cite[Theorems\, 3 and 5]{Sha72}.

\begin{lemma}\label{lem.marg-vectors}\rm
Let $\gam:\caP\to {\dv R}$ be a supermodular game over $N$. Then the set of vertices
of the core $\cor(\gam)$ coincides with the set of (distinct) marginal vectors
for enumerations of $N$. Formally: $\ext(\cor (\gam)) \,=\, \{\, \mrg^{\gam}(\permA)\, :\
\permA\in\enuN\,\}$.
\end{lemma}

This allows one to derive the claim about the permutohedron $\perN$ mentioned in the end of Section~\ref{ssec.enum-permutohedron}.
Consider a supermodular game $\gam$ defined by $\gam(S):={|S|+1\choose 2}$ for $S\subseteq N$.
Then Lemma~\ref{lem.marg-vectors} allows one to conclude that $\ext(\cor(\gam))=\{\, \rankv_{\permA}\,:\ \permA\in\enuN\,\}$. Hence, by Krein-Milman theorem, $\perN =\conv(\{\, \rankv_{\permA}\,:\ \permA\in\enuN\,\})= \conv( \ext(\cor(\gam)))= \cor(\gam)$,
and $\cor(\gam)$ yields a polyhedral description of the permutohedron over $N$.

\subsection{Relation of supermodularity and the lattice based on posets}\label{ssec.supermo-posets}
In this section we relate the marginal-vector mappings $\permA\in\enuN\mapsto\mrg^{\gam}(\permA)\in {\dv R}^{N}$
for $\gam\in\supmoN$ defined in Section~\ref{ssec.marginal-vectors} and the lattice based on posets from Section~\ref{ssec.poset-lattice}.

\begin{lemma}\label{lem.supermo-marg-1}\rm
Let $\gam\in\supmoN$ be a supermodular game and $y\in\mrg^{\gam}(\enuN)$.
Then $\tight^{\gam}_{y}$ is closed under intersection and union and
one has $\tight^{\gam}_{y}=\bigcup_{\permA\in S} \chain_{\permA}$ with $S=\mrg^{\gam}_{-1}(\{y\})$.
\end{lemma}

Thus, the tightness classes $\tight^{\gam}_{y}$ are topologies on $N$
distinguishing points (see Section~\ref{sssec.preposet-topologies}).

\begin{proof}
Lemma~\ref{lem.marg-vectors} says $y\in\ext(\cor(\gam ))$. Consider $D,E\in\tight^{\gam}_{y}$ and write
\begin{eqnarray*}
\lefteqn{\hspace*{-22mm}0\geq \gam(D) + \gam(E) - \gam(D\cup E) - \gam(D\cap E) 
= \sum_{u\in D} y_{u} +  \sum_{v\in E} y_{v}  - \gam(D\cup E) - \gam(D\cap E)} \\
&=& [\, \underbrace{\sum_{u\in D\cup E} y_{u}  - \gam(D\cup E)}_{\geq 0}\,] \,+\,
[\, \underbrace{\sum_{\ell\in D\cap E} y_{\ell} - \gam(D\cap E)}_{\geq 0}\,] \geq 0\,,
\end{eqnarray*}
where the first inequality follows from supermodularity of $\gam$, the first equality
from the tightness assumptions on $D$ and $E$, and the remaining inequalities from the fact $y\in\cor(\gam )$.
All the inequalities must thus be equalities, which enforces $D\cup E, D\cap E\in\tight^{\gam}_{y}$.

By \eqref{eq.marg-chain}, the existence of $\permB\in\enuN$ with $\mrg^{\gam}(\permB)=y$ gives
$\chain_{\permB}\subseteq \tight^{\gam}_{y}$. Hence, $\emptyset,N\in\tight^{\gam}_{y}$ and the
class $\calD:=\tight^{\gam}_{y}$ is a topology on $N$ which distinguishes points.
Apply Corollary~\ref{cor.set-systems}(i) to $\calD$ with $S:=\mrg^{\gam}_{-1}(\{y\})\stackrel{\eqref{eq.marg-chain}}{=}\{\, \permA\in\enuN\,:\, \chain_{\permA}\subseteq \tight^{\gam}_{y} \}$, which yields the desired conclusion $\calD=\bigcup_{\permA\in S} \chain_{\permA}$.
\end{proof}

Nonetheless, Lemma~\ref{lem.supermo-marg-1} need not hold without the supermodularity assumption.

\begin{figure}
\setlength{\unitlength}{0.65mm}
\begin{center}
\scalebox{0.7}{\begin{picture}(90,95)
%
\put(30,79){\makebox(0,-8){\tiny $\{a,b,c\}$}}
\put(30,79){\makebox(0,4){\rm\normalsize $+3$}}
\put(30,79){\makebox(0,-4.5){\line(1,0){17}}}
\put(6,55){\oval(18,12)}
\put(6,55){\makebox(0,-8){\tiny $\{a,b\}$}}
\put(6,55){\makebox(0,4){\rm\normalsize $0$}}
\put(6,55){\makebox(0,-4.5){\line(1,0){17}}}
%
\put(30,55){\makebox(0,-8){\tiny $\{a,c\}$}}
\put(30,55){\makebox(0,4){\rm\normalsize $+2$}}
\put(30,55){\makebox(0,-4.5){\line(1,0){17}}}
%
\put(54,55){\makebox(0,-8){\tiny $\{b,c\}$}}
\put(54,55){\makebox(0,4){\rm\normalsize $+2$}}
\put(54,55){\makebox(0,-4.5){\line(1,0){17}}}
%
\put(6,31) {\makebox(0,-8){\tiny $\{a\}$}}
\put(6,31) {\makebox(0,4){\rm\normalsize $+1$}}
\put(6,31){\makebox(0,-4.5){\line(1,0){17}}}
\put(30,31){\oval(18,12)}
\put(30,31){\makebox(0,-8){\tiny $\{b\}$}}
\put(30,31){\makebox(0,4){\rm\normalsize $-1$}}
\put(30,31){\makebox(0,-4.5){\line(1,0){17}}}
\put(54,31){\oval(18,12)}
\put(54,31){\makebox(0,-8){\tiny $\{c\}$}}
\put(54,31){\makebox(0,4){\rm\normalsize $0$}}
\put(54,31){\makebox(0,-4.5){\line(1,0){17}}}
%
\put(30,7){\makebox(0,-8){\tiny $\emptyset$}}
\put(30,7){\makebox(0,4){\rm\normalsize $0$}}
\put(30,7){\makebox(0,-4.5){\line(1,0){17}}}
\put(27,73){\line(-3,-2){18}}
\put(33,73){\line( 3,-2){18}}
\put(6,49){\line( 0,-1){12}}
\put(9,49){\line( 3,-2){18}}
%
\put(33,49){\line( 3,-2){18}}
\put(51,49){\line(-3,-2){18}}
\put(54,49){\line( 0,-1){12}}
\put(51,25){\line(-3,-2){18}}
\put(30,25){\line( 0,-1){12}}
\thicklines
\put(30,79){\oval(18,12)}
\put(30,55){\oval(18,12)}
\put(54,55){\oval(18,12)}
\put(30,7){\oval(18,12)}
\put(6,31) {\oval(18,12)}
\put(30,73){\line( 0,-1){12}}
\put(9,25){\line( 3,-2){18}}
\put(27,49){\line(-3,-2){18}}
\end{picture}
}
\end{center}
\caption{A picture illustrating Example~\ref{exa.supermo-marg}.\label{fig.5}}
\end{figure}
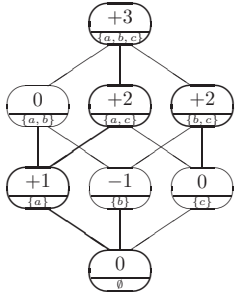

\begin{example}\label{exa.supermo-marg}\rm
Put $N:=\{a,b,c\}$ and $\gam:= 3\cdot\delta_{N}+2\cdot\delta_{ac}+2\cdot\delta_{bc}
+\delta_{a}-\delta_{b}$. Consider the constant vector $y=[1,1,1]\in {\dv R}^{N}$;
one has $S:=\mrg^{\gam}_{-1}(\{y\})=\{\permB\}$ where $\permB :=|a|c|b|$.
Observe that $\tight^{\gam}_{y}=\{\, \emptyset,a,ac,bc,abc \,\}$
is not closed under intersection because $ac\cap bc=c\not\in \tight^{\gam}_{y}$.
See Figure~\ref{fig.5} for illustration, where the elements of the tightness class $\tight^{\gam}_{y}$ are in bold.
Moreover, $bc\in\tight^{\gam}_{y}\setminus\bigcup_{\permA\in S} \chain_{\permA}$.
The reason why $\gam$ is not supermodular is $\diff\gam\, (a,b|c)=-1$.
\end{example}

Elements of the lattice $({\calX}^{\preced},\subseteq)$ from Section~\ref{ssec.poset-lattice}
appear to be just the inverse images of singletons with the marginal-vector mappings for supermodular games.

\begin{lemma}\label{lem.supermo-marg-2}\rm
If $n\geq 2$ and $S\subseteq\enuN$ then $S\in {\calX}^{\preced}$ iff
[\,$\exists\, \gam\in\supmoN ~\exists\, y\in {\dv R}^{N}\, :~\, S= \mrg^{\gam}_{-1}(\{y\})$\,].
\end{lemma}

\begin{proof}
For the sufficiency assume $S\neq\emptyset$ for otherwise the claim is evident.
Apply Lemma~\ref{lem.supermo-marg-1} to $S= \mrg^{\gam}_{-1}(\{y\})$ and observe that $\calD:=\bigcup_{\permA\in S} \chain_{\permA}$ is closed under intersection and union.
Because of $S=\mrg^{\gam}_{-1}(\{y\})\stackrel{\eqref{eq.marg-chain}}{=}\{\, \permA\in\enuN\,:\, \chain_{\permA}\subseteq \calD \}$ one has $S\in {\calX}^{\outside}$.
By Corollary~\ref{cor.set-systems}(ii) get $S\in {\calX}^{\preced}$.

For the necessity it is enough to verify that the class of sets $S\subseteq\enuN$ satisfying
\begin{equation}
\exists\, \gam\in\supmoN ~~\exists\, y\in {\dv R}^{N}~ :\quad S= \mrg^{\gam}_{-1}(\{y\})
\label{eq.supermo-marg-2}
\end{equation}
is a Moore family of subsets of $\enuN$ involving coatoms of ${\calX}^{\preced}$. We prove this in three steps.

\noindent {\bf I.} The set $S:=\enuN$ satisfies \eqref{eq.supermo-marg-2}.\\[0.2ex]
To this end put $\gam(S):=0$ for any $S\subseteq N$ and $y=[y_{\ell}]_{\ell\in N}$ with $y_{\ell}:=0$
for any $\ell\in N$.

\noindent {\bf II.} Any set $S_{u\prec v}:= \{\,\permA\in\enuN\,:\ u\prec_{\permA}v\,\}$
with distinct $u,v\in N$ satisfies \eqref{eq.supermo-marg-2}.\\[0.2ex]
To this end put $\gam(S):=\delta(\,\{u,v\}\subseteq S\,)$ for any $S\subseteq N$.
It is a supermodular game and one has $\mrg^{\gam}(\permA)=\chi_{v}$ for $\permA\in S_{u\prec v}$
while $\mrg^{\gam}(\permB)=\chi_{u}$ for $\permB\in S_{v\prec u}$. Hence,
$S_{u\prec v}=\mrg^{\gam}_{-1}(\{\chi_{v}\})$.

\noindent {\bf III.} If $S_{1},S_{2}\subseteq\enuN$ satisfy \eqref{eq.supermo-marg-2}
then $S_{1}\cap S_{2}$ satisfies \eqref{eq.supermo-marg-2}.\\[0.2ex]
Indeed, the assumption means that, for $i=1,2$, there exist games $\gam_{i}\in\supmoN$ and $y^{i}\in {\dv R}^{N}$ such that
$S_{i}=\{\, \permA\in\enuN\,:\ y^{i}=\mrg^{\gam_{i}}(\permA)\,\}$. We choose $\beta>0$ such that
$$
\beta \,>\, \max\, \left\{\,
\frac{y^{1}_{u}-\mrg^{\gam_{1}}(\permB)_{u}}{\mrg^{\gam_{2}}(\permB)_{u}-y^{2}_{u}} \,:\ \permB\in\enuN ~~\mbox{and $u\in N$ that satisfy}~~
y^{2}_{u}\neq \mrg^{\gam_{2}}(\permB)_{u}
\right\}\,.
$$
The existence of $\beta$ follows from the fact that the set above is finite. We put $\gam:=\gam_{1}+\beta\cdot\gam_{2}$ and $y:=y^{1}+\beta\cdot y^{2}$. Clearly, $\gam\in\supmoN$ is also a supermodular game and the linearity of the transformation $\gam\mapsto \mrg^{\gam}(\permA)$ defined in \eqref{eq.marg-vector}
(for any fixed $\permA$) yields
$$
\permA\in S_{1}\cap S_{2}~\Rightarrow~
\mrg^{\gam}(\permA)=\mrg^{\gam_{1}}(\permA)+\beta\cdot\mrg^{\gam_{2}}(\permA)=y^{1}+\beta\cdot y^{2}=y\,.
$$
It remains to show that $\permB\in \enuN\setminus S_{1}\cap S_{2}~\Rightarrow~
y\neq \mrg^{\gam}(\permB)$.
\begin{itemize}
\item If $\permB\not\in S_{2}$ then $u\in N$ exists with $y^{2}_{u}\neq \mrg^{\gam_{2}}(\permB)_{u}$.
Assume for a contradiction that $y=\mrg^{\gam}(\permB)$ which implies $y_{u}=\mrg^{\gam}(\permB)_{u}$
and then
\begin{eqnarray*}
\lefteqn{\hspace*{-46mm}y^{1}_{u}+\beta\cdot y^{2}_{u}=\mrg^{\gam_{1}}(\permB)_{u}+\beta\cdot \mrg^{\gam_{2}}(\permB)_{u}
~~\Rightarrow~~
[\, y^{1}_{u}-\mrg^{\gam_{1}}(\permB)_{u}\,] = \beta\cdot [\,\mrg^{\gam_{2}}(\permB)_{u}-y^{2}_{u}\,]}\\[0.3ex]
&\Rightarrow&
\beta = \frac{y^{1}_{u}-\mrg^{\gam_{1}}(\permB)_{u}}{\mrg^{\gam_{2}}(\permB)_{u}-y^{2}_{u}}\,,
\end{eqnarray*}
which is in contradiction with the definition of $\beta$.
\item If $\permB\not\in S_{1}$ then $u\in N$ exists with
$y^{1}_{u}\neq \mrg^{\gam_{1}}(\permB)_{u}$. Assume for a contradiction that $y=\mrg^{\gam}(\permB)$
and observe
$[\, y^{1}_{u}-\mrg^{\gam_{1}}(\permB)_{u}\,] = \beta\cdot [\,\mrg^{\gam_{2}}(\permB)_{u}-y^{2}_{u}\,]$
by the same consideration as in the previous case. Since the left-hand side is non-zero it implies
$\mrg^{\gam_{2}}(\permB)_{u}-y^{2}_{u}\neq 0$ and the contradiction can be derived as in the previous
case.
\end{itemize}
These facts together give $S_{1}\cap S_{2}=\{\,\permA\in\enuN\,:\ y=\mrg^{\gam}(\permA)\,\}$ and
$S_{1}\cap S_{2}$ satisfies \eqref{eq.supermo-marg-2}.

By Lemma~\ref{lem.poset-lattice}, the lattice $({\calX}^{\preced},\subseteq)$ is coatomistic, with sets $S_{u\prec v}$ as coatoms, and
the infimum in it is the intersection, the observations above imply that any set $S\in{\calX}^{\preced}$
satisfies \eqref{eq.supermo-marg-2}.
\end{proof}

\subsection{Geometric view on enumerations}\label{ssec.normal-cone}
To read this section, the reader is advised to recall the concepts defined in Appendix~\ref{ssec.app.normal.fan}.

\begin{lemma}\label{lem.enum-braid}\rm
Let us denote by ${\dv R}^{N}_{\neq}$ the set of\/ $x\in {\dv R}^{N}$ with pairwise distinct components.
\begin{itemize}
\item[(i)] Given a full-dimensional polyhedral cone $\nor\subseteq {\dv R}^{N}$,
the set ${\dv R}^{N}_{\neq}$ is dense in $\nor$ relative to Euclidean topology on ${\dv R}^{N}$, that is, one has $\nor=\overline{\nor\cap {\dv R}^{N}_{\neq}}$.
\item[(ii)] Given a poset $T\subseteq N\times N$ on $N$, one has
$\nor_{T}=\bigcup_{\permA\in {\cal L}(T)} \nor^{\permA}$. (see Section~\ref{sssec.preposet-cones})
\end{itemize}
\end{lemma}

\begin{proof}
For (i) realize that $\nor$ is the closure of its (relative)
interior (Appendix~\ref{sec.app.polyhedral}),
and, for any $x\in\nor$, a sequence $x_{k}\in\inte(\nor)$ exists which converges to $x$. Consider
balls $B(x_{k},\epsilon_{k})\subseteq \inte(\nor)$ in the Euclidean metric
with shrinking diameters $\epsilon_{k}>0$, that is, $\epsilon_{k}\to 0$.
The complement of ${\dv R}^{N}_{\neq}$ is the union of finitely many linear subspaces of
the dimension $|N|-1$. Thus, for each $k$, there is at least one vector
$\tilde{x}_{k}\in B(x_{k},\epsilon_{k})\cap {\dv R}^{N}_{\neq}\subseteq \nor\cap {\dv R}^{N}_{\neq}$.
As $\epsilon_{k}\to 0$, the sequence $\tilde{x}_{k}$ converges to $x$.

As concerns (ii), the inclusion $\nor^{\permA}\subseteq \nor_{T}$ for $\permA\in {\cal L}(T)$ is
immediate. 
For the other inclusion consider $x\in \nor_{T}\cap {\dv R}^{N}_{\neq}$ and realize that unique $\permB\in\enuN$ exists with $x\in \nor^{\permB}$. Moreover, for any $(u,v)\in N\times N$, $u\neq v$, one has $x_{u}<x_{v}$ iff
$u\prec_{\permB} v$, which allows one to conclude that $\permB\in {\cal L}(T)$.
Hence, $\nor_{T}\cap {\dv R}^{N}_{\neq}\subseteq \bigcup_{\permA\in {\cal L}(T)} \nor^{\permA}$ and the application of (i) yields the other inclusion.
\end{proof}

The following are basic observations on Weyl chambers $\nor^{\permA}$ defined by \eqref{eq.enum-cone} in Section~\ref{sssec.preposet-cones}.

\begin{lemma}\label{lem.braid-cone}\rm
Given $\permA\in\enuN$, one has $\nor^{\permA}\,=\,\cone
(\,\{\,\chi_{N}\}\cup\{-\chi_{D}\,:\ D\in\chain_{\permA}\}\,)$.
The topological interior of $\nor^{\permA}$
(relative to Euclidean topology in ${\dv R}^{N}$) has the form
$$
\inte(\nor^{\permA}) \,:=~ \{\, x\in {\dv R}^{N}\,:\ x_{\permA(1)} < x_{\permA(2)} < \ldots < x_{\permA(n)}\,\,\}\,.
$$
The cone $\nor^{\permA}$ is the outer normal cone to the permutohedron $\perN$ at the rank vector $\rankv_{\permA}$.\\
The braid fan is the normal fan of $\perN$, its maximal cones have the form $\nor^{\permA}$ for $\permA\in\enuN$.\\
If $n\geq 2$, then, given $S\subseteq\enuN$, one has $S\in{\calX}^{\preced}$ iff $\bigcup_{\permA\in S} \nor^{\permA}$ is convex.
\end{lemma}

\begin{proof}
The first claim follows from Lemma~\ref{lem.topol-lattice}(iv) applied to the toset $T_{\permA}$ and topology $\calD:=\chain_{\permA}$;
indeed, the respective braid cone is $\nor^{\permA}$ (see Section~\ref{sssec.preposet-cones}).

The claim about the interior follows from the fact that the (relative) interior of a closed convex set $\nor^{\permA}$ is the relative complement of the union of its non-empty proper faces, which is nothing
but the union of its facets (see Appendix~\ref{ssec.app.vertices-facets}).

Recall from Section~\ref{ssec.face-lattice} that
edges of the permutohedron $\perN$ correspond to ordered partitions of $N$ whose blocks are singletons except
one two-element set. Hence, the geometric neighbours of a vertex $\rankv_{\permA}$ in $\perN$
are just the rank vectors $\rankv_{\permB}$, where $\permB$ differs from $\permA$ by an adjacent transposition.
Specifically, if $\permA$ and $\permB$ differ in the $i$-th two positions, where $i\in[n-1]$,
as described in \eqref{eq.adjacent-transpos}, then one has
$\rankv_{\permA}-\rankv_{\permB}=\chi_{\permA(i+1)}-\chi_{\permA(i)}$.
The polyhedral description of the normal cone at a vertex
described in Appendix~\ref{ssec.app.normal.fan} implies the following:
\begin{eqnarray*}
\lefteqn{\hspace*{-15mm}\nor_{\perNs}(\rankv_{\permA}) =\{\, x\in {\dv R}^{N}\,:~~
\forall\, \rankv_{\permB}\in\ne_{\,\perNs}(\rankv_{\permA})\quad
\langle x,\rankv_{\permA}-\rankv_{\permB}\rangle\geq 0\,\}}\\
&=&
\{\, x\in {\dv R}^{N}\,:~~ \forall\, i\in[n-1]\quad
\langle x,\chi_{\permA(i+1)}-\chi_{\permA(i)}\rangle\geq 0\,\} ~\stackrel{\eqref{eq.enum-cone}}{=}~ \nor^{\permA}\,.
\end{eqnarray*}

The observation about the form of the interiors of $\nor^{\permA}$ for $\permA\in\enuN$ implies that these sets
are not cut by the hyperplanes defining the braid fan. This fact allows one to derive the claim about the form of maximal cones of the braid fan. The rest is then evident.

The necessity in the last claim in Lemma~\ref{lem.braid-cone} follows from Lemma~\ref{lem.poset-lattice}(ii): given $\emptyset\neq S\in {\calX}^{\preced}$ take a poset $T\subset N\times N$ with $S={\cal L}(T)$. By Lemma~\ref{lem.enum-braid}(ii), the set $\bigcup_{\permA\in S} \nor^{\permA}=\nor_{T}$ is convex.
\smallskip

For the sufficiency assume non-empty $S\subseteq\enuN$ such that $\nor:=\bigcup_{\permA\in S} \nor^{\permA}$ is convex; note that if $S=\emptyset$ then $n\geq 2$ gives immediately $S\in{\calX}^{\preced}$.\\
{\bf I.}~ The first step is to observe that $\nor$ is a polyhedral cone.\\[0.3ex]
Indeed, this basically follows from the first claim of Lemma~\ref{lem.braid-cone} saying that
any $\nor^{\permA}$ is the conic hull of a finite collection of vectors in ${\dv R}^{N}$.
Hence, $\cone (\bigcup_{\permA\in S} \nor^{\permA})$ is also a polyhedral cone. Nevertheless, the convex hull and the conic hull of the union of a finite collection of polyhedral cones
coincide, which gives $\nor=\bigcup_{\permA\in S} \nor^{\permA}=\conv(\bigcup_{\permA\in S} \nor^{\permA})=\cone (\bigcup_{\permA\in S} \nor^{\permA})$.

Thus, $\nor$ is a full-dimensional polyhedral cone. Assume $\nor\neq {\dv R}^{N}$ as otherwise $S=\enuN$ and we are done. Then $\nor$ has unique minimal polyhedral description by means of its facet-defining
inequalities (see Appendix~\ref{ssec.app.vertices-facets}).\\
{\bf II.}~ Observe that any facet-defining inequality for $\nor$ has the form $x_{u}\leq x_{v}$ for $u,v\in N$, $u\neq v$.\\[0.3ex]
Indeed, since $\nor$ is a cone, any facet-defining inequality for $x\in\nor$ has the form $\langle x,y\rangle\geq 0$, with non-zero $y\in {\dv R}^{N}$, and the dimension of\/ $\{ x\in\nor\,:\ \langle x,y\rangle =0\}$ is $|N|-1$.
Hence, $\permA\in S$ exists such that $\{ x\in\nor^{\permA}\,:\ \langle x,y\rangle =0\}$ has the dimension $|N|-1$,
implying that $\langle x,y\rangle\geq 0$ is facet-defining for $x\in \nor^{\permA}$. This implies that the inequality has the required form; use \eqref{eq.enum-cone}.

Therefore, $\nor$ is a full-dimensional braid cone, and, by application of Lemma~\ref{lem.braid-lattice}(iii), a poset $T\subset N\times N$ over $N$ exists with $\nor=\nor_{T}$.\\
{\bf III.}~ The third step is to observe that one has
$S={\cal L}(T)$.\\[0.3ex]
Indeed, if $\permA\in S$ then choose $x\in\inte (\nor^{\permA})$ and observe that, for any
$(u,v)\in T$, one must have $x_{u}<x_{v} ~\Rightarrow~ u\prec_{\permA}v$, giving $\permA\in {\cal L}(T)$.
Conversely, assuming $\permA\in {\cal L}(T)$, for any $(u,v)\in T$ and $x\in\nor^{\permA}$, one has $u\preceq_{\permA} v$ and
$x_{u}\leq x_{v}$, implying $\nor^{\permA}\subseteq\nor_{T}$.
Therefore, $\nor^{\permA}\subseteq\nor_{T}=\bigcup_{\permB\in S} \nor^{\permB}$.
Because $\inte(\nor^{\permA})$ is disjoint with $\nor^{\permB}$ for $\permB\neq\permA$
(see Appendix~\ref{ssec.app.normal.fan}) this forces $\permA\in S$.

Altogether, $S={\cal L}(T)$ means, by Lemma~\ref{lem.poset-lattice}(ii), that $S\in{\calX}^{\preced}$.
\end{proof}

We can now gather various characterizations of poset-based sets of enumerations.

\begin{corol}\rm\label{cor.braid-cone}
The following conditions are equivalent for a \underline{non-empty} set $S\subseteq\enuN$.
\begin{itemize}
\item[(i)] There exists $T\subseteq N\times N$ such that $S={\cal L}(T)$,\quad
(that is, $S\in {\calX}^{\preced}$)
\item[(ii)] there exists a poset $T\subseteq N\times N$ on $N$ such that $S={\cal L}(T)$,
\item[(iii)] $S$ is geodetically convex (in the permutohedral graph),
\item[(iv)] $\forall\, \permA,\permB\in S~\forall\,\permC\in\enuN$  if\,
$[\,\forall\, u,v\in N ~~ u\prec_{\permA}v\,\,\&\,\, u\prec_{\permB}v ~\Rightarrow~ u\prec_{\permC}v \,]$ then $\permC\in S$,
\item[(v)] $\bigcup_{\permA\in S} \nor^{\permA}$ is a convex set in ${\dv R}^{N}$.
\end{itemize}
\end{corol}

Note that the combinatorial 
condition (iv) was named ``pre-convexity" in \cite[\S\,2.2]{Mor07};
the argument there was that it is implied by the geometric convexity condition (v).

\begin{proof}
Assume $n\geq 2$ because all conditions hold trivially in case $n=1$.
The equivalence (i)$\Leftrightarrow$(ii) was shown in Lemma~\ref{lem.poset-lattice}(ii),
(i)$\Leftrightarrow$(iii) in Theorem~\ref{thm.poset-characterization};
(iii)$\Leftrightarrow$(iv) follows from Lemma~\ref{lem.permut-graph}(iii) and (i)$\Leftrightarrow$(v) from the last claim of Lemma~\ref{lem.braid-cone}.
\end{proof}

The next lemma relates the marginal-vector mapping from Section~\ref{ssec.marginal-vectors}
and the (maximal) normal cones to the core polytope of a supermodular game.

\begin{lemma}\label{lem.refine-fan}\rm
Given a supermodular game $\gam\in\supmoN$ and $y\in\ext(\cor(\gam))$, one has
\begin{itemize}
\item[(a)] if $\mrg^{\gam}(\permA)=y$ for $\permA\in\enuN$
then~ $\nor^{\permA}\subseteq\nor_{\cors(\gam)}(y)$,
\item[(b)] if $\mrg^{\gam}(\permB)\neq y$ for $\permB\in\enuN$
then~ $\inte(\nor^{\permB})\cap\nor_{\cors(\gam)}(y)=\emptyset$.
\end{itemize}
These two observations together yield $\nor_{\cors(\gam)}(y)=\bigcup_{\permA\in\mrg^{\gam}_{-1}(\{y\})} \nor^{\permA}$, which further implies that $\nor_{\cors(\gam)}(y)=T^{\lGbraid}$ for the
poset $T:=\mrg^{\gam}_{-1}(\{y\})^{\rGprec}$.\\[0.2ex]
Moreover, the poset $T$ coincides with the transitive closure $\tr(R)$ of the relation

\vspace*{-5mm}
$$
R:=\{\, (u,v)\in N\times N\, :\ \exists\, z\in\ext (\cor(\gam)) ~\exists\, k>0 ~\mbox{with}~ z-y=k\cdot(\chi_{u}-\chi_{v})\,\}\,.
$$
\end{lemma}

\begin{proof}
Assume $n\geq 2$ as otherwise all claims trivially hold.
To verify (a) we realize that, for any $D\in\chain_{\permA}$, one has, by \eqref{eq.marg-chain}, $\gam(D)=\sum_{\ell\in D} y_{\ell}$, which
allows one to write, for arbitrary $z\in\cor(\gam)$,
$$
\langle -\chi_{D},z\rangle \,=\, - \sum_{\ell\in D} z_{\ell} \,\leq\, -w(D) \,=\, -\sum_{\ell\in D} y_{\ell} \,=\,
\langle -\chi_{D},y\rangle\,.
$$
This basically means that $-\chi_{D}\in\nor_{\cors(\gam)}(y)$. The definition of the core yields $\chi_{N}\in\nor_{\cors(\gam)}(y)$ and the application of (the first claim in) Lemma~\ref{lem.braid-cone} gives what is desired.

To verify (b) we put $z:=\mrg^{\gam}(\permB)$ and observe $z\in\ext(\cor(\gam))$ by Lemma~\ref{lem.marg-vectors}.
The implication (a) applied to $z$ and $\permB$ yields\,
$\inte(\nor^{\permB})\subseteq \nor^{\permB}\subseteq \nor_{\cors(\gam)}(z)$. Since $\inte(\nor^{\permB})$
is an open set (in the Euclidean topology on ${\dv R}^{N}$) contained in $\nor_{\cors(\gam)}(z)$ one even has\,
$\inte(\nor^{\permB})\subseteq \inte(\nor_{\cors(\gam)}(z))$. Since $y\neq z$, the latter set is disjoint with $\nor_{\cors(\gam)}(y)$, owing to a particular relation between distinct maximal cones in a fan (see Appendix~\ref{ssec.app.normal.fan}).

Abbreviate $\nor:=\nor_{\cors(\gam)}(y)$ and $S:=\mrg^{\gam}_{-1}(\{y\})$ for the rest of the proof.
Note that $S\in{\calX}^{\preced}$ by Lemma~\ref{lem.supermo-marg-2}.
Then (a) implies $\bigcup_{\permA\in S} \nor^{\permA}\subseteq \nor$. Consider again the set ${\dv R}^{N}_{\neq}$ of vectors in ${\dv R}^{N}$ with distinct components from Lemma~\ref{lem.enum-braid}.
Using (the second claim in) Lemma~\ref{lem.braid-cone}
observe that ${\dv R}^{N}_{\neq}=\bigcup_{\permA\in\enuN} \inte(\nor^{\permA})$. Hence, (b)
gives $\nor\cap {\dv R}^{N}_{\neq}\subseteq \bigcup_{\permA\in S} \inte(\nor^{\permA})$.
As $\nor=\overline{\nor\cap {\dv R}^{N}_{\neq}}$
by Lemma~\ref{lem.enum-braid}(i), by applying the topological closure to that inclusion one gets $\nor\subseteq \bigcup_{\permA\in S} \nor^{\permA}$.

We thus know $\nor=\bigcup_{\permA\in S} \nor^{\permA}$ and put $T:=S^{\rGprec}$. As $S\in{\calX}^{\preced}$, one has $S=S^{\rGprec\lGprec}=T^{\lGprec}={\cal L}(T)$ and Lemma~\ref{lem.poset-lattice}(ii) implies that $T$ is a poset.
Given $x\in\nor$, choose $\permA\in S$ with $x\in\nor^{\permA}$.
As $\permA\in S=T^{\lGprec}$, for any $(u,v)\in T$, one has $u\preceq_{\permA} v \,\Rightarrow\, x_{u}\leq x_{v}$, by \eqref{eq.enum-cone}. So, $x\in T^{\lGbraid}$ and $\nor\subseteq T^{\lGbraid}$.
Conversely, given $x\in T^{\lGbraid}\cap {\dv R}^{N}_{\neq}$,
fix unique $\permA\in\enuN$ with $x\in\inte(\nor^{\permA})$. For any
$(u,v)\in T\setminus\diag$, one has $x_{u}<x_{v}$ and $\permA_{-1}(u)<\permA_{-1}(v)$, using \eqref{eq.enum-cone}.
Hence, $(u,v)\in T \,\Rightarrow\, u\preceq_{\permA} v$ means $\permA\in T^{\lGprec}=S$, implying $x\in\nor$; we have shown $T^{\lGbraid}\cap {\dv R}^{N}_{\neq}
\subseteq\nor$.
By Lemma~\ref{lem.braid-lattice}(iii),
$T^{\lGbraid}$ is a full-dimensional polyhedral cone and, using Lemma~\ref{lem.enum-braid}(i),  $T^{\lGbraid}\subseteq\nor$. Thus, $\nor=T^{\lGbraid}$.

For the last claim we notice a geometric fact about the polytope $\cor(\gam)$. By \mbox{\cite[Corollary\,11]{SK16},} it is a generalized permutohedron in ${\dv R}^{N}$, and, by one of its equivalent definitions \mbox{\cite[Theorem\,2.3]{JR22},} all its edges are parallel to $\chi_{u}-\chi_{v}$ for distinct $u,v\in N$. This implies
\begin{equation}
\forall\, z\in\ne_{\cors(\gam)}(y) ~\, \exists\,!~
(u(z),v(z))\in N\times N\setminus\diag
~:~ z-y=k\cdot(\chi_{u(z)}-\chi_{v(z)}),~ \mbox{with $k>0$}\,.
\label{eq.paral-edges}
\end{equation}
Put $R_{\min} :=\{\, (u(z),v(z))\,\, :\ z\in\ne_{\cors(\gam)}(y)\,\}\subseteq N\times N$ and write,
by the inequality description of the normal cone at a vertex from Appendix~\ref{ssec.app.normal.fan},

\vspace*{-8mm}
\begin{eqnarray*}
\lefteqn{\nor=\nor_{\cors(\gam)}(y)= \{\, x\in {\dv R}^{N}\, :\
\forall\, z\in\ne_{\cors(\gam)}(y) \quad \langle x,y-z\rangle\geq 0\,\}}\\
&\stackrel{\eqref{eq.paral-edges}}{=}&  \{\, x\in {\dv R}^{N}\, :\
\forall\, (u,v)\in R_{\min} ~~ \langle x,\chi_{v}-\chi_{u}\rangle\geq 0\,\} =
\{\, x\in {\dv R}^{N}\, :\ \forall\, (u,v)\in R_{\min} ~~ x_{u}\leq x_{v}\,\}\,.
\end{eqnarray*}
Thus, $\nor=R_{\min}^{\,\,\lGbraid}$, and, as $T\in\calY^{\braid}$, one has $T=T^{\lGbraid\rGbraid}=\nor^{\rGbraid}=R_{\min}^{\,\,\lGbraid\rGbraid}=\tr(R_{\min}\cup\diag)$. Because $R_{\min}\cup\diag\subseteq R$, to show\/ $T=\tr(R)$ it remains to evidence $R\subseteq \tr(R_{\min}\cup\diag)$.

To verify that consider $(u,v)\in R\setminus\diag$ and $z^{\prime}\in\ext (\cor(\gam))$ with $z^{\prime}-y$ being a positive multiple of $\chi_{u}-\chi_{v}$. Then Lemma~3.6 in \cite{Zie95} says that $z^{\prime}-y$ is in the conical hull of vectors $z-y$, $z\in\ne_{\cors(\gam)}(y)$. Hence, using \eqref{eq.paral-edges} observe that
\begin{equation}
\chi_{u}-\chi_{v} \,=\, \sum_{z\in I}\,\, \beta_{z}\cdot (\chi_{u(z)}- \chi_{v(z)})
\qquad \mbox{where $I:=\ne_{\cors(\gam)}(y)$ and $\beta_{z}\geq 0$ for $z\in I$}\,.
\label{eq.conic-comb}
\end{equation}
Because $y$ is a vertex of\/ $\cor(\gam)$, there is $x\in {\dv R}^{N}$
with $\langle x,y\rangle=0$ and $\langle x,z\rangle>0$ for $z\in\cor(\gam)\setminus\{y\}$,
and \eqref{eq.paral-edges} gives $x_{u(z)}>x_{v(z)}$ for any $z\in I$.
The formula \eqref{eq.conic-comb} then allows one to construct inductively a sequence
$u=a_{1},a_{2},\ldots $ of elements of $N$ with $(a_{j},a_{j+1})=(u(z),v(z))$ for $z\in I$, which
must necessarily end with $a_{r}=v$.
The existence of such a sequence implies $(u,v)\in\tr (R_{\min})$.
\end{proof}

Note in this context 
that non-neighbouring vertices $y,z$ of\/ $\cor(\gam)$ may exist which differ in two components only:
take $N=\{ a,b,c \}$, $\gam:=2\cdot\delta_{N}+\delta_{ab}+\delta_{ac}$. Then
$\cor(\gam)$ has 4 vertices: $y=(y_{a},y_{b},y_{c})=(1,0,1)$,
its neighbours $(2,0,0)$ and $(0,1,1)$, and its counterpart $z=(1,1,0)$.

\subsection{Conditional independence interpretation}\label{ssec.CI-iterpret}
Every triplet $(A,B|C)$ of pairwise disjoint subsets of $N$ (= a triplet over $N$) can be assigned a formal
{\em conditional independence\/} (CI) statement that ``$A$ is conditionally independent of $B$ given $C$" \cite{Stu05}.
Given a supermodular game $\gam\in\supmoN$, the set of triplets 
over $N$ satisfying $\diff\gam\,(A,B|C)=0$ is interpreted as the CI structure induced by $\gam$. The next lemma characterizes the induced elementary CI statements in terms of the tightness classes (defined in Section~\ref{ssec.marginal-vectors}); recall that these are topologies
on $N$ (see Section~\ref{ssec.supermo-posets}).

\begin{lemma}\label{lem.CI-interpret}\rm
Given a supermodular game $\gam\in\supmoN$ and $(a,b|C)\in\elemN$ (see p.\,\pageref{ssec.supermodular}), one has
\begin{equation}
\diff\gam\,(a,b|C)=0 ~~\Leftrightarrow~~ [\,\,\exists\, y\in\ext(\cor(\gam))\,:~~ aC,bC\in\tight^{\gam}_{y} \,]\,.
\label{eq.CI-interpret}
\end{equation}
\end{lemma}

\begin{proof}
For the necessity we fix an enumeration $|c_{1}|\ldots |c_{k}|$ of $C$ in case $k:=|C|\geq 1$
and also an enumeration $|d_{1}|\ldots |d_{l}|$ of $D:=N\setminus abC$ in case $l:=|D|\geq 1$.
Then we put
$\permA := |c_{1}|\ldots |c_{k}|a|b|d_{1}|\ldots |d_{l}|$ and
$\permB := |c_{1}|\ldots |c_{k}|b|a|d_{1}|\ldots |d_{l}|$.
Additionally, we put $y:=\mrg^{\gam}(\permA)$ and $z:=\mrg^{\gam}(\permB)$ and the first step is to observe that $y=z$.\\[0.2ex]
Indeed, because $\chain_{\permA}$ and $\chain_{\permB}$ differ by exchange of $aC$ for $bC$, \eqref{eq.marg-vector}
gives $y_{\ell}=z_{\ell}$ for $\ell\in N\setminus ab$. As $y_{a}\stackrel{\eqref{eq.marg-vector}}{=}\gam(aC)-\gam(C)$ while $z_{a}\stackrel{\eqref{eq.marg-vector}}{=}\gam(abC)-\gam(bC)$ the assumption $\diff\gam\,(a,b|C)=0$ yields $y_{a}=z_{a}$. Analogously,
$y_{b}=\gam(abC)-\gam(aC)=\gam(bC)-\gam(C)=z_{b}$.\\[0.2ex]
Hence, $\permA,\permB\in\mrg^{\gam}_{-1}(\{y\}) ~ \stackrel{\eqref{eq.marg-chain}}{\Longrightarrow}~ \chain_{\permA},\chain_{\permB}\subseteq\tight^{\gam}_{y}$, which implies
$aC,bC\in\tight^{\gam}_{y}$.
Recall that Lemma~\ref{lem.marg-vectors} implies $y\in\ext(\cor(\gam))$.

For the sufficiency consider $y\in\ext(\cor(\gam))$ with $aC,bC\in\tight^{\gam}_{y}$.
Thus, $y\in\mrg^{\gam}(\enuN)$ by Lemma~\ref{lem.marg-vectors}, which allows one to
apply Lemma~\ref{lem.supermo-marg-1} saying that $\tight^{\gam}_{y}$ is closed under
intersection and union. This gives $C,abC\in\tight^{\gam}_{y}$ and justifies the following writing
\begin{eqnarray*}
\lefteqn{\hspace*{-25mm}\diff\gam\,(a,b|C)=\gam(abC)+\gam(C)-\gam(aC)-\gam(bC)}\\
&\stackrel{\eqref{eq.tight-class}}{=}&
\sum_{\ell\in abC} y_{\ell} +\sum_{\ell\in C} y_{\ell}
-\sum_{\ell\in aC} y_{\ell}-\sum_{\ell\in bC} y_{\ell}=0\,,
\end{eqnarray*}
concluding the sufficiency proof.
\end{proof}

\section{Main result}\label{sec.main-result}
To formulate our equivalence result we need to introduce four different combinatorial concepts
assigned to a supermodular game $\gam$ over $N$.
\begin{itemize}
\item The {\em enumeration partition\/} (or the {\em rank test\,}) induced by $\gam\in\supmoN$ is
$$
\EnPart(\gam) ~:=~ \{\, \mrg^{\gam}_{-1}(\{y\})\, :\ y\in\ext(\cor(\gam))\,\},
\quad\mbox{being a partition of the set $\enuN$},
$$
\item the {\em fan of posets\/} induced by $\gam\in\supmoN$ is
$$
\FanPos(\gam) ~:=~ \{\, (\mrg^{\gam}_{-1}(\{y\}))^{\rGprec}\, :\ y\in\ext(\cor(\gam))\,\},
~~\mbox{being a class of posets on $N$},
$$
(where the forward Galois connection $\rGprec$ corresponds to the relation $\preced$ from Section~\ref{ssec.poset-lattice})
\item the {\em tightness structure\/} (or the {\em core structure\/}) induced by $\gam$ is
$$
\TiStr(\gam) ~:=~ \{\, \tight^{\gam}_{y}\, :\ y\in\ext(\cor(\gam))\,\},
\quad\mbox{being a covering of the power set $\caP$},
$$
\item and the (conditional) {\em independence structure\/} induced by $\gam$ is
$$
\InStr(\gam) ~:=~ \{\,\, (a,b|C)\in\elemN\, :\   \diff\gam\, (a,b|C)=0\,\},
\quad\mbox{being a subset of $\elemN$}.
$$
\end{itemize}
The fifth option is a particular subgraph of the permutohedral graph
from Section~\ref{sssec.poset-lattice-characterization}. This subgraph
has the same set $\enuN$ of nodes as the full permutohedral graph
and its edges are determined by $\gam$ through the independence structure $\InStr(\gam)$.
Specifically,
\begin{itemize}
\item the {\em permutohedral subgraph\/} induced by $\gam$ is 
as follows: if $\permA,\permB\in\enuN$ only differ 
in the $i$-th two positions, $i\in[n-1]$, as described in \eqref{eq.adjacent-transpos}, then
$$
\PerSG(\gam) ~:=~ ~~\permA \leftrightarrow \permB ~~\mbox{are adjacent in $\PerSG(\gam)$~~ iff}~~
\diff\gam\, (\,\permA(i),\permA(i+1)|\bigcup_{j<i} \permA(j)\,)=0\,.
$$
\end{itemize}
The permutohedral subgraph has a close relation to the enumeration partition.

\begin{lemma}\label{lem.permut-subgraph}\rm
Given a supermodular game $\gam\in\supmoN$,
the connected components of the permutohedral subgraph $\PerSG(\gam)$ are precisely
the elements of the enumeration partition
$\EnPart(\gam)$.\\
The induced subgraphs of $\PerSG(\gam)$ for its components
coincide with the respective induced subgraphs of the whole permutohedral graph.
\end{lemma}

\begin{proof}
Assume $n\geq 2$ as otherwise the claims are trivial.
The first observation is that if $\permA,\permB\in\enuN$ satisfy
\eqref{eq.adjacent-transpos} then they are adjacent in $\PerSG(\gam)$ iff
$\mrg^{\gam}(\permA)=\mrg^{\gam}(\permB)$, which fact follows easily from
the definition of $\mrg^{\gam}$ in \eqref{eq.marg-vector} (as in the proof of Lemma~\ref{lem.CI-interpret}).  This implies that the mapping $\permA\mapsto\mrg^{\gam}(\permA)$ is constant on connectivity components of $\PerSG(\gam)$.

To show that components are distinguished by different values of $\mrg^{\gam}$ consider $\permA,\permB\in\enuN$ satisfying $\mrg^{\gam}(\permA)=\mrg^{\gam}(\permB)$. Take $y:=\mrg^{\gam}(\permA)$ and $S:=\mrg^{\gam}_{-1}(\{y\})$.
By Lemma~\ref{lem.supermo-marg-2} one has $S\in{\calX}^{\preced}$ and, by Theorem~\ref{thm.poset-characterization}, $S$ is geodetically convex in the permutohedral graph. Since $\permA,\permB\in S$, there is a geodesic between them which is fully in $S$. In particular, every consecutive pair of nodes on this geodesic is adjacent in $\PerSG(\gam)$, by the first observation. This implies that $\permA$ and $\permB$
belong to the same connectivity component of $\PerSG(\gam)$.

The first claim in Lemma~\ref{lem.permut-subgraph} thus follows from Lemma~\ref{lem.marg-vectors} saying $\ext(\cor(\gam))=\mrg^{\gam}(\enuN)$.
As concern the second claim, by the first observation, adjacent enumerations in the whole permutohedral graph belonging to the same component of $\PerSG(\gam)$  are adjacent in $\PerSG(\gam)$.
\end{proof}

\subsection{Characterization of inclusion of faces}\label{ssec.faces-inclusion}

Now, the main equivalence result can be formulated. Recall that $\face_{\supmoN}(\gam)$ denotes the {\em face of the cone\/} $\supmoN$ of {\em supermodular games\/} generated by $\gam$
(see Appendix~\ref{sec.app.polyhedral} for definitions of the geometric concepts used below).

\begin{thm}\label{thm.main}\rm
The following conditions on $\gamA,\gamB\in\supmoN$ are equivalent:
\begin{itemize}
\item[(i)] $\face_{\supmoN}(\gamA)\subseteq\face_{\supmoN}(\gamB)$,
\item[(ii)] either $m=r$ or $\exists\, \gamC\in\supmoN ~:~~ \gamB\in\,]\gamA,\gamC[\,$,
\item[(iii)] $\exists\, \gamC\in\supmoN ~ \exists\, \alpha\in (0,1)
 ~:~~ \cor(\gamB) \,=\, (1-\alpha)\cdot\cor(\gamA)\oplus \alpha\cdot\cor(\gamC)$,
\item[(iv)] the normal fan $\norm_{\cors(\gamB)}$ refines the normal fan $\norm_{\cors(\gamA)}$,
\item[(v)] the enumeration partition $\EnPart(\gamB)$ refines $\EnPart(\gamA)$,\, i.e.\,
$\forall\,S\in\EnPart(\gamB) ~\,\exists\,T\in\EnPart(\gamA):~ S\subseteq T$,
\item[(vi)] the tightness structure $\TiStr(\gamB)$ refines $\TiStr(\gamA)$,\, i.e.\,
$\forall\,\tight\in\TiStr(\gamB) ~\,\exists\,\tightB\in\TiStr(\gamA):~
\tight\subseteq\tightB$,
\item[(vii)] $\InStr(\gamB)\subseteq \InStr(\gamA)$, that is,
$\diff\gamB\, (a,b|C)=0 \,\Rightarrow\, \diff\gamA\, (a,b|C)=0$ for any $(a,b|C)\in\elemN$,
\item[(viii)] the permutohedral subgraph $\PerSG(\gamB)$ is a subgraph of $\PerSG(\gamA)$,
\item[(ix)] the fan of posets $\FanPos(\gamA)$ is sparser than $\FanPos(\gamB)$,\, i.e.\,
$\forall\,R\in\FanPos(\gamB) ~\,\exists\,Q\in\FanPos(\gamA):~ Q\subseteq R$.
\end{itemize}
\end{thm}

Note in this context that the conditions (i) and (ii) are geometric and concern the space ${\dv R}^{\caP}$ while two other geometric conditions (iii) and (iv) concern the space ${\dv R}^{N}$. On the other hand, the conditions (v)-(ix) are combinatorial ones, (viii) being specifically graphical.

\begin{proof}
Without loss of generality assume $n\geq 2$, as otherwise all conditions
(i)-(ix) trivially hold.
This is the proof scheme: (i)$\Rightarrow$(ii)$\Rightarrow$\ldots $\Rightarrow$(vii)$\Rightarrow$(i) and
(vii)$\Rightarrow$(viii)$\Rightarrow$(v)$\Leftrightarrow$(ix).
\smallskip

To derive the implication (i)$\Rightarrow$(ii) one can use the explicit characterization of vector-generated faces of a polyhedral cone ${\mathsf C}$ mentioned in Appendix~\ref{ssec.app.face-lattice}.
Recall from Section~\ref{ssec.supermodular} that the facet-defining inequalities for $\supmoN$
correspond to elementary triplets over $N$:
$$
\supmoN \,=\, \{\,\gam\in {\dv R}^{\caP}\,:\
\gam(\emptyset)=0 ~~\&~~
\diff\gam\,(a,b|C)\geq 0\quad \mbox{for any $(a,b|C)\in\elemN$}\,\}\,.
$$
The condition (i) implies  $\gamA\in\face_{\supmoN}(\gamB)$,
which means (see Appendix~\ref{ssec.app.face-lattice})
that
\begin{equation}
\forall\, (a,b|C)\in\elemN\qquad  \diff\gamB\,(a,b|C)=0 ~~\Rightarrow~~ \diff\gamA\,(a,b|C)=0\,.
\label{eq.aux-main}
\end{equation}
Let us put $\gamC_{\beta}:= (1-\beta)\cdot\gamA +\beta\cdot\gamB\in {\dv R}^{\caP}$ for every $\beta\geq 0$. Assuming $\gamA\neq\gamB$, this
is a ray in the space of games over $N$. As $\gamB$ is supermodular, one has
$\diff\gamB\,(a,b|C)\geq 0$ for any $(a,b|C)\in\elemN$. In case $\diff\gamB\,(a,b|C)=0$ one has
$\diff\gamC_{\beta}\,(a,b|C)=0$ for any $\beta\geq 0$ by \eqref{eq.aux-main}.
If $\diff\gamB\,(a,b|C)>0$ then, by continuity argument, $\diff\gamC_{\beta}\,(a,b|C)>0$ for $\beta$
close to $1$. Because $\elemN$ is finite, there exists $\beta>1$ with $\diff\gamC_{\beta}\,(a,b|C)>0$
for all $(a,b|C)\in\elemN$ with $\diff\gamB\,(a,b|C)>0$. Thus, one has $\gamC_{\beta}\in\supmoN$
for such $\beta>1$, which yields $\gamB\in\,]\,\gamA,\gamC_{\beta}\,[$, implying (ii).

For (ii)$\Rightarrow$(iii) in case $\gamA=\gamB$ put $\gamC:=\gamA$ and $\alpha:=\frac{1}{2}$. Otherwise
assume $\gamB=(1-\alpha)\cdot\gamA+\alpha\cdot\gamC$
with $\alpha\in (0,1)$ and $\gamA,\gamC\in\supmoN$. The definition of the core implies
$$
(1-\alpha)\cdot\cor(\gamA)\oplus\alpha\cdot\cor(\gamC)\subseteq \cor((1-\alpha)\cdot\gamA +
\alpha\cdot\gamC)=\cor(\gamB)\,.
$$
For the other inclusion use Lemma~\ref{lem.marg-vectors}:
given $y\in\ext(\cor(\gamB))$ there exists $\permA\in\enuN$ with $y=\mrg^{\gamB}(\permA)$ and by the linearity of the marginal-vector mapping defined in \eqref{eq.marg-vector}, we obtain
$y=(1-\alpha)\cdot\mrg^{\gamA}(\permA)+\alpha\cdot\mrg^{\gamC}(\permA)$, where $\mrg^{\gamA}(\permA)\in\cor(\gamA)$ and $\mrg^{\gamC}(\permA)\in\cor(\gamC)$.
Therefore, $y\in (1-\alpha)\cdot\cor(\gamA)\oplus\alpha\cdot\cor(\gamC)$.
The fact $\ext(\cor(\gamB))\subseteq(1-\alpha)\cdot\cor(\gamA)\oplus\alpha\cdot\cor(\gamC)$
implies, by Krein-Milman theorem, the other inclusion
$\cor(\gamB)\subseteq (1-\alpha)\cdot\cor(\gamA)\oplus\alpha\cdot\cor(\gamC)$
because the considered Minkowski sum is a convex set.

The implication (iii)$\Rightarrow$(iv) follows from a general claim about
polytopes (see Appendix~\ref{ssec.app.normal.fan}): if ${\mathsf P}={\mathsf Q}\oplus{\mathsf R}$ for non-empty polytopes ${\mathsf P},{\mathsf Q},{\mathsf R}\subseteq {\dv R}^{N}$ then
$\norm_{{\mathsf P}}$ refines $\norm_{{\mathsf Q}}$. Note that one has
$\norm_{\alpha\cdot{\mathsf Q}}=\norm_{{\mathsf Q}}$ for any positive scalar $\alpha>0$
and a polytope ${\mathsf Q}\subseteq {\dv R}^{N}$.

As concerns the implication (iv)$\Rightarrow$(v), the assumption (iv) means
$$
\forall\, y\in\ext(\cor(\gamB)) ~~\exists\, z\in\ext(\cor(\gamA)) ~\,:\quad
\nor_{\cors(\gamB)}(y)\subseteq \nor_{\cors(\gamA)}(z)\,.
$$
We apply Lemma~\ref{lem.refine-fan} both to $y$ with $\gamB$ and $z$ with $\gamA$ to get
$\bigcup_{\permA\in\mrg^{\gamB}_{-1}(\{y\})} \nor^{\permA}\subseteq \bigcup_{\permB\in\mrg^{\gamA}_{-1}(\{z\})} \nor^{\permB}$. Since distinct maximal cones $\nor^{\permA},\nor^{\permB}$
of the braid fan have the relation\, $\inte(\nor^{\permA})\cap\nor^{\permB}=\emptyset$
this necessitates
$\mrg^{\gamB}_{-1}(\{y\})\subseteq \mrg^{\gamA}_{-1}(\{z\})$; use
Lemma~\ref{lem.braid-cone}.

As concerns (v)$\Rightarrow$(vi), the assumption (v) says
$$
\forall\, y\in\ext(\cor(\gamB)) ~~\exists\, z\in\ext(\cor(\gamA)) ~\,:\quad
\mrg^{\gamB}_{-1}(\{y\})\subseteq \mrg^{\gamA}_{-1}(\{z\})\,.
$$
By combining Lemma~\ref{lem.marg-vectors} with (the second claim in) Lemma~\ref{lem.supermo-marg-1} one easily observes the inclusion
$\tight^{\gamB}_{y}  =\bigcup_{\permA\in\mrg^{\gamB}_{-1}(\{y\})} \chain_{\permA}
\subseteq \bigcup_{\permB\in\mrg^{\gamA}_{-1}(\{z\})} \chain_{\permB} =\tight^{\gamA}_{z}$.
Thus, $\tight^{\gamB}_{y} \subseteq \tight^{\gamA}_{z}$, that is,
$\TiStr(\gamB)$ refines $\TiStr(\gamA)$.

As concerns the implication (vi)$\Rightarrow$(vii), the assumption (vi) means
$$
\forall\, y\in\ext(\cor(\gamB)) ~~\exists\, z\in\ext(\cor(\gamA)) ~\,:\quad
\tight^{\gamB}_{y} \subseteq \tight^{\gamA}_{z}\,.
$$
Given $(a,b|C)\in\elemN$, the equality $\diff\gamB\,(a,b|C)=0$ implies,
by Lemma~\ref{lem.CI-interpret}, the existence of $y\in\ext(\cor(\gamB))$
with $aC,bC\in\tight^{\gamB}_{y}$. Hence, there exists a vector $z\in\ext(\cor(\gamA))$ such that
one has $aC,bC\in\tight^{\gamB}_{y}\subseteq\tight^{\gamA}_{z}$ and, again by
Lemma~\ref{lem.CI-interpret}, one gets $\diff\gamA\,(a,b|C)=0$.
In particular, $\diff\gamB\,(a,b|C)=0 ~\Rightarrow~ \diff\gamA\,(a,b|C)=0$,
which is another way of saying $\InStr(\gamB)\subseteq \InStr(\gamA)$.

The implication (vii)$\Rightarrow$(i) also follows from the characterization of facets of\/
$\supmoN$ which is recalled in Section~\ref{ssec.supermodular}.
They all have the form
$$
\face(a,b|C) \,:=\, \{\,\gamC\in\supmoN\,:\
\diff\gamC\,(a,b|C)=0\,\}\qquad \mbox{for some $(a,b|C)\in\elemN$}\,.
$$
As the face generated by a vector is the intersection of facets  (Appendix~\ref{ssec.app.face-lattice}) one has
\begin{eqnarray*}
\lefteqn{\hspace*{-5mm}\face_{\supmoN}(\gamA) \,=\, \bigcap\, \{\,\face(a,b|C):\, \diff\gamA\,(a,b|C)=0\,\}}\\[0.2ex]
&=&  \{\,\gamC\in\supmoN\,:\ \forall\, (a,b|C)\in\elemN\quad  \diff\gamA\,(a,b|C)=0 ~\Rightarrow~ \diff\gamC\,(a,b|C)=0\,\}\\[0.2ex]
&\subseteq& \bigcap\, \{\,\face(a,b|C):\, \diff\gamB\,(a,b|C)=0\,\} \,=\, \face_{\supmoN}(\gamB)\,,
\end{eqnarray*}
where the middle inclusion follows from the assumption (vii).

The implication (vii)$\Rightarrow$(viii) is immediate from the definition of the permutohedral subgraph, while
the implication (viii)$\Rightarrow$(v) follows from Lemma~\ref{lem.permut-subgraph}.
For the equivalence (v)$\Leftrightarrow$(ix) realize that (ix) says
$$
\forall\, y\in\ext(\cor(\gamB)) ~~\exists\, z\in\ext(\cor(\gamA)) ~\,:\quad
(\mrg^{\gamA}_{-1}(\{z\}))^{\rGprec}\subseteq  (\mrg^{\gamB}_{-1}(\{y\}))^{\rGprec}\,.
$$
Thus, the equivalence follows from the basic (anti-isomorphic) properties of the respective Galois connections, in this case
$S\in{\calX}^{\preced}\mapsto S^{\rGprec}\in {\calY}^{\preced}$ and its inverse
$T\in{\calY}^{\preced}\mapsto T^{\lGprec}\in {\calX}^{\preced}$. Note in this context that
$\mrg^{\gamB}_{-1}(\{y\}),\mrg^{\gamA}_{-1}(\{z\})\in {\calX}^{\preced}$ by Lemma~\ref{lem.supermo-marg-2}.
\end{proof}

\subsection{Combinatorial and geometric description of supermodular faces}\label{ssec.equiv-char}
It follows from our result that all the combinatorial
objects defined in the beginning of Section~\ref{sec.main-result} together with the geometric object of a normal fan (of the core polytope) are in one-to-one correspondence
with non-empty faces of the supermodular cone $\supmoN$.

\begin{corol}\label{cor.main}\rm
Given a supermodular game $\gam\in\supmoN$, the following
mathematical objects uniquely characterize
the face $\face_{\supmoN}(\gam)$ of the supermodular cone:
\begin{itemize}
\item[(a)] the independence structure\, $\InStr(\gam)$,
\item[(b)] the subgraph $\PerSG(\gam)$ of the permutohedral graph,
\item[(c)] the partition $\EnPart(\gam)$ of\/ $\enuN$,
\item[(d)] the fan of posets $\FanPos(\gam)$ induced by $\gam$,
\item[(e)] the tightness structure $\TiStr(\gam)$ induced by $\gam$,
\item[(f)] the normal fan $\norm_{\cors(\gam)}$ of the corresponding
core polytope $\cor(\gam)$.
\end{itemize}
\end{corol}

\begin{proof}
Given supermodular games $\gamA,\gamB\in\supmoN$, one easily observes, by double application of Theorem~\ref{thm.main}, that the following conditions are equivalent:
\begin{description}
\item[{\rm ~~(i)~~}] ~~$\face_{\supmoN}(\gamB)=\face_{\supmoN}(\gamA)$,
\item[{\rm (a)$\equiv$(vii)}] ~~$\InStr(\gamB)=\InStr(\gamA)$,
\item[{\rm (b)$\equiv$(viii)}] the permutohedral subgraph
$\PerSG(\gamB)$ coincides with the subgraph $\PerSG(\gamA)$,
\item[{\rm (c)$\equiv$(v)}] the enumeration partition
$\EnPart(\gamB)$ coincides with the partition $\EnPart(\gamA)$,
\item[{\rm (d)$\equiv$(ix)}] the fan of posets $\FanPos(\gamA)$
coincides with the fan of posets $\FanPos(\gamB)$,
\item[{\rm (e)$\equiv$(vi)}] the tightness structure
$\TiStr(\gamB)$ coincides with the tightness structure $\TiStr(\gamA)$,
\item[{\rm (f)$\equiv$(iv)}] the normal fan
$\norm_{\cors(\gamB)}$ coincides with the normal fan $\norm_{\cors(\gamA)}$.
\end{description}
Indeed, realize that all the considered binary relations of objects are anti-symmetric. This is easy in cases (a)-(c). In case (d) note that,
for every $\gam\in\supmoN$, posets in $\FanPos(\gam)$ for different
$y\in\ext(\cor(\gam))$ are inclusively incomparable, for otherwise the
respective sets of linear extensions are in inclusion (use the Galois connections from Section~\ref{ssec.poset-lattice}). Analogous observations in cases of
(e) and (f) follow by applying Galois connections from Sections~\ref{sssec.preposet-topologies} and \ref{sssec.preposet-cones}.
\end{proof}

Note is this context that Theorem~\ref{thm.main} allows one to conclude analogously that (non-empty) supermodular
faces also uniquely correspond to the following geometric objects:
\begin{itemize}
\item an equivalence class of supermodular games over $N$, where games $\gamA$ and $\gamB$ are equivalent iff there exist $\gamC_{1},\gamC_{2}\in\supmoN$ such that
$\gamA,\gamB\in\, ]\gamC_{1},\gamC_{2}[$\,,
\item an equivalence class of generalized permutohedra in
${\dv R}^{N}$, where two polytopes are equivalent iff they are each others weak Minkowski summands
(= summands of dilates).
\end{itemize}

\subsection{Dead ends on the way to a combinatorial description}\label{ssec.dead-ends}
The central theme of this treatise is the description (and intended later characterization) of atomic (non-empty)
faces of the supermodular cone in combinatorial terms, which problem is equivalent
to the task of combinatorial description of the extreme rays of the 
polymatroidal cone. The first attempts in this direction led to the idea of classifying these faces by means of
special {\em decompositions\/} (= partitions) of the {\em power set of\/ $N$}, that is, by means of special {\em equivalences on $\caP$}. This was based on the observation that, for small $|N|$, the canonical representatives of these atomic faces (= standardized versions) take on only a few distinct values and, thus, yield regular and simple decompositions of\/ $\caP$.

More specifically, any (level) equivalence on $\caP$ discussed in \cite[\S\,9.1.2]{Stu05} is
induced by a partition $\wp=\{A_{1},\ldots,A_{\ell}\}$, $\ell\geq 1$, of a subset $M\subseteq N$
in the following way: given $S,T\subseteq N$, one has $S\sim_{\wp}T$ iff $|S\cap A_{i}|=|T\cap A_{i}|$ for
$i=1,\ldots , \ell$. The same type of equivalences on $\caP$ was considered in the context of polymatroids in
\cite{CY16}, where the authors associated a certain group action on $\caP$ with $\wp$ 
and interpreted the equivalence classes of\/ $\sim_{\wp}$ as its orbits. These decompositions of\/ $\caP$
offer an elegant classification/interpretation of atomic faces despite they do not
allow one to distinguish between them: for example, if $|N|=3$ and $\wp:=\{ N\}$ then
$\sim_{\wp}$ corresponds to two different atomic faces. The problem is that for $|N|=5$
the number of available $\wp$-based equivalences on $\caP$ is 203 while the number of faces is 117978.
Therefore, this approach gives too rough classification of the faces.

On the other hand, similar equivalences on the set $\caP$ can be associated with extreme polymatroids over $N$ through their {\em lattices of flats\/} and these equivalences differentiate them in case $|N|\leq 4$.
More specifically, any (rank function of a) polymatroid $h:\caP\to {\dv R}$
defines a closure operation on subsets of $N$: given $S\subseteq N$, its closure is the set
$$
\cl_{\!h}(S):=\{\, t\in N\,:\ h(S\cup t)=h(S)\,\}\,.
$$
The respective equivalence on $\caP$ is then as follows: given $S,T\subseteq N$, one has $S\sim_{h}T$ iff
$\cl_{\!h}(S)=\cl_{\!h}(T)$. Nevertheless, this does not work in general because there exist
two distinct extreme polymatroids over $N$, $|N|=5$, inducing the same lattice of flats and, thus,
the same corresponding equivalence on $\caP$.  L\'{a}szl\'{o} Csirmaz kindly
gave an example:
take $N:=\{a,b,c,d,e\}$, the function $h:\caP\to {\dv Z}$ given by
$$
h(S) ~:=~
\left\{
\begin{array}{cl}
0 & \mbox{if $S=\emptyset$,}\\
2 & \mbox{if $S=\{a\}$,}\\
4 & \mbox{for $S\subseteq N$ with $|S|=1$, $S\neq\{a\}$,}\\
5 & \mbox{if $S=\{a,b\}$,}\\
6 & \mbox{for $S\subseteq N$ with $|S|=2$, $S\neq\{a,b\}$, $S\neq\{b,c\}$,}\\
7 & \mbox{for $S\in\{\,\, \{a,b,d\}, \{a,b,e\}, \{a,c,e\}, \{a,d,e\}\,\,\}$,}\\
8 & \mbox{for remaining $S\subseteq N$, including $S=\{b,c\}$,}
\end{array}
\right.
$$
and $h^{\prime}:\caP\to {\dv Z}$ given by $h^{\prime}(\{a\}):=3$
and $h^{\prime}(S):=h(S)$ for remaining $S\subseteq N$. Then
$h$ and $h^{\prime}$ yield distinct
extreme polymatroidal rays inducing the same lattice of flats.

\section{Conclusions: plans for extensions}\label{sec.conclusions}
In this paper, combinatorial representatives for the faces of the supermodular
cone (in ${\dv R}^{\caP}$) were introduced. The main result is Corollary~\ref{cor.main},
which offers 5 different combinatorial representatives of such a face and 1 geometric one
(in a less-dimensional space ${\dv R}^{N}$).

Despite the main objective of this long manuscript has been met, the author feels that the mission
is not complete and plans to expand this report later with additional results. The nearest goal is, in cooperation with co-workers, to provide a web platform offering computational transfers between different combinatorial representatives from Corollary~\ref{cor.main} (for $|N|\leq 6$).

Another plan is to come up with a complete proof of the claim from \cite{MPSSW09}
saying that fans of braid cones in ${\dv R}^{N}$ correspond to certain
independence structures over $N$, known as {\em semi-graphoids} \cite{Pea88}.\footnote{Let us clarify that the arguments in \cite{MPSSW09} supporting this particular claim are incomplete, so it is still possible that the claim under discussion is not true.}
If this statement is confirmed then one can establish a wider class of combinatorial representatives of fans of braid cones (in style indicated by Corollary~\ref{cor.main}) and define a sensible concept of {\em dimension\/} for any semi-graphoid. This dimension would coincide, in case of the independence structure\, $\InStr(\gam)$
induced by a standardized supermodular game $\gam$ over $N$, with the dimension of the face of the standardized supermodular cone generated by $\gam$.

Perhaps one can also offer a (linear polyhedral-geometric) criterion to recognize (strong) consistency of a combinatorial representative, by which is meant the existence of a supermodular game $\gam\in\supmoN$
inducing the combinatorial representative in sense of Corollary~\ref{cor.main}.

\appendix

\section{Posets, tosets, preposets, and lattices}\label{sec.app.lattice}
A partially ordered set $(\latt , \preceq)$, briefly a {\em poset}, is a
non-empty set $\latt$ endowed with a partial order (= ordering), that is,
a binary relation $\preceq$ on $\latt$ which is reflexive, anti-symmetric, and transitive.
The symbol $l\prec u$ for $l,u\in \latt$ will denote $l\preceq u$ with $l\neq u$.
A total order is a partial order in which every two elements are
comparable: $\forall\, l,u\in \latt$ either $u\preceq l$ or $l\preceq u$.
We will use the abbreviation {\em toset\/} for a totally ordered set.

A {\em preposet}, also {\em preorder\/} or {\em quasi-order}, is a non-empty set $\latt$ endowed with a reflexive and transitive relation $R\subseteq \latt\times\latt$.
The {\em opposite\/} of $R$ is then $R^{op}:=\{\, (u,l)\,:\, (l,u)\in R\}$, being also a preposet.
A preposet which is additionally symmetric is an {\em equivalence\/} (on $\latt$). An elementary
fact is that, given a preposet $R$, one can consider the equivalence $E:=R\cap R^{op}$ and interpret
$R$ as a partial order $\preceq$ on the set ${\cal E}$ of equivalence classes of $E$, defined by
$$
A\preceq B\quad \mbox{for $A,B\in {\cal E}$}
\quad:=\quad
[\,\exists\, l\in A~~\exists\, u\in B \,:\ (l,u)\in R\,]\,.
$$
Indeed, then $(l,u)\in R$ iff $A\preceq B$ for $A,B\in {\cal E}$ determined by $l\in A$ and $u\in B$.

A partially ordered set $(\latt ,\preceq)$ is called a {\em lattice\/} \cite[\S\,I.4]{Bir95} if every two-element subset of $\latt$ has both the least upper bound, also named the {\em supremum}, and the greatest lower bound, also named the {\em infimum}. A join/meet {\em semi-lattice\/} is a poset in which every two-element subset has the supremum/infimum.
A finite lattice, that is, a lattice $(\latt ,\preceq)$
with $\latt$ finite, is necessarily {\em complete}, which means that the requirement above holds for any subset of\/ $\latt$. Every complete lattice has the {\em least element},
denoted by $\mbox{\bf 0}$, which satisfies $\mbox{\bf 0}\preceq u$ for any $u\in \latt$ and the {\em greatest elements}, that is, an element $\mbox{\bf 1}\in \latt$ satisfying $l\preceq \mbox{\bf 1}$ for any $l\in \latt$.
One of standard examples of a finite lattice is the face lattice of a polytope: it is the set
of (all) its faces ordered by the inclusion relation $\subseteq$\,; for details see Appendix~\ref{ssec.app.face-lattice}.

A poset $(\latt ,\preceq)$ is (order) {\em isomorphic\/}
to a poset $(\latt^{\prime},\preceq^{\prime})$
if there is a bijective mapping $\iota$ from $\latt$ onto $\latt^{\prime}$ which preserves the order,
that is, $l\preceq u$ iff $\iota(l)\preceq^{\prime}\iota(u)$ for $l,u\in \latt$. If this is the case
for a lattice $(\latt ,\preceq)$
then both suprema and infima are preserved by the isomorphism $\iota$ and the lattices 
are considered to be identical mathematical structures.

On the other hand, $(\latt ,\preceq)$ is {\em anti-isomorphic\/} to a poset
$(\latt^{\prime},\preceq^{\prime})$ if there is a bijective mapping $\varrho$ from $\latt$ onto $\latt^{\prime}$ which reverses the ordering: for $l,u\in \latt$, one has $l\preceq u$ iff $\varrho(u)\preceq^{\prime}\varrho (l)$. Then the suprema are transformed to infima and conversely and we say that the posets/lattices $(\latt ,\preceq)$ and $(\latt^{\prime},\preceq^{\prime})$ are (order) {\em dual\/} each other.

\subsection{Atomistic, coatomistic, and graded lattices}\label{ssec.app.graded}

Given a poset $(\latt ,\preceq)$ and $l,u\in \latt$ such that $l\prec u$ and there is no $e\in \latt$ with
$l\prec e\prec u$, we will write $l\cover u$ and say that $l$ is {\em covered\/} by $u$ or, alternatively, that
$u$ {\em covers\/} $l$.

The covering relation is the basis for a common method to visualize finite posets, that is, posets $(\latt ,\preceq)$ with $\latt$ finite. In the so-called {\em Hasse diagram\/} of a finite poset
$(\latt ,\preceq)$, elements of $\latt$ are represented by ovals, possibly containing identifiers of
the elements. If $l\prec u$ then the oval for $u$ is placed higher than the oval for $l$.
If, moreover, $l\cover u$ then a segment is drawn between the ovals.
Clearly, any finite poset is defined up to isomorphism by its diagram.

The {\em atoms\/} of a lattice $(\latt ,\preceq)$ are its elements that cover the least element in the lattice, that is, $a\in \latt$ satisfying $\mbox{\bf 0}\cover a$.
Analogously, {\em coatoms\/} are the elements of the lattice covered by the greatest element, that is, elements $c\in \latt$ satisfying $c\cover \mbox{\bf 1}$.
We will call a lattice $(\latt ,\preceq)$ {\em atomistic\/} if every element of
$\latt$ is the supremum of a set of atoms (in $\latt$). Note that some authors
use the term ``atomic" instead of ``atomistic", see \cite[\S\,VIII.9]{Bir95} or \cite[\S\,2.2]{Zie95},
but other authors use the term ``atomic" to name a weaker condition on $(\latt,\preceq)$, namely
that, for every non-least element $e\neq \mbox{\bf 0}$ of the lattice
$\latt$, an atom $a$ exists such that $a\preceq e$.
Analogously, a lattice $(\latt ,\preceq)$
is {\em coatomistic\/} if every element of $\latt$ is the infimum of a set of coatoms.
It is evident that a lattice anti-isomorphic to an atomistic lattice is coatomistic and conversely.

A poset $(\latt ,\preceq)$ is called {\em graded\/} \cite[\S\,I.3]{Bir95}
if there exists a function $h:\latt\to {\dv Z}$, called a {\em height function},
which satisfies, for any $l,u\in \latt$,
\begin{itemize}
\item if $l\prec u$ then $h(l)<h(u)$,
\item if $l\cover u$ then $h(u)=h(l)+1$.
\end{itemize}
Note that, in case of a finite poset, the former condition is superfluous because it follows from the latter one.
The graded lattices can be characterized in terms of finite {\em chains\/}, which are (non-empty)
finite subsets of $\latt$ that are totally ordered; the {\em length\/} of such a chain is the number
of its elements minus 1.
Note that the maximal chains (in sense of inclusion) in a finite lattice with ${\bf 0}\neq {\bf 1}$ have necessarily the form ${\bf 0}\cover \ldots \cover {\bf 1}$.
A well-known fact is that a finite lattice is graded iff all maximal chains in it have the same length.
A formally stronger, but equivalent, is the {\em Jordan-Dedekind chain condition\/} requiring that all maximal chains between the same endpoints have the same length \cite[Lemma~1 in \S\,I.3]{Bir95}.
The function $h$ which assigns to every $u\in\latt$ the (shared) length of maximal chains between ${\bf 0}$ and $u$ is then a (standardized) height function for the lattice $(\latt ,\preceq)$. It is evident that the dual lattice to a graded lattice is a graded lattice as well.

\subsection{Finite posets viewed as directed acyclic graphs}\label{ssec.app.trans-DAG}
There are two (different) standard  ways to represent partial orders $\preceq$ on a non-empty finite set $N$
by means of directed acyclic graphs having $N$ as the set of nodes.

We say that a binary relation $R\subseteq N\times N$ is {\em acyclic\/} if there is no sequence $u_{1},u_{2},\ldots, u_{k}$, $k\geq 2$, of elements of $N$ with $u_{k}=u_{1}$ such that $(u_{i},u_{i+1})\in R$ for $i=1,\ldots k-1$. Every such a relation can be depicted in the form of a {\em directed acyclic graph\/} having $N$ as the set of nodes, where an {\em arrow\/} (= directed edge) $u\to v$ is drawn iff $(u,v)\in R$. A well-known equivalent definition (of a directed acyclic graph over $N$) is that all elements of $N$ can be ordered into a sequence $u_{1},\ldots,u_{n}$, $n=|N|$, such that $u_{i}\to u_{j}$ implies $i<j$.

The first way to represent a partial order $\preceq$ on $N$ is based on its {\em strict version $\prec$\,}.
Evidently, a binary relation $R\subseteq N\times N$ is the strict version $\prec$ of 
a partial order $\preceq$ on $N$ iff it is transitive and {\em irreflexive}, which means that there is no $u\in N$ with $(u,u)\in R$. Another equivalent condition is that $R$ is transitive and {\em asymmetric}, which means
that $(u,v)\in R$ implies $(v,u)\not\in R$.
Because a transitive relation is irreflexive iff it is acyclic, $\prec$\, can be depicted using a {\em transitive directed acyclic\/} graph $(N,\to)$, which adjective means that the relation $\to$ is transitive.
This implies that $T\subseteq N\times N$ is a poset iff the {\em diagonal\/}
$\diag:=\{\, (u,u)\,:\ u\in N\}$  for $N$  is contained in $T$ and the rest
$T\setminus\diag$ is both transitive and acyclic.

The second way comes from the {\em covering relation $\cover$} induced by the respective partial order $\preceq$ on $N$. Since $\cover$ is a sub-relation of $\prec$, it is also acyclic and can be represented by a directed acyclic graph $(N,\hookrightarrow)$ where $u\hookrightarrow v$ iff $u\cover v$.
This {\em directed Hasse diagram\/} of\, $\preceq$ is
(strongly) {\em anti-transitive\/} in sense that $u_{1}\hookrightarrow u_{2}\,\ldots \hookrightarrow u_{k}$, $k\geq 3$, implies $\neg [\,u_{1}\hookrightarrow u_{k}\,]$. In fact,
it is the sparsest graph $(N,\hookrightarrow)$ such that the transitive closure of\/ $\hookrightarrow$ is $\prec$.

\subsection{Galois connections for representation of a complete lattice}\label{ssec.app.Galois}
There is a general method for generating examples of complete lattices
in which method a lattice together with its dual lattice is defined on basis
of a given binary relation between two sets. The method is universal
in sense that every complete lattice is isomorphic with a lattice
defined in this way. Moreover, every finite lattice is isomorphic to
a lattice defined on basis of a given binary relation between elements of two finite sets.

This method uses the so-called {\em Galois connections}, also called {\em polarities} \cite[\S\,V.7]{Bir95}, and allows one to represent every finite lattice in the form of a {\em concept lattice}, which notion has an elegant
interpretation from the point of view of formal ontology \cite{GW99}.
To present the method we first need to recall a couple of relevant notions.

Given a non-empty set $X$, its power set will be denoted by ${\cal P}(X):=\{\, S\,:\ S\subseteq X\,\}$.
A {\em Moore family\/} of subsets of $X$, also called a {\em closure system\/} (for $X$), is a collection ${\calF}\subseteq {\cal P}(X)$ of subsets of $X$ which is closed under (arbitrary) set intersection and includes the set $X$ itself\/: $\calD\subseteq {\calF} ~\Rightarrow ~ \bigcap\calD\in {\calF}$ with a convention $\bigcap\calD=X$ for $\calD=\emptyset$. An elementary observation is that any
Moore family $({\calF},\subseteq)$, equipped with the set inclusion relation as a partial order,
is a complete lattice \cite[Theorem\,2 in \S\,V.1]{Bir95}.

A {\em closure operation\/} on subsets of $X$ is
a mapping $\cl :{\cal P}(X)\to {\cal P}(X)$ which is {\em extensive}, which means $S\subseteq \cl (S)$ for $S\subseteq X$, {\em isotone}, which means
$S\subseteq T\subseteq X \Rightarrow \cl (S)\subseteq \cl (T)$, and {\em idempotent}, which means $\cl (\cl (S))=\cl (S)$ for $S\subseteq X$. A set $S\subseteq X$ is then called {\em closed with respect to $\cl$} if $S=\cl (S)$.
The basic fact is that the collection of subsets of $X$ closed
with respect to $\cl$ is a Moore family and any Moore family can be defined in this way \cite[Theorem\,1 in \S\,V.1]{Bir95}. Specifically, given a Moore family $\calF$ of subsets of $X$, the formula
$$
\cl_{{\calF}}(S) \,:=\, \bigcap\, \{\,F\subseteq X\,:\, S\subseteq F\in
{\calF}\,\} \quad \mbox{\rm for} ~ S\subseteq X\,,
$$
defines a closure operation on subsets of $X$ having ${\calF}$ as
the collection of closed sets with respect to $\cl_{{\calF}}$;
see \cite[Theorem\,1]{GW99}.

In other words, there is a one-to-one correspondence between closure operations on subsets of $X$ and Moore families of subsets of $X$, which are
examples of complete lattices relative to set inclusion relation $\subseteq$.
Thus, any closure operation $\cl$ on (subsets of) $X$ yields an example of
a complete lattice $({\calF},\subseteq)$, where ${\calF}$ is the collection
of sets closed with respect to $\cl$.
\smallskip

We now present the method of Galois connections itself. Let $X$ and $Y$ be non-empty
sets and $\mathfrak{r}$ a binary relation between elements of $X$ and $Y$: $\mathfrak{r}\subseteq X\times Y$.
We are going to write $x\,\mathfrak{r}\,y$ to denote that $x\in X$ and $y\in Y$ are in this
{\em incidence\/} relation: $(x,y)\in \mathfrak{r}$. Then every subset $S$ of $X$ can be assigned
a subset of $Y$ as follows (the {\em forward\/} direction):
$$
S\subseteq X \quad\mapsto\quad S^{\rGalo} \,:=\,  \{\, y\in Y\,:\ x\,\mathfrak{r}\,y \quad
\mbox{for every $x\in S$}\,\}~ \subseteq Y\,.
$$
Analogously, every subset $T$ of\/ $Y$ can be assigned a subset of $X$ (the {\em backward\/} direction):
$$
T\subseteq Y \quad\mapsto\quad T^{\lGalo} \,:=\,  \{\, x\in X\,:\ x\,\mathfrak{r}\,y \quad
\mbox{for every $y\in T$}\,\}~ \subseteq X\,.
$$
The mappings $S\mapsto S^{\rGalo}$ and $T\mapsto T^{\lGalo}$ are then
{\em Galois connections\/} between the lattices $({\cal P}(X),\subseteq)$ and $({\cal P}(Y),\subseteq)$, also
called {\em polarities\/} \cite[\S\,V.7]{Bir95}. The point is that the mapping $S\mapsto
S^{\rGalo\lGalo}:=(S^{\rGalo})^{\lGalo}$ is a closure operation on subsets of $X$ and the mapping $T\mapsto
T^{\lGalo\rGalo}:=(T^{\lGalo})^{\rGalo}$ is a closure operation on subsets of $Y$. Thus, they define a complete lattice
$({\calX}^{\mathfrak{r}},\subseteq)$ of closed subsets of $X$ and a complete lattice $({\calY}^{\mathfrak{r}},\subseteq)$ of
closed subsets of $Y$. In addition to that,
\begin{itemize}
\item one has $S\in {\calX}^{\mathfrak{r}}$ iff
$S=T^{\lGalo}$ for some $T\subseteq Y$,
\item analogously $T\in {\calY}^{\mathfrak{r}}$ iff $T=S^{\rGalo}$ for some $S\subseteq X$, and
\item the inverse of the mapping $S\in {\calX}^{\mathfrak{r}}\mapsto S^{\rGalo}$ is the mapping
$T\in {\calY}^{\mathfrak{r}}\mapsto T^{\lGalo}$, defining an anti-isomorphism between lattices
$({\calX}^{\mathfrak{r}},\subseteq)$ and $({\calY}^{\mathfrak{r}},\subseteq)$ \cite[Theorem\,19 in \S\,V.7]{Bir95}.
\end{itemize}
In particular, the complete lattice $({\calX}^{\mathfrak{r}},\subseteq)$ and
its dual lattice $({\calY}^{\mathfrak{r}},\subseteq)$ have together been defined on basis of a given binary relation $\mathfrak{r}\subseteq X\times Y$. Note that in case of finite $X$, or of finite $Y$, the lattice
$({\calX}^{\mathfrak{r}},\subseteq)$ is finite, too.
Useful observation is that $S^{\rGalo}=S^{\rGalo\lGalo\,\rGalo}$ for any $S\subseteq X$.

The arguments why this method is universal for representing complete lattices are based on results from \cite[\S\,V.9]{Bir95}. More specifically, given a complete lattice $(\latt ,\preceq)$, one can put $X:=\latt$,
$Y:=\latt$ and consider a binary relation $\mathfrak{r}:= \{\, (x,y)\in X\times Y\,:\  x\preceq y\}$.
Then, for any $S\subseteq X=\latt$, $S^{\rGalo}$ is the set of upper bounds for $S$ in $(\latt ,\preceq)$
while, for any $T\subseteq Y=\latt$, $T^{\lGalo}$ is the set of lower bounds for $T$. Then
\cite[Theorem\,23 in \S\,V.9]{Bir95} implies that every closed subset $S\subseteq X$, that is,
every set $S\in {\calX}^{\mathfrak{r}}$, is a principal ideal in  $(\latt ,\preceq)$, which means
it has the form $S=\{z\in \latt\,:\ z\preceq a\,\}$ for some $a\in \latt$.
Thus, $({\calX}^{\mathfrak{r}},\subseteq)$ coincides with the lattice of principal ideals in $(\latt ,\preceq)$,
which is known to be isomorphic with the lattice $(\latt ,\preceq)$ itself by \cite[Theorem\,3 in \S\,II.3]{Bir95}.
Note that in case of a finite lattice the sets $X$ and $Y$ are both finite.

As a by-product of the above construction we observe that any complete lattice is isomorphic
to a Moore family of subsets of a non-empty set $X$ (ordered by set inclusion), and, therefore,
it can be introduced by means of a closure operation on subsets of $X$.

\section{Concepts from polyhedral geometry}\label{sec.app.polyhedral}
Standard concepts and basic facts from polyhedral geometry are explained in detail in textbooks on this topic \cite{Bro83,Sch98,Zie95}. Nevertheless, for reader's convenience,
we recall those of them that are commonly used in the paper. Before reading this section one
should refresh the notational conventions introduced in the beginning of Section~\ref{sec.preliminary}.

Our geometric considerations concern real (finite-dimensional) Euclidean spaces ${\dv R}^{\indset}$, where $\indset$ is a \underline{non-empty} finite index set. There are two different contexts we consider.
The basic level is to have $\indset :=N$, where $N$ is our
(unordered) {\em basic set\/} (of variables). In this case
the index set $\indset$ has no additional mathematical structure and
the components of vectors are denoted by means of lower indices: $x=[x_{\ell}]_{\ell\in N}$.
The advanced level is to have\/ $\indset :=\caP$, in which case the index set\/ $\indset$ has a deeper mathematical structure. It is the binary relation $\subseteq$ of inclusion among subsets of $N$; thus, $\indset$ is equipped with the structure of a Boolean algebra then.
The vectors in ${\dv R}^{\caP}$ are interpreted as {\em set functions\/} $\gam :\caP\to {\dv R}$
and their components are denoted using functional arguments: $\gam(A)$ for $A\subseteq N$.
Nevertheless, the concepts and facts from polyhedral geometry are needed in both contexts.
Recall that the {\em scalar product\/} of two vectors $x,y\in {\dv R}^{\indset}$
is denoted by 
$\langle x,y\rangle$.
The {\em Minkowski sum\/} of two sets ${\mathsf S},{\mathsf T}\subseteq {\dv R}^{\indset}$
is the set ${\mathsf S}\oplus {\mathsf T}:=\{\, x+y\,:\ x\in {\mathsf S},~ y\in {\mathsf T}\,\}\subseteq {\dv R}^{\indset}$. The {\em scalar multiple\/}
of a set ${\mathsf S}\subseteq {\dv R}^{\indset}$ by $\alpha\in {\dv R}$
is $\alpha\cdot {\mathsf S}:=\{\, \alpha\cdot x\,:\ x\in {\mathsf S}\,\}$.
\smallskip

An affine combination of vectors in ${\dv R}^{\indset}$ is a real (finite) linear combination
$\sum_{j} \alpha_{j}\cdot x^{j}$, $\alpha_{j}\in {\dv R}$, $x^{j}\in {\dv R}^{\indset}$,
where $\sum_{j} \alpha_{j}=1$. The affine hull of a set ${\mathsf S}\subseteq {\dv R}^{\indset}$
is the collection of all affine combinations of vectors from ${\mathsf S}$.
It is always an {\em affine subspace\/} of\/ ${\dv R}^{\indset}$, that is,
a subset ${\mathsf A}\subseteq {\dv R}^{\indset}$ closed under affine combinations.
A non-empty affine subspace is always a translate of a linear subspace of ${\dv R}^{\indset}$,
that is, a set of the form ${\mathsf A}=\{x\}\oplus {\mathsf L}$, 
where $x\in {\dv R}^{\indset}$ and ${\mathsf L}\subseteq {\dv R}^{\indset}$ is a linear subspace
\cite[\S\,1]{Bro83}; ${\mathsf L}$ is then uniquely determined by ${\mathsf A}$ while $x$ is not unique.
Recall that a linear subspace is always non-empty by definition.
The dimension of the affine space ${\mathsf A}$ is then the dimension of ${\mathsf L}$.
This allows one to introduce the {\em dimension\/} of any set ${\mathsf S}\subseteq {\dv R}^{\indset}$
as the dimension of its affine hull. By a convention, the dimension of the empty set $\emptyset\subseteq {\dv R}^{\indset}$ is $-1$.
A {\em hyperplane\/} in ${\dv R}^{\indset}$ is an affine subspace ${\mathsf H}$ of\/ ${\dv R}^{\indset}$
of the dimension $|\indset|-1$. An equivalent condition is that it is the set of solutions
$x\in {\dv R}^{\indset}$ to the equation $\langle x,y\rangle=\beta$, where $\beta\in {\dv R}$
and $y$ is a non-zero vector in ${\dv R}^{\indset}$; see \cite[\S\,1]{Bro83}.
A set ${\mathsf S}\subseteq {\dv R}^{\indset}$ is {\em bounded\/} if
there are constants $\alpha,\beta\in {\dv R}$ such that $\alpha\leq x_{\ell}\leq\beta$
for any component $x_{\ell}$, $\ell\in\indset$, of any $x\in {\mathsf S}$.

A set ${\mathsf S}\subseteq {\dv R}^{\indset}$ is {\em convex\/} if it is closed under convex combinations.
We also consider the usual Euclidean topology on the space ${\dv R}^{\indset}$ and a set
${\mathsf S}\subseteq {\dv R}^{\indset}$ is named briefly {\em closed\/} if it is
closed under limits relative to this topology. Most of the sets considered in this treatment
are {\em closed convex\/} subsets of\/ ${\dv R}^{\indset}$.
The {\em relative interior\/} of a (non-empty) closed convex set ${\mathsf C}\subseteq{\dv R}^{\indset}$, denoted by
$\relint({\mathsf C})$, is the topological interior of ${\mathsf C}$ relative to (the induced) Euclidean
topology on the affine hull of ${\mathsf C}$. In case ${\mathsf C}$ is full-dimensional, that is, if the
affine hull of ${\mathsf C}$ is ${\dv R}^{\indset}$, it reduces to the usual topological
{\em interior}, denoted by $\inte({\mathsf C})$. A well-known fact is that any closed convex  set
${\mathsf C}\subseteq{\dv R}^{\indset}$ is the topological closure of its relative interior; use
\cite[Theorem\,3.4(c)]{Bro83}.

\subsection{Polyhedrons, polytopes, and polyhedral cones}\label{ssec.app.polyhedrons}
A {\em polyhedron\/} in ${\dv R}^{\indset}$ is a set of vectors $x\in {\dv R}^{\indset}$ specified by
finitely many linear inequalities $\langle x,y\rangle\geq \beta$, where
$y\in {\dv R}^{\indset}$ and $\beta\in {\dv R}$ are parameters of the inequalities.
It is evident that every polyhedron is a closed convex set in the above sense.

A {\em polytope\/} (in ${\dv R}^{\indset}$) is the convex hull of a finite set of vectors in ${\dv R}^{\indset}$.
Note explicitly that we consider the empty set to be a polytope (this is a point where some authors might differ in their conventions).
A fundamental result in polyhedral geometry says that a subset of\/ ${\dv R}^{\indset}$ is a polytope iff it is a {\em bounded polyhedron}; see \cite[Theorem\,1.1]{Zie95} or \cite[Theorem\,9.2]{Bro83}.

A {\em polyhedral cone\/} in\/ ${\dv R}^{\indset}$ is a subset ${\mathsf C}$ of\/ ${\dv R}^{\indset}$ defined as the conic hull of a non-empty finite set ${\mathsf S}$ of vectors from ${\dv R}^{\indset}$.
An equivalent definition is that ${\mathsf C}$ is specified by finitely many inequalities $\langle x,y\rangle\geq 0$ for $x\in {\dv R}^{\indset}$, with $y\in {\dv R}^{\indset}$ determining the inequality
\cite[Theorem\,1.3]{Zie95}. Thus, because ${\mathsf C}$ always contains the zero vector in ${\dv R}^{\indset}$,
denoted here by ${\mathsf 0}$, it is a non-empty polyhedron.
The {\em linearity space\/} of a polyhedral cone ${\mathsf C}\subseteq {\dv R}^{\indset}$ is the set
${\mathsf L}:=-{\mathsf C}\cap{\mathsf C}$, where $-{\mathsf C}:=\{-x\,:\ x\in {\mathsf C}\,\}$; it is indeed
a linear subspace of ${\dv R}^{\indset}$.

A polyhedral cone ${\mathsf C}\subseteq {\dv R}^{\indset}$ is called {\em pointed}, if $-{\mathsf C}\cap{\mathsf C}=\{{\mathsf 0}\}$.
An equivalent condition is that there exists a non-zero
$y\in {\dv R}^{\indset}$ such that
$\langle x,y\rangle>0$ for any $x\in {\mathsf C}\setminus\{{\mathsf 0}\}$; see \cite[Proposition\,2]{Stu93}.
%
It makes no problem to observe that, if this is the case then,
for each $\beta>0$, the intersection of ${\mathsf C}$ with the hyperplane ${\mathsf H}:=\{\, x\in {\dv R}^{\indset}\,:\
\langle x,y\rangle =\beta\,\}$ is a polytope. This is because then one can find finite
${\mathsf S}^{\prime}\subseteq {\mathsf H}$ with ${\mathsf C}=\cone ({\mathsf S}^{\prime}\cup \{{\mathsf 0}\})$,
where $\cone ({\mathsf S})$ denotes the conical hull of ${\mathsf S}\subseteq {\dv R}^{\indset}$; hence,
${\mathsf C}\cap{\mathsf H}$ is the convex hull of\/ ${\mathsf S}^{\prime}$.

\subsection{Faces, vertices, edges, and facets}\label{ssec.app.vertices-facets}
Given distinct $x,y\in {\dv R}^{\indset}$, the convex hull of $\{x,y\}$ is the
{\em closed segment}, denoted by $[x,y]$, while the {\em open segment\/}, denoted by $]x,y[$,
consists of convex combinations of $x$ and $y$ which have both coefficients non-zero:
$$
]x,y[ ~:=~ \{\, (1-\alpha)\cdot x+\alpha\cdot y\,:\ 0<\alpha<1\,\}\,.
$$
An {\em exposed face\/} of a closed convex set ${\mathsf C}\subseteq {\dv R}^{\indset}$ \cite[\S\,5]{Bro83}
is a subset ${\mathsf F}\subseteq {\mathsf C}$ consisting of vectors $x\in {\mathsf C}$ satisfying $\langle x,y\rangle=\beta$ for some
$y\in {\dv R}^{\indset}$, $\beta\in {\dv R}$, such that $\langle z,y\rangle\geq\beta$ is a valid inequality
for all $z\in {\mathsf C}$.
Observe that the empty set is always an exposed face of (any)
${\mathsf C}$: take $y:={\mathsf 0}$ and $\beta:=-1$.
A {\em face\/} of a closed convex set ${\mathsf C}\subseteq {\dv R}^{\indset}$ \cite[\S\,5]{Bro83}
is a convex subset ${\mathsf F}$ of ${\mathsf C}$ such that one has $[x,z]\subseteq {\mathsf F}$
whenever $x,z\in{\mathsf C}$ and $]x,z[\,\cap\, {\mathsf F}\neq\emptyset$.
Clearly, every exposed face is a face in this sense and, in case of a polyhedron ${\mathsf C}$,
these two conditions are equivalent, which implies that a face of a polyhedron is a polyhedron.
These facts follow from a particular characterization of non-empty faces of a polyhedron ${\mathsf P}\subseteq {\dv R}^{\indset}$. Specifically, consider an inequality description of a polyhedron
${\mathsf P}\,=\,\{\, x\in {\dv R}^{\indset}\,:\ \langle x,y_{i}\rangle\geq \beta_{i}~~\mbox{for $i\in I$}\,\}$
where $y_{i}\in {\dv R}^{\indset}$, $\beta_{i}\in {\dv R}$, $I$ finite, are parameters.
Then a \underline{non-empty} subset ${\mathsf F}\subseteq {\mathsf P}$ is a face of ${\mathsf P}$
iff there exists (possibly empty) $J\subseteq I$ such that
${\mathsf F}\,=\,\{\, x\in {\mathsf P}\,:\ \langle x,y_{j}\rangle =\beta_{j}~~\mbox{for $j\in J$}\,\}$;
use \cite[Theorem\,7.51]{BK92} which has an implicit assumption of non-emptiness of the set ${\mathsf F}$.

One of basic facts is that every proper face ${\mathsf F}\subset {\mathsf C}$ of a closed convex set ${\mathsf C}\subseteq {\dv R}^{\indset}$ is disjoint with its relative interior: ${\mathsf F}\cap \relint({\mathsf C})=\emptyset$  \cite[Theorem\,5.3]{Bro83}. In fact,
the complement of the union of all proper faces of ${\mathsf C}$ is precisely the relative interior $\relint({\mathsf C})$; use \cite[Corollary\,5.7]{Bro83}.

The number of faces of a polyhedron ${\mathsf P}$ is finite; see \cite[Corollary\,8.5]{Bro83}.
Faces of a polyhedron ${\mathsf P}\subseteq {\dv R}^{\indset}$ can be classified according to their dimensions.
A {\em vertex\/} of a polyhedron ${\mathsf P}$ is such a vector $y\in {\mathsf P}$ that
the singleton $\{y\}$ is a face of ${\mathsf P}$ (of the dimension $0$).
An alternative terminology is that $y$ is an {\em extreme point\/} of\/ ${\mathsf P}$,
which is a vector $y\in {\mathsf P}$ such that there is no open segment $]x,z[$ with distinct
$x,z\in {\mathsf P}$ and $y\in\,]x,z[$.
The set of vertices of a polytope ${\mathsf P}$ will be denoted by $\ext ({\mathsf P})$.
A well-known consequence of famous {\em Krein-Milman theorem\/} is that every polytope ${\mathsf P}$
has finitely many vertices and equals to the convex hull of the set of its vertices:
${\mathsf P}= \conv (\ext ({\mathsf P}))$; see \cite[Theorem\,7.2(c)]{Bro83} or \cite[Proposition\,2.2(i)]{Zie95}.

An {\em edge\/} of a polytope ${\mathsf P}$ is a closed segment $[y,z]\subseteq {\mathsf P}$ which is a face of
${\mathsf P}$ (of the dimension~$1$); then necessarily $y,z\in\ext ({\mathsf P})$.
The set of ({\em geometric\/}) {\em neighbours\/} of a vertex $y$ of ${\mathsf P}$ is
$$
\ne_{{\mathsf P}}(y) ~:=~ \{\,z\in\ext ({\mathsf P})\,:\ \mbox{$[y,z]$ is an edge of ${\mathsf P}$}\,\}\,.
$$
In case of a general polyhedron ${\mathsf P}$ other kinds of (unbounded) faces of the dimension $1$ may exist:
these could be extreme rays or lines.

Given a non-empty polyhedron ${\mathsf P}$, any of its faces of the dimension $\dim ({\mathsf P})-1$ is called a {\em facet\/} of\/ ${\mathsf P}$.
A \underline{non-empty} face ${\mathsf F}$ of ${\mathsf P}$ is its facet iff it is an inclusion-maximal proper face of\/ ${\mathsf P}$ (=
a face distinct from ${\mathsf P}$); use \cite[Theorem\,5.2 and Corollaries~5.5 and 8.6]{Bro83}.
A basic fact is that each {\em full-dimensional proper\/} polyhedron ${\mathsf P}\subset {\dv R}^{\indset}$, that is,
each polyhedron with $\dim ({\mathsf P})=|\indset|$ and ${\mathsf P}\neq {\dv R}^{\indset}$, is specified by those valid inequalities for ${\mathsf P}$
which define facets and the specification of\/ ${\mathsf P}$ by {\em facet-defining\/} inequalities is the
{\em unique\/} inclusion-minimal inequality description of ${\mathsf P}$ within ${\dv R}^{\indset}$ (up to positive multiples of inequalities); see \cite[Theorem\,8.2]{Bro83}.
In case ${\mathsf P}$ is a polytope this is what is called the (minimal) {\em polyhedral description\/} of ${\mathsf P}$.
\smallskip

Since the intersection of any collection of faces of a polyhedron ${\mathsf P}$ is a face of ${\mathsf P}$, for
any vector $y\in {\mathsf P}$, there exists the smallest face of\/ ${\mathsf P}$ containing $y$. This (non-empty) {\em face generated by $y$} will be denoted by $\face_{{\mathsf P}}(y)$ throughout this paper.
Given a vector $y\in {\mathsf P}$ and a face ${\mathsf F}\subseteq {\mathsf P}$,
one has ${\mathsf F}=\face_{{\mathsf P}}(y)$ iff $y\in\relint ({\mathsf F})$; see \cite[Proposition\,5.6]{Bro83}.
In particular, for $y,z\in {\mathsf P}$, one has $\face_{{\mathsf P}}(y)=\face_{{\mathsf P}}(z)$ iff there is a face
${\mathsf F}$ of ${\mathsf P}$ with $y,z\in\relint ({\mathsf F})$. Thus, relative interiors of non-empty
faces of\/ ${\mathsf P}$ are equivalence classes of this particular equivalence of vectors $y,z\in {\mathsf P}$
defined by $\face_{{\mathsf P}}(y)=\face_{{\mathsf P}}(z)$; compare with \cite[Corollary\,5.7]{Bro83}.
\smallskip

A (non-empty) polytope ${\mathsf P}\subseteq {\dv R}^{\indset}$ of the dimension $d\geq 0$ is called
{\em simple\/} if any of its vertices belongs precisely to $d$ facets of\/ ${\mathsf P}$.
A well-known equivalent condition is that every vertex has precisely $d$ neighbours; see \cite[Theorems~12.11 and 12.12]{Bro83}.

\subsection{Face lattice of a polytope}\label{ssec.app.face-lattice}
A fundamental fact about the collection of faces of a polytope ${\mathsf P}$, ordered by inclusion relation $\subseteq$, is
that it is a lattice which is both atomistic and coatomistic; moreover, it a graded lattice, where the dimension of a face
plays the role of a height function; see \cite[Theorem~2.7(i)(v)]{Zie95}.
That is why it is called the {\em face lattice\/} of\/ ${\mathsf P}$.

The same lattice-theoretical properties hold for the collection of non-empty faces of a polyhedral cone ${\mathsf C}\subseteq {\dv R}^{\indset}$. Since a non-empty face of a polyhedral cone is also a polyhedral cone,
it is a collection of polyhedral cones in ${\dv R}^{\indset}$.
To observe that it is order-isomorphic to a face lattice of a polytope we first realize that every polyhedral cone ${\mathsf C}$
has the form ${\mathsf C}={\mathsf L}+{\mathsf C}^{\prime}$, where ${\mathsf L}$ is the linearity space of ${\mathsf C}$ and ${\mathsf C}^{\prime}$ a pointed polyhedral cone \cite[Lemma\,7.4]{BK92}.
Given a (non-empty) face ${\mathsf F}$ of ${\mathsf C}$, the set ${\mathsf F}^{\prime}:={\mathsf F}\cap {\mathsf C}^{\prime}$
is a (non-empty) face of ${\mathsf C}^{\prime}$ such that ${\mathsf F}={\mathsf L}+{\mathsf F}^{\prime}$, which defines an order isomorphism
of the considered posets \cite[Lemma\,7.18]{BK92}. In particular, the least non-empty face of any polyhedral cone
${\mathsf C}$ is its linearity space ${\mathsf L}$.

Thus, without loss of generality, ${\mathsf C}$ can be assumed to be pointed. In that case, a hyperplane  ${\mathsf H}\subseteq {\dv R}^{\indset}$ exists such that ${\mathsf P}:={\mathsf C}\cap {\mathsf H}$ is a polytope (see Appendix~\ref{ssec.app.polyhedrons}).
Thus, given a non-empty face ${\mathsf F}$ of ${\mathsf C}$, the set ${\mathsf F}^{\prime\prime}:={\mathsf F}\cap {\mathsf H}$ is a face of ${\mathsf P}$ such that ${\mathsf F}=\cone ({\mathsf F}^{\prime\prime}\cup \{{\mathsf 0}\})$, which defines an order isomorphism between non-empty faces of ${\mathsf C}$ and (all) faces of ${\mathsf P}$.

The lattice-theoretical facts above imply that every non-empty face ${\mathsf F}$ of a polyhedral cone ${\mathsf C}\subseteq {\dv R}^{\indset}$ is the intersection of (a class of) facets of ${\mathsf C}$ containing ${\mathsf F}$ (where the intersection of the empty class is ${\mathsf C}$ by a convention).
This can further be extended as follows. Given an inequality description
${\mathsf C}\,=\,\{\, x\in {\dv R}^{\indset}\,:\ \langle x,y_{i}\rangle\geq 0~~\mbox{for $i\in I$}\,\}$
of a polyhedral cone and a vector $z\in {\mathsf C}$ in it, the
face generated by $z$ has the following explicit description:
$$
\face_{{\mathsf C}}(z) ~=~ \{\, x\in {\mathsf C}\,:\ \forall\, i\in I\quad \langle z,y_{i}\rangle =0
~\Rightarrow~ \langle x,y_{i}\rangle =0\,\,\}\,,
$$
which conclusion can be easily derived from \cite[Theorem\,7.27]{BK92}.

Further elementary observations about being-a-face relation concern its transitivity:
if, for polyhedral cones ${\mathsf C}_{1},{\mathsf C}_{2},{\mathsf C}_{3}\subseteq {\dv R}^{\indset}$,
${\mathsf C}_{1}$ is a face of ${\mathsf C}_{2}$, which itself is a face of ${\mathsf C}_{3}$, then
${\mathsf C}_{1}$ is a face of ${\mathsf C}_{3}$. Conversely, if ${\mathsf C}_{1}$ is a face of ${\mathsf C}_{3}$
and ${\mathsf C}_{1}\subseteq{\mathsf C}_{2}\subseteq{\mathsf C}_{3}$ then ${\mathsf C}_{1}$ is a face of ${\mathsf C}_{2}$.

\subsection{Normal fan of a polytope}\label{ssec.app.normal.fan}
Recall from \cite[\S\,7.1]{Zie95} that a \underline{finite} collection ${\cal F}$ of \underline{polyhedral}
cones in ${\dv R}^{N}$ is called a {\em fan\/} 
in ${\dv R}^{N}$ if it satisfies the following two conditions:
\begin{itemize}
\item[(i)] if ${\mathsf F}$ is a non-empty face of\/ ${\mathsf C}\in {\cal F}$ then ${\mathsf F}\in {\cal F}$,
\item[(ii)] if ${\mathsf C}_{1},{\mathsf C}_{2}\in {\cal F}$ then  ${\mathsf C}_{1}\cap{\mathsf C}_{2}$ is
a face of both ${\mathsf C}_{1}$ and ${\mathsf C}_{2}$.
\end{itemize}
A fan ${\cal F}$ is {\em complete\/} 
if the union of cones in ${\cal F}$ is ${\dv R}^{N}$. Elementary considerations about the dimension of cones allow one to observe that the inclusion-maximal cones in a complete fan ${\cal F}$
are just the full-dimensional cones in ${\cal F}$. A complete fan ${\cal F}$ is then determined by the collection
${\cal F}^{\max}$ of its full-dimensional cones, characterized by $\bigcup_{{\mathsf C}\in {\cal F}^{\max}} {\mathsf C}={\dv R}^{N}$
and (ii) for ${\cal F}^{\max}$.

Note that ${\mathsf C}_{1},{\mathsf C}_{2}\in {\cal F}^{\max}$, ${\mathsf C}_{1}\neq {\mathsf C}_{2}$, implies\,
$\inte({\mathsf C}_{1})\cap {\mathsf C}_{2}=\emptyset$.
Indeed, by (ii), ${\mathsf C}_{1}\cap{\mathsf C}_{2}$ is a proper face of ${\mathsf C}_{1}$, and, therefore, as explained in Appendix~\ref{ssec.app.vertices-facets}, ${\mathsf C}_{1}\cap{\mathsf C}_{2}$ is disjoint with the interior $\inte({\mathsf C}_{1})$ (recall that ${\mathsf C}_{1}$ is full-dimensional).

Given two complete fans ${\cal F}_{1},{\cal F}_{2}$ in ${\dv R}^{N}$ we say that
${\cal F}_{1}$ {\em refines\/} ${\cal F}_{2}$, or  that ${\cal F}_{2}$ {\em coarsens\/} ${\cal F}_{1}$,
if, for every every ${\mathsf C}_{1}\in {\cal F}_{1}$, a cone ${\mathsf C}_{2}\in {\cal F}_{2}$ exists
with ${\mathsf C}_{1}\subseteq {\mathsf C}_{2}$. This appears to be equivalent to the condition
that every ${\mathsf C}\in {\cal F}^{\max}_{2}$ is the union of cones from ${\cal F}^{\max}_{1}$.

Given a non-empty finite collection of complete fans $\{ {\cal F}_{i} : i\in I\}$ in ${\dv R}^{N}$,
their {\em common refinement\/} is the collection ${\cal F}$ consisting of the cones of the form
${\mathsf C}:=\bigcap_{i\in I}{\mathsf C}_{i}$, where ${\mathsf C}_{i}\in {\cal F}_{i}$.
It appears to be a complete fan in ${\dv R}^{N}$ as well; see \cite[Definition\,7.6]{Zie95}.

Given a non-empty polytope ${\mathsf P}\subseteq {\dv R}^{N}$,
the ({\em outer\/}) {\em normal cone\/} to\/ ${\mathsf P}$ at $y\in {\mathsf P}$ is the set
$$
\nor_{{\mathsf P}}(y) ~:=~ \{\, x\in {\dv R}^{N}\, :\ \forall\, z\in {\mathsf P}
\quad \langle x, y-z\rangle\geq 0\,\}\,,
$$
which is another way of describing the set of $x\in {\dv R}^{N}$ such that the linear function
$z\mapsto \langle x,z\rangle$ achieves its maximum within
${\mathsf P}$ at $y$, that is,
$\nor_{{\mathsf P}}(y) \,=\, \{\, x\in {\dv R}^{N}\,:\ \langle x,y\rangle =
\max_{\,z\in {\mathsf P}} \langle x,z\rangle\,\}$.
One can give an explicit inequality description of the normal cones at vertices of ${\mathsf P}$
in terms of their (geometric) neighbours:
$$
y\in\ext({\mathsf P}) ~~\Rightarrow~~ \nor_{{\mathsf P}}(y) \,=\, \{\, x\in {\dv R}^{N}\,:~~ \forall\, z\in\ne_{{\mathsf P}}(y)\quad
\langle x,y-z\rangle\geq 0\,\},
$$
which fact can be derived from \cite[Lemma\,3.6]{Zie95}.

The (outer) {\em normal fan\/} of a non-empty polytope ${\mathsf P}\subseteq {\dv R}^{N}$ can then be introduced as the collection of all such cones:
$$
\norm_{{\mathsf P}} ~:=~ \{\, \nor_{{\mathsf P}}(y)\, :\
y\in {\mathsf P}\,\}\,.
$$
This collection is, in fact, finite. This observation can be done by means of the method of Galois connections
presented in Appendix~\ref{ssec.app.Galois}. Specifically, we put $X:={\dv R}^{N}$, $Y:={\mathsf P}$
and define an incidence relation $\norfan$ between $x\in X={\dv R}^{N}$ and $y\in Y={\mathsf P}$
as follows:
$$
x\norfan y \quad :=~\quad \langle x,y\rangle =\max_{z\in {\mathsf P}} \,\langle x,z\rangle\,.
$$
This allows one to define formally the normal cone to ${\mathsf P}$ for arbitrary subset ${\mathsf F}\subseteq {\mathsf P}$ by means of the backward Galois connection $\lGalo$ based on this
particular incidence relation $\norfan$\,:
$$
\nor_{{\mathsf P}}({\mathsf F}) ~:=~ {\mathsf F}^{\lGalo} =\{\, x\in {\dv R}^{N}\, :\
\langle x,y\rangle =\max_{z\in {\mathsf P}} \,\langle x,z\rangle
\quad \mbox{for every $y\in {\mathsf F}$}\,\}\,.
$$
If ${\mathsf P}$ is a singleton then the induced lattices are also singletons:
${\calX}^{\norfan}=\{\,{\dv R}^{N}\}$ and ${\calY}^{\norfan}=\{{\mathsf P}\}$. In the non-degenerate
case $|{\mathsf P}|\geq 2$, however, ${\calY}^{\norfan}$ appears to be the whole face lattice of ${\mathsf P}$,
${\calX}^{\norfan}=\{\,  \nor_{{\mathsf P}}({\mathsf F})\, :\ \mbox{${\mathsf F}$ is a face of\/ ${\mathsf P}$}\,\}$,
while $\nor_{{\mathsf P}}({\mathsf F})=\nor_{{\mathsf P}}(y)$ for any face ${\mathsf F}$ of\/ ${\mathsf P}$ and
$y\in\relint ({\mathsf F})$. In particular, $\norm_{{\mathsf P}}$ can equivalently be introduced as a finite set:
$$
\norm_{{\mathsf P}} ~=~ \{\, \nor_{{\mathsf P}}({\mathsf F})\, :\
\mbox{${\mathsf F}$ is a non-empty face of\/ ${\mathsf P}$}\,\}\,.
$$
This allows one to observe that $\norm_{{\mathsf P}}$ is a complete fan of polyhedral
cones in ${\dv R}^{N}$  \cite[Example\,7.3]{Zie95}; thus, the maximal cones in $\norm_{{\mathsf P}}$ are precisely the full-dimensional cones in $\norm_{{\mathsf P}}$.
The method of introducing of\/ $\norm_{{\mathsf P}}$ by means of Galois connections also implies that $(\norm_{{\mathsf P}},\subseteq)$ is anti-isomorphic to the poset of non-empty faces of\/ ${\mathsf P}$.
In fact, the observation that these two posets are graded, with the dimension playing the role
of the height function, allows one to observe that, for any non-empty face ${\mathsf F}$ of ${\mathsf P}$,
one has $\dim({\mathsf F}) +\dim(\nor_{{\mathsf P}}({\mathsf F}))=n=|N|$.
Hence, the maximal cones in $\norm_{{\mathsf P}}$ are just the normal cones
at the vertices of\/ ${\mathsf P}$, that is, the cones $\nor_{{\mathsf P}}(y)$, where $y\in\ext ({\mathsf P})$.
One can also show that $\nor_{{\mathsf P}}(y)=\nor_{{\mathsf P}}(z)$ for $y,z\in {\mathsf P}$ iff
they generate the same face, that is, $\face_{{\mathsf P}}(y)=\face_{{\mathsf P}}(z)$ (see Appendix~\ref{ssec.app.vertices-facets}).

Given a Minkowski sum ${\mathsf P}={\mathsf Q}\oplus{\mathsf R}$ of non-empty polytopes in ${\dv R}^{N}$,
$\ext({\mathsf P})\subseteq\ext({\mathsf Q})\oplus\ext({\mathsf R})$ implies
that the normal fan $\norm_{{\mathsf P}}$ refines both $\norm_{{\mathsf Q}}$
and $\norm_{{\mathsf R}}$. In fact, $\norm_{{\mathsf P}}$ is nothing but the common
refinement of $\norm_{{\mathsf Q}}$ and $\norm_{{\mathsf R}}$ \cite[Proposition\,7.12]{Zie95}.

Every non-zero vector ${\mathsf 0}\neq y\in {\dv R}^{N}$ defines a hyperplane $\{\, x\in {\dv R}^{N}:  \langle x,y\rangle= 0\}$ and two half-spaces $\{\, x\in {\dv R}^{N}: \langle x,y\rangle\leq 0\}$ and
$\{\, x\in {\dv R}^{N}:  \langle x,y\rangle\geq 0\}$ in ${\dv R}^{N}$. These three cones give together an
elementary example of a complete fan in ${\dv R}^{N}$.
Every {\em vector configuration}, that is, a non-empty finite collection $\{\, y_{i}\in {\dv R}^{N}\setminus\{{\mathsf 0}\} : i\in I\}$ of non-zero vectors then induces a complete fan in ${\dv R}^{N}$ as the common refinement of the respective three-cone fans, interpreted as the fan of the {\em hyperplane arrangement\/}
given by the vector configuration $\{\, y_{i} : i\in I\}$.

One of classic examples of a complete fan in ${\dv R}^{N}$ is the so-called braid arrangement fan \cite[\S\,3.2]{PRW08}, induced by the configuration
of vectors $y\in {\dv R}^{N}$ of the form $y:=\chi_{u}-\chi_{v}$ for $u,v\in N$, $u\neq v$. To simplify the terminology, we will call it briefly the {\em braid fan}, following to \cite{CL20}.

\subsubsection*{Acknowledgements}
I am thankful to Federico Castillo and Jos\'{e} Samper for their help when searching the simplest arguments in the proof of Corollary~\ref{cor.braid-cone}. My thanks also belong to L\'{a}szl\'{o} Csirmaz who
provided me with the example presented in Section~\ref{ssec.dead-ends}.

\end{document}